\documentclass{tac}
\usepackage{xy}
\usepackage{MnSymbol}
\usepackage{tikz}

\xyoption{all}

\def\mathcaldef#1{\expandafter\def\csname#1\endcsname{{\cal#1}}}

\def\"{``}
\def\q{\quad}
\def\qq{\quad\quad}
\def\qv{\qq ;\qq}

\def\iso{\,\cong\,}
\def\imp{\Rightarrow}
\def\la{\langle}
\def\ra{\rangle}
\def\adj{\dashv}
\def\op{^{\rm op}}
\def\ov{\overline}

\def\tm{\times}
\def\otm{\otimes}
\def\al{\alpha}

\def\si{\sigma}
\def\t{_\blacktriangleright}
\def\fp{_{fp}}
\def\cd{\cdots}

\mathcaldef{C}
\mathcaldef{D}
\mathcaldef{M}
\mathcaldef{N}
\mathcaldef{L}
\mathcaldef{S}
\mathcaldef{P}
\mathcaldef{V}

\mathbfdef{Set}
\mathbfdef{Pos}
\mathbfdef{cMon}
\mathbfdef{Mon}
\mathbfdef{Mod}
\mathbfdef{Mlt}
\mathbfdef{fpMlt}
\mathbfdef{sMlt}
\mathbfdef{Rep}
\mathbfdef{fpRep}
\mathbfdef{sRep}
\mathbfdef{Seq}
\mathbfdef{fpSeq}
\mathbfdef{Cat}
\mathbfdef{fpCat}
\mathbfdef{Sum}
\mathbfdef{fpSum}
\mathbfdef{mltfb}
\mathbfdef{mltdf}
\mathbfdef{Prof}
\mathbfdef{MProf}
\mathbfdef{preAdd}
\mathbfdef{Add}

\mathrmdef{id}
\mathrmdef{I}
\mathrmdef{m}
\mathrmdef{v}
\mathrmdef{s}
\mathrmdef{d}

\def\obj{{\rm obj\,}}

\newtheorem{prop}{Proposition}
\newtheorem{corol}{Corollary}
\let\pf\proof
\let\epf\endproof
\def\eq{\begin{equation}}
\def\eeq{\end{equation}}

\author{Claudio Pisani}
\address{via Saluzzo 67,\\ 10125 Torino, Italy.}
\title{Sequential multicategories}
\copyrightyear{2014}
\eaddress{pisclau@yahoo.it}

\begin{document}

\maketitle

\begin{abstract}

We study the monoidal closed category of symmetric multicategories,
especially in relation with its cartesian structure and with sequential multicategories
(whose arrows are sequences of concurrent arrows in a given category).
%
Then we consider cartesian multicategories in a similar perspective and develop 
some peculiar items such as algebraic products. 
Several classical facts arise as a consequence of this analysis when some of
the multicategories involved are representable.

\end{abstract}

\section{Introduction}
\label{intro}

The overall aim of the present work is to show how symmetric and cartesian multicategories 
offer a natural frame encompassing several aspects of the theory of symmetric monoidal
and finite product categories. 

It is well-known for instance that \"preadditive" categories, 
whose hom-sets are commutative monoids and composition is distributive, 
occupy a special place among enriched categories.
Here they are characterized as the categories of commutative monoids 
in a cartesian multicategory; 
consequently, the (semi)additive categories are characterized as 
the categories of commutative monoids in a finite product category.

A crucial role is played by the \"Boardman-Vogt" monoidal closed structure
on the category $\sMlt$ of symmetric multicategories (see \cite{weiss} and \cite{trova}).
The internal hom $[\M,\N]$ has as objects the functors $F:\M\to\N$ in $\sMlt$ and as arrows
$\al:F_1,\cd,F_n\to F$ the families of arrows $\al_A:F_1A,\cd,F_nA\to FA$ in $\N$ ($A\in\M$), 
such that the following \"naturality" condition holds for any arrow $f:A_1,\cd,A_m\to A$ in $\M$:
\[
Ff(\al_{A_1},\cd,\al_{A_m}) = \si \al_A(F_1f,\cd,F_nf)
\]
where $\si$ is the obvious permutation.
The composition and the symmetric structure are inherited pointwise
from $\N$ (see also \cite{tronin}).
If $\N$ has a symmetric monoidal structure (that is, if it is representable) or a cartesian
structure, these are also inherited pointwise by $[\M,\N]$. 

In particular, if $1\t\in\sMlt$ is the terminal multicategory,
then $[1\t,\M]\in\sMlt$ is the multicategory of commutative monoids in $\M$.
We call \"sequential" these multicategories $[1\t,\M]$ of commutative monoids, 
since they can be also characterized as those of the form $\C\t$, where 
$(-)\t:\Cat\to\sMlt $ is the fully faithful \"discrete cocone" functor: 
 \[ 
\C\t(X_1,\cdots , X_n ; X) \iso \C(X_1;X)\tm \cdots \tm \C(X_n;X)
\]
So $n$-ary arrows in $\C\t$ are sequences $\la f_1,\cd,f_n \ra$ of concurrent arrows in $\C$
with the obvious composition, for instance 
\[
\la f,g,h\ra(\la l,t\ra,\la\,\ra,\la q \ra) = \la f l,f t,h q \ra 
\]
and for any functor $F:\C\to\D$
\[
F\t \la f_1,\cd, f_n\ra = \la F f_1,\cdots, F f_n \ra
\]

We denote by $\Seq\subset\sMlt$ the full subcategory of sequential
multicategories, so that there is an equivalence $\Seq \simeq \Cat$.
In fact, $[1\t,\M]$ is the coreflection of $\M$ in $\Seq$, while $1\t\otm\M$ is
its reflection.

The relation between (symmetric) monoidal categories and (symmetric) multicategories is well-known 
(see for instance \cite{lambek}, \cite{hermida} and \cite{leinster});  
if $(\C,\I,\otm)$ is a monoidal category, one gets a multicategory $\C_\otm$ by
\[
\C_\otm(X_1,\cdots , X_n ; X) \iso \C(X_1\otm \cdots \otm X_n;X)
\]
(where the particular choice of bracketing is omitted). 
In fact, there is an equivalence between the category of (symmetric) monoidal categories 
and lax monoidal functors and the category $\Rep \subset \Mlt$ ($\sRep \subset \sMlt$)
of representable (symmetric) multicategories.
Working directly in $\Rep$ has the advantage that coherence issues are replaced by 
the more natural universal issues.

Since $\C\t$ is representable iff $\C$ has finite coproducts,
$\Sum := \Seq \cap \sRep \subset \sMlt$ is equivalent to the category
of categories with finite coproducts (and all functors).
Thus $[1\t,-]$ gives also the coreflection of $\sRep$ in $\Sum$,
and finite sum categories can be characterized as the categories
of commutative monoids in a symmetric monoidal category. 

If $\M\in\fpMlt$ is a cartesian multicategory (see for instance \cite{gould} and the references therein), 
there are in particular \"contraction" mappings
\[
\gamma_n:\M(A,\cd,A;B)\to \M(A;B)
\]
We will see that giving a cartesian structure on a sequential multicategory $\C\t$ 
is equivalent to giving a preadditive structure on $\C$
(where contractions correspond to sums of maps).
Thus the full subcategory $\fpSeq \subset \fpMlt$ of sequential cartesian multicategories
is equivalent to the category of preadditive categories and additive functors,
and $[1\t,-]$ gives also the coreflection of $\fpMlt$ in $\fpSeq$.
In particular, we get the characterization of preadditive categories mentioned 
at the beginning of this introduction;
indeed, the preadditive category corresponding to $[1\t,\M]$ is the category $\cMon(\M)$ 
of commutative monoids in $\M$, with the sum of the monoid morphisms 
$f_1,\cd,f_n:M\to N$ given by the contraction
\[
\gamma_n m_n(f_1,\cd,f_n)
\]
(the arrows $m_n:N,\cd,N\to N$ being the monoid structure).

Furthermore, the well-known interplay in additive categories between finite sums, 
products and (algebraic) biproducts arises as a consequence of the interplay between 
representability and universal and algebraic products in a cartesian multicategory.


A useful perspective is to consider symmetric and cartesian multicategories
as two doctrines (in the sense of Lawvere's categorical logic):
those of \"linear" theories and those of algebraic theories respectively. 
In fact finite product categories, which are the usual many-sorted version of Lawvere theories
(see for instance \cite{adamek}), are included in $\fpMlt$ as the representable cartesian multicategories.

In this perspective, the category of models of $\M$ in $\N$,
$\sMlt(\M,\N)$ or $\fpMlt(\M,\N)$, 
has itself the structure of a symmetric or cartesian multicategory:
\[
[\M,\N]  \qv  [\M,\N]\fp 
\]  
(where $[\M,\N]_{fp}$ is the internal hom for a monoidal closed structure on $\fpMlt$).
Therefore it makes sense to consider models valued in $[\M,\N]$ or $[\M,\N]\fp$
and the Boardman-Vogt tensor product (on $\sMlt$ or $\fpMlt$) and sequential multicategories 
play the same role that the Kronecker product 
and annular theories play in Lawvere theories (see \cite{freyd}, \cite{lawvere} and \cite{wraith}).

On the symmetric level, categories $\C\in\Cat$ parameterize
two important sorts of linear theories: the unary ones $\C^-$ 
(which have only unary arrows) and the sequential ones $\C\t$:
\[
(-)^- :\Cat \to \sMlt  \qv  (-)\t:\Cat \to \sMlt
\]
For instance (considering the basic background $\Set$), the models for $1\t$ 
are commutative monoids, the models for $\C^-$ are presheaves on $\C$,
and the isomorphism $\C\t \iso \C^-\otm 1\t$ says 
that the models for $\C\t$ are monoids in $\C$-presheaves (or $\C$-presheaves of monoids).

On the cartesian level, preadditive categories parameterize \"annular" theories
\[
\qq\qq (-)\t:\cMon{\rm-}\Cat \to \fpMlt
\]
and a consequence of the coreflection $[1\t,-]:\fpMlt\to\fpSeq$ is that
the models in any $\M$ for an annular theory, namely the functors $\C\t\to\M$ in $\fpMlt$,
are indeed (generalized) modules $\C\t\to [1\t,\M]$, that is additive functors $\C \to \cMon(\M)$.


Summarizing, the category $\sMlt$ of symmetric multicategories includes
the category $\Cat$ of categories in two ways
(as sequential multicategories and as unary multicategories)
as well as the category of symmetric monoidal categories (as represanteble multicategories).
Similarly, the category $\fpMlt$ of cartesian multicategories includes 
preadditive categories (as sequential multicategories)
and finite product categories (as represanteble multicategories).
The two levels are themselves related by an adjunction and
in both cases the closed structure (in particular the monoid construction)
plays a prominent role;
it restricts, on sequential multicategories, to the cartesian closed structure
of $\Cat$ and to the monoidal closed structure of $\cMon$-$\Cat$ respectively.
This unifying role of multicategories, along with the perspective of categorical logic and a feasible
diagrammatic calculus, provides an effective point of view on some aspects of
category theory.



Some of the items studied here have been considered also in \cite{pisani}.


\subsection{Summary}

In Section \ref{cart} we shortly review some basic aspects of multicategories
and investigate the cartesian structure of $\Mlt$.
In particular, we observe that the exponentiable multicategories
coincide with the promonoidal ones.
In this perspective, Day convolution appears as the monoidal structure
on the exponential $\N^\M$ in $\Mlt$ when $\M$ is promonoidal and $\N$ is monoidal 
(that is representable) and cocomplete.

If $\M$ is sequential, we have a particularly simple formula for the exponential $\N^\M$:

\[
\N^\M(F_1,\cd,F_n; F) = \int_A \N(F_1A,\cd,F_nA;FA)
\]
that is an arrow $\al:F_1,\cd,F_n\to F$ in $\N^\M$ consists of arrows
$\al_A:F_1A,\cd,F_nA \to FA$ ($A\in\M$) such that for any unary arrow $f:A\to B$ in $\M$
\[
\al_B(F_1f,\cd,F_nf) = (Ff)\al_A
\]

In Section \ref{mon} we study the monoidal closed structure of $\sMlt$ 
and compare it with the cartesian one.
The main points are:

\begin{itemize}
\item
$[\C^-,\N] \iso \N^{\C\t}$ and $\C^-\otm\N \iso \C\t\tm\N$,
so that the cartesian and the monoidal structures in fact coincide on $\Seq\subset\sMlt$.
\item
Sequential multicategories are characterized as those with a \"central monoid".
\item
The functors $1\t\otm - $ and $[1\t,-]$ 
give respectively a reflection and a coreflection of $\sMlt$ in $\Seq$.
\item
Sequential multicategories form an ideal with respect to $\otm$ and $[-,-]$
(on both sides).
\item
If $\N$ is representable then so it is also $[\M,\N]$, for any $\M\in\sMlt$.
\end{itemize}

As a consequence, we derive under a unified perspective some important facts
(which seem to be well-known at least as folklore): 
\begin{itemize}
\item
The characterization of cocartesian symmetric monoidal categories as those 
with a \"central monoid".
\item
If $\M$ is symmetric monoidal, the category $\cMon(\M)$ of commutative
monoids in $\M$ is cocartesian with the tensor product inherited by $\M$,
and the functor $\cMon(-)$ gives a right adjoint to the inclusion of 
cocartesian categories in the monoidal ones.
\end{itemize}

In Section \ref{fp} we move to the level of cartesian multicategories.

If $\C$ is a finite product category, then $\C_\tm$ has a cartesian structure.
In fact, the full subcategory $\fpRep\subset\fpMlt$ of representable cartesian multicategories 
is equivalent to the category of finite product categories and finite product preserving functors.

If $\C$ is preadditive, that is enriched in commutative monoids, 
then $\C\t$ has a cartesian structure. 
In fact, the full subcategory $\fpSeq\subset\fpMlt$ of sequential cartesian multicategories  
is equivalent to the category of preadditive categories and additive functors.

The view of cartesian multicategories as a common generalization of
finite product categories and of preadditive categories
is not devoid of consequences:
\begin{itemize}
\item
\(
\fpRep \cap  \fpSeq  \subset \fpMlt
\)
is equivalent to the category of additive categories
(that is preadditive categories with (bi)products) and additive functors.
\item
If $\M\in\fpSeq$ and $\N\in\fpRep$, the morphisms $\M\to\N$ in $\fpMlt$ can be seen
as generalized modules. In particular, if $\M$ is an operad
(that is, has just one object) the morphisms $\M\to\Set_\tm$ 
coincide with the modules over the rig $\M$.
\item
A cartesian multicategory is representable iff it has \"algebraic products"
(which in the sequential case reduce to ordinary algebraic biproducts)
In fact, the well-known relations between $\cMon$-enrichments, products, coproducts and biproducts are
particular cases of more general properties of cartesian multicategories.
\end{itemize}




\section{The cartesian structure of the category of multicategories}
\label{cart}

For an introduction to plain and symmetric multicategories
(also known as coloured operads) the reader may consult
for instance \cite{leinster} and \cite{trova}.
Note that a symmetric structure on a multicategory can be defined in the same way 
as a cartesian structure (see Section \ref{fp}), except that only bijective
mappings are supposed to act on the hom-sets. 
We do not consider here generalized or enriched multicategories.

We denote by $\Mlt$ the category of multicategories and functors,
and by $\sMlt$ the category of symmetric multicategories and (symmetric) functors.
Recall that given $F,G:\M\to\N$, a natural transformation $\al:F\to G$ consists of
a family of unary arrows $\al_A:FA\to GA$ ($A\in\M$) such that for any arrow
$f:A_1,\cd,A_n \to A$ in $\M$
\[
Gf(\al_{A_1},\cd,\al_{A_n}) = \al_A Ff
\]
Thus $\Mlt$ and $\sMlt$  are in fact 2-categories, with natural transformations as 2-cells.


\subsection{The unary embedding and the underlying functor}

A category $\C$ gives rise to a (symmetric) \"unary" multicategory $\C^-$ consisting only of unary arrows:
\[
\C^-(X;Y) = \C(X;Y) \qv \C^-(X_1,\cdots , X_n ; Y) = \emptyset \q (n\neq 1)
\]
and the construction clearly extends to full and faithful 2-functors:
\[ 
(-)^- :\Cat\to\Mlt  \qv  (-)^- : \Cat\to\sMlt
\]
In the other direction, there are \"underlying" 2-functors which take any multicategory $\M$   
to the category $\M_-$ with $\M_-(X;Y) = \M(X;Y)$.
It is immediate to verify that there are adjunctions
\[ 
(-)^- \adj (-)_-:\Mlt\to\Cat  \qv  (-)^- \adj (-)_-:\sMlt\to\Cat
\]
\begin{remark}
\label{fr}
The adjunction $(-)^- \adj (-)_-$ satisfies the Frobenius law:
\[ 
\M\tm\C^- \iso (\M_-\tm\C)^-
\]
In particular, $\M\tm 1^- \iso (\M_-)^-$.
\end{remark}


\subsection{The sequential embedding}

The discrete cocone functor $(-)\t:\Cat\to\Mlt$ mentioned in the introduction
(see also \cite{hermida}) is in fact a full and faithful 2-functor.
For fullness, let $F:\C\t \to \D\t$ be a functor in $\Mlt$,
$\la f,g\ra \in \C\t(X,Y;Z)$, and $F \la f,g\ra = \la f',g'\ra \in \D\t(FX,FY;FZ)$;
then
\[
f' = F \la f,g\ra (\id_{FX},\la\,\ra) = F \la f,g\ra (F\id_X,F\la\,\ra) = 
F (\la f,g\ra (\id_X,\la\,\ra)) = F f
\]
so that in fact $F = (F_-)\t$, where $F_- : \C \to \D$ is the \"underlying" functor.
We say that $\M$ is \"sequential" if it is isomorphic to some $\C\t$.

Sequential multicategories have obvious natural symmetric structures
which are preserved by any functor. Thus we also have an embedding $(-)\t:\Cat\to\sMlt$ which restricts
to an equivalence $\Cat\simeq\Seq$, where $\Seq\subseteq\sMlt$ is the full subcategory of sequential
multicategories.

\begin{remark}
Note that $1\t$ is terminal in $\Mlt$ and in $\sMlt$. 
In fact, we will see in Corollary \ref{lradj} that $(-)\t:\Cat\to\sMlt$ 
has both a left and a right adjoint, so that it preserves (co)limits.
Recall that the category $\Mlt(1\t,\N)$  and $\sMlt(1\t,\N)$ can be identified with the 
categories $\Mon(\N)$ and $\cMon(\N)$ of monoids  and commutative monoids in $\N$, respectively.
In particular, when $\N$ is representable we find again the usual notion of (commutative)
monoid in a (symmetric) monoidal category.

Note also that any full sub-multicategory $\N\subseteq\M$ of a sequential multicategory 
is itself sequential (on the corresponding full subcategory of $\M_-$).
\end{remark}


\subsection{The monoidal embedding}

For a detailed treatment of the relations between multicategories and monoidal
categories we refer to \cite{hermida} and \cite{leinster}.
Let us just recall the main points.
An arrow $A_1,\cd, A_n \to A$ in a multicategory $\M$ is \"preuniversal" if it gives a
representation for the functor $\M(A_1,\cd, A_n,-):\M_-\to\Set$;
$\M$ is \"representable" if for any $A_1,\cd, A_n \in \M$ there is a preuniversal arrow
$A_1,\cd, A_n \to A$ and if preuniversal arrows are closed with respect to composition.

Equivalently, $\M$ is representable if it has a representation in the following sense:
to any $A_1,\cd, A_n \in \M$ it is assigned a \"universal" arrow $u_{A_1,\cd, A_n}:A_1,\cd, A_n \to A$
such that, for any double sequence $A_{i1},\cd, A_{im_i}$,
composition with $u_{A_{i1},\cd, A_{im_i}}:A_{i1},\cd, A_{im_i} \to A_i$ yields (for any $B\in\M$)
a bijection $\M(A_1,\cd, A_n, B) \to \M(A_{11},\cd, A_{nm_n}, B)$.

We denote by $\Rep\subset\Mlt$ and $\sRep\subset\sMlt$ the full subcategories
of (symmetric) representable multicategories.
When we will consider representable multicategories, we assume that a representation
in the above sense is given.

\begin{prop}
$\Rep$ ($\sRep$) is equivalent to the category of (symmetric) monoidal categories 
and lax monoidal functors.
\end{prop}
\pf
We just give an idea of the correspondence.
As mentioned in the introduction, to any monoidal category $\C$ there corresponds 
a multicategory $\C_\otm$. 
Conversely, a representation for $\M$ yields a tensor product for $\M_-$.

Since the arrows in a representable multicategory are generated by unary arrows and universal arrows,
the conditions for $F:\M\to \N$ to be a functor in $\Mlt$ correspond to
the conditions for $F_-:(\M_-,\I,\otm)\to(\N_-,\I,\otm)$ to be a monoidal functor:
\begin{enumerate}
\item
preservation of composition of unary arrows corresponds the functoriality of $F_-$;
\item
the assignment of an immage in $\N$ to the universal arrows $A,B\to A\otm B$ in $\M$
corresponds to the assignment of the arrows $FA\otm FB \to F(A\otm B)$;
\item
preservation of composition of unary with universal arrows corresponds to
the naturality of $FA\otm FB \to F(A\otm B)$;
\item
preservation of composition of universal arrows
corresponds essentially to the coherence conditions for a monoidal functor.
\end{enumerate}
We so get in fact a 2-equivalence: the condition for $\al_A:FA\to GA$ to be natural becomes, when $\M,\N\in\Rep$,
the commutativity condition with respect to unary arrows (ordinary naturality) and with respect to
universal arrows, giving the usual definition of monoidal natural transformation.
\epf


\subsection{Multicategories over a multicategory}

A good deal of the analysis of categories over a category can be extended
to multicategories.
In particular, recall the following well-known facts (see for instance \cite{benabou}):
\begin{enumerate}
\item
Categories over $\C$ correspond to normalized lax functor $\C\to\Prof$.
\item
Exponentiable multicategories over $\C$ correspond to pseudofunctors $\C\to\Prof$.
\item
Fibrations over $\C$ correspond to pseudofunctors $\C\to\Prof$
which furthermore factor through the inclusion $\Cat\to\Prof$.
\end{enumerate}
(Of course, the sense of these \"correspondences" is made precise by suitable equivalences.)

In order to generalize to multicategories, we need to replace bicategories
with \"bi-multicategories": they have objects, hom categories $\M(X_1,\cd,X_n;X)$ and
a composition which is unitary and associative up to a coherent isomorphism;
there are obvious notions of lax and pseudo functors of bi-multicategories.
The main instance is $\MProf$, whose objects are categories and arrows are \"multiprofunctors"
$f:X_1\op\tm\cd\tm X_n\op\tm X \to \Set$ with the obvious coend composition.
(Note that $\MProf$ is in fact a monoidal bicategory, since profunctors
$X_1,\cd, X_n\to X$ correspond to profunctors $X_1\tm\cd\tm X_n\to X$.)

Then one can prove:
\begin{enumerate}
\item
Multicategories over $\M$ correspond to (lax) functors $\M\to\MProf$.
\item
Exponentiable multicategories over $\M$ correspond to pseudofunctors $\M\to\MProf$.
\item
Fibrations of multicategories over $\M$ (see \cite{hermida2}) correspond to pseudofunctors $\M\to\MProf$
which furthermore factor through the inclusion $\Cat_\tm\to\MProf$.
\end{enumerate}

When the base multicategory is the terminal one, we get the corresponding absolute notions:
the fibrations over $1\t$ are the representable multicategories
while the exponentiable multicategories are
the pseudomonoids in $\MProf$, that is the promonoidal (multi)categories
(see for instance \cite{day}).
We will prove directly in Proposition \ref{pro} 
that exponentiable and promonoidal multicategories indeed coincide.


\subsection{Powers of multicategories}

If the power $\N^\M$ does exist in $\Mlt$, then its objects are the functors $\M_-\to\N_-$, while
arrows $F_1,\cd,F_n\to F$ are mappings 
\[
\M(A_1,\cd,A_n;A) \to \N(F_1A_1,\cd,F_nA_n;FA)
\]
that are natural in all the variables:
\eq 
\label{exp}
\N^\M(F_1,\cd,F_n; F) = \int_{A_1,\cd,A_n,A}\Set(\M(A_1,\cd,A_n;A) , \N(F_1A_1,\cd,F_nA_n;FA))
\eeq
Indeed, for objects we have
\[
\Mlt(1^-,\N^\M) \iso \Mlt(1^-\tm\M,\N) \iso \Mlt((\M_-)^-,\N) \iso \Mlt(\M_-,\N_-)
\]
while the $n$-ary arrows $F_1,\cd,F_n\to F$ are given by those functors $t^n\tm\M \to \N$
which restrict to the given functors on the objects of $t^n$ (the generic $n$-ary arrow).
Thus (\ref{exp}) follows from the description of $t^n\tm\M$ as given (for $n=2$)
in the proof of Proposition \ref{pro}.

In particular, for $n=1$ we get $(\N^\M)_- = \N_-^{\M_-}$, because
\[
\N^\M(G; F) = \int_{B,A}\Set(\M(B;A) , \N(GB;FA)) = \int_{A\in\M}\N(GA;FA)
\]
or also because $t^1\tm \M \iso t^-\tm \M \iso  (t\tm \M_-)^- $, where $t$ is the arrow category,
so that $\Mlt(t^1\tm\M , \N) \iso \Mlt((t\tm \M_-)^- ,\N) \iso \Cat(t\tm \M_-,\N_-)$.

So the undelying functor $(-)_-:\Mlt\to\Cat$ preserves exponentials, when they exist.
In fact this is a well-known direct consequence of the Frobenius law for the adjunction $(-)^- \adj (-)_-$
(see Remark \ref{fr}).


\begin{prop}
\label{pro}
The exponentiable multicategories are the promonoidal ones.
In particular, representable multicategories and sequential multicategories are exponentiable,
while no unary non-empty multicategory $\C^-$ is exponentiable.
\end{prop}
\pf
Suppose first that $\M$ is exponentiable and consider the pushout $q$ of the generic 
2-ary arrow $t:X,Y\to Z$ with itself along $X:1^- \to t$ and $Z:1^- \to t$. 
Thus $q$ has three non-identity arrows: $t$, $t'$ and their composite $t''=t(t',Y):U,V,Y\to Z$
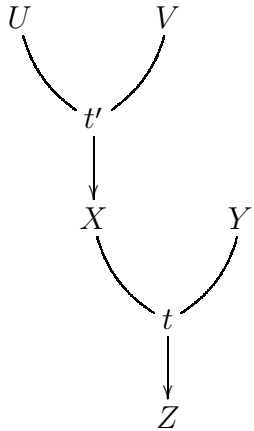
\[
\xymatrix@R=2pc@C=1pc{
U  \ar@/_/@{-}[dr]                &         &   V  \ar@/^/@{-}[dl]   &           \\
                 & t' \ar[d]              &                                \\
                 & X \ar@{-}@/_/[dr]    &                      & Y \ar@{-}@/^/[dl]                            \\
                 &                               &  t \ar[d]                              &      &                    \\
                 &        &  Z                                          &      &
}
\]
The product $\M\tm t$ consists of three copies of the unary $(\M_-)^-$, say $\M_X$, $\M_Y$ and $\M_Z$,
an arrow $f_t:A_X,B_Y\to C_Z$ for any 2-ary arrow $f:A,B\to C$ in $\M$ and the composition rule
$c_Zf_t(a_X,b_Y) = [cf(a,b)]_t$ for $a:A'\to A$, $b:B'\to B$ and $c:C\to C'$.

Similarly, the product $\M\tm q$ consists of five copies of the unary $(\M_-)^-$, say $\M_X$, $\M_Y$, $\M_Z$,
$\M_U$ and $\M_V$, arrows $f_t:A_X,B_Y\to C_Z$, $f_{t'}:A_U,B_V\to C_X$ for any 2-ary arrow $f:A,B\to C$ in $\M$  
and $h_{t''}:A_U,B_V,C_Y\to D_Z$ for any 3-ary arrow $h:A,B,C\to D$ in $\M$.
The composition rule are as above and moreover $f_t(g_{t'},B_Y) = [f(g,B)]_{t''}$,
whenever the compositions are meaningful.

On the other hand, it is easy to see that the pushout $\M\tm t +_{\M\tm 1^-}\M\tm t$ 
(along the inclusion of $\M\tm 1^- = (\M_-)^-$ as $\M_X$ and $\M_Z$)
consists of five copies of the unary $(\M_-)^-$, say $\M_X$, $\M_Y$, $\M_Z$,
$\M_U$ and $\M_V$, arrows $f_t:A_X,B_Y\to C_Z$ and $f_{t'}:A_U,B_V\to C_X$ for any 2-ary arrow $f:A,B\to C$ in $\M$
and formal composites $f_t\circ g_{t'}$ (for $f$ and $g$ composable) with the associativity constraints
$f(a,B)_t\circ g_{t'} = f_t\circ (ag)_{t'}$ for any suitable unary arrow $a$.

Since $\M\tm-$ preserves colimits, the canonical $\M\tm t +_{\M\tm 1^-}\M\tm t \to \M\tm q$ 
(sending $f_t\circ g_{t'}$ in $f_t(g_{t'},B_Y)$ and obvious elsewhere) is an isomorphism.
This amounts to say that any $h_{t''}$ actually has the form $f_t(g_{t'},B_Y)$ for $f$ and $g$ unique
up to the above associativity constraints.
Repeating the argumentation for the pushout of $t$ with the generic 0-ary arrow, we get the other
condition for a multicategory to be promonoidal (see \cite{day}).

In the other direction, suppose that $\M$ is promonoidal. 
Let the multigraph $\N^\M$ have the functors $\M_-\to\N_-$ as objects 
and arrows in $\N^\M(F_1,\cd,F_n; F)$ given by (\ref{exp}).
Then we can define the composition of, say, $\al:F,G\to H$ and $\beta:H,L\to M$ as 
the mapping which takes a 3-ary arrow $f:A,B,C\to D$ in $\M$
to the arrow $\beta_h(\al_g,LC):FA,GB,LC\to MD$, where $h(g,C) = f$
is a decomposition of $f$ given by the promonoidal structure of $\M$.
Then it is straightforward to see that 
\begin{enumerate}
\item
the mapping does not depend on the particular decomposition
(since two decompositions related by a unary arrow give the same result);
\item
the mapping is natural;
\item
the composition so defined on the multigraph $\N^\M$ is associative and unitary;
\item
the multicategory so obtained is indeed the exponential in $\Mlt$.
\end{enumerate}
\epf

Next we consider some relevant particular cases of exponentials.


\subsection{Monoidal exponents}

If $\M$ is representable then for any $\N$ the power $\N^\M$ is given by
\[
\N^\M(F_1,\cd,F_n; F) = \int_{A_1,\cd,A_n,A}\Set(\M(A_1,\cd,A_n;A) , \N(F_1A_1,\cd,F_nA_n;FA)) = 
\]
\[
\int_{A_1,\cd, A_n,A}\Set(\M(A_1\otm\cd\otm A_n;A) , \N(F_1A_1,\cd,F_nA_n;FA)) = 
\]
\[
\int_{A_1,\cd, A_n}\N(F_1A_1,\cd,F_nA_n;F(A_1\otm\cd\otm A_n)) 
\]
which for $\N$ also representable becomes
\[
\N^\M(F_1,\cd,F_n; F) =\int_{A_1,\cd, A_n}\N(F_1A_1\otm\cd\otm F_nA_n;F(A_1\otm\cd\otm A_n))  
\]
Indeed, monoids $1\t \to \N^\M$ on the object (functor) $F:\M_-\to\N_-$ correspond to lax monoidal 
structures for $F$.

\begin{remark}
Monoids are included in $\Cat$ as one object categories $\C$; 
thus they are also included in two way in $\Mlt$: as $\C^-$ and as $\C\t$.
In fact, monoids are included in $\Mlt$ in a third way, namely as representable
discrete multicategories (that is with a discrete underlying category).
If $\M$ and $\N$ are monoids in the third sense, $\N^\M\in\Mlt$ has all mappings of the underlying
sets as objects and a (unique) arrow $F,G\to H$ when $(Fx)\cdot (Gy) = H(x\cdot y)$
(and $H$ is the codomain of a 0-ary arrow iff $H1 = 1$).
\end{remark}


\subsection{Monoidal base}

Now we show that Day convolution is nothing but the monoidal structure on $\N^\M$,
when $\N$ is representable and cocomplete.
\begin{prop}
If $\N$ is representable and cocomplete, then $\N^\M$ is also representable (for any promonoidal  $\M$)
by the usual convolution tensor product
\end{prop}
\pf
\[
\N^\M(F_1,\cd,F_n; F) = \int_{A_1,\cd,A_n,A}\Set(\M(A_1,\cd,A_n;A) , \N(F_1A_1,\cd,F_nA_n;FA)) = 
\]
\[
\int_{A_1,\cd,A_n,A}\Set(\M(A_1,\cd,A_n;A) , \N(F_1A_1\otm\cd\otm F_nA_n;FA)) =
\]
\[
\int_{A_1,\cd,A_n,A}\N(\M(A_1,\cd,A_n;A)\cdot F_1A_1\otm\cd\otm F_nA_n;FA)) =
\]
\[
\int_A\N(\int^{A_1,\cd,A_n}\M(A_1,\cd,A_n;A)\cdot F_1A_1\otm\cd\otm F_nA_n;FA)) =
\]
Thus the functor $\int^{A_1,\cd,A_n}\M(A_1,\cd,A_n;-)\cdot F_1A_1\otm\cd\otm F_nA_n$ 
is a representing object for $\N^\M(F_1,\cd,F_n; -)$:
\epf



\subsection{Sequential exponents}
\label{seqexp}

When $\M$ is sequential, the formula (\ref{exp}) becomes 
\[
\int_{A_1,\cd,A_n,A}\Set(\M(A_1;A)\tm\cd\tm\M(A_n;A) , \N(F_1A_1,\cd,F_nA_n;FA))
\]
which by Yoneda reduction gives
\eq
\label{seqexp1}
\N^\M(F_1,\cd,F_n; F) = \int_A \N(F_1A,\cd,F_nA;FA)
\eeq
Thus, an arrow $\al:F_1,\cd,F_n\to F$ in $\N^\M$ amounts to a family of arrows 
\[
\al_A:F_1A,\cd,F_nA\to FA \q
(A\in\M)
\]
such that for any unary $f:A\to B$ in $\M$
\eq
\label{seqexp2}
(Ff)\al_A = \al_B(F_1f,\cd,F_nf)
\eeq
Indeed, in this case (as it happens for ordinary categories) the natural
mappings $\al:\M(A_1,\cd,A_n;A) \to \N(F_1A_1,\cd,F_nA_n;FA)$ are in fact determined by
the image of those of form $\la\id_X,\cd,\id_X\ra :X,\cd,X\to X$, since
for $f_i:A_i\to A$ we have 
\[
\la f_1,\cd,f_n \ra = \la \id_A,\cd,\id_A \ra(f_1,\cd,f_n)
\]
and these have to satisfy condition (\ref{seqexp2}) since for any $f:A\to B$
\[
f\la\id_A,\cd,\id_A\ra = \la\id_B,\cd,\id_B \ra(f,\cd,f)  
\]

If furthermore $\N$ is representable, one thus gets the following well-known particular case 
of the convolution product:
\begin{prop}
If $\N$ is representable then $\N^{\C\t}$ is also representable, in a pointwise way:
\[
\N^{\C\t}(F_1,\cd,F_n; F) = \int_A \N(F_1A\otm\cd\otm F_nA;FA)
\]
\epf
\end{prop}

Since $\Seq$ is equivalent to $\Cat$, it is cartesian closed. 
Moreover, product and exponentials in $\Seq$ can be computed as in $\Mlt$:

\begin{corol}
The discrete cocone functor $(-)\t:\Cat\to\Mlt$ preserves products and exponentials.
Thus $\Seq$ is closed with respect to products and exponentials in $\Mlt$.
\end{corol}
\pf
Since $(-)\t$ clearly preserves products, we need to show that $\D\t^{\,\C\t} \iso (\D^{\,\C})\t$.
By the above description, an arrow $F_1,\cd,F_n \to F$
in $\D\t^{\,\C\t}$ is a sequence of arrows $\al_A^i:F_i A\to FA$, for any $A\in\C$,
such that 
\[
(Ff)\al_A^i = \al_B^i(F_if)
\]
for any $f:A\to B$ in $\C$.
But this amounts exactly to a sequence of natural transformations $\al^i:F_i\to F$, 
that is to an arrow in $(\D^{\,\C})\t$.
\epf

\subsection{Powers of symmetric multicategories}

Powers in $\sMlt$ are computed essentially as in $\Mlt$.
Let $\M,\N\in\sMlt$ and suppose that the power $\N^\M$ of the underlying multicategories
exists in $\Mlt$.
By (\ref{exp}), an arrow $\al:F_1,\cd,F_n\to F$ is given by natural mappings 
$\al:\M(A_1,\cd,A_n;A) \to \N(F_1A_1,\cd,F_nA_n;FA)$.
If we define 
\eq
\label{symexp}
\si \al : f \mapsto \si\al(\si^{-1}f)
\eeq
it is easy to see that $\si\al$ is also natural and that we so get a symmetric structure on $\N^\M$.

\begin{prop}
In the above hypothesis, $\N^\M$ with the symmetric structure given by {\rm(\ref{symexp})}
is also the power in $\sMlt$.
\end{prop}
\pf
A functor $\al:\L\to\N^\M$ in $\sMlt$ is a functor $\al:\L\to\N^\M$ in $\Mlt$ such that,
for any $g$ in $\L$ and any permutation $\si$ of its domain, $\al(\si g) = \si(\al g)$, that is
\[
\al(\si g): f \mapsto \si\al(\si^{-1}f)   \iff   \al(\si g): \si f \mapsto \si(\al f)
\]
When transposed, $\al$ corresponds thus to a functor $\ov\al:\L\tm\M\to\N$ such that
\[
\ov\al(\si g,\si f) = \si\ov\al(g, f)
\]
that is to a functor $\ov\al:\L\tm\M\to\N$ in $\sMlt$.
\epf

Thus, promonoidal symmetric multicategories are exponentiable in $\sMlt$.
We have seen in Section \ref{seqexp} that the 
$\al:\M(A_1,\cd,A_n;A) \to \N(F_1A_1,\cd,F_nA_n;FA)$ in $\N^{\C\t}$
are in fact determined by the $\al_A = \al\la \id_A,\cd,\id_A \ra$.
Since the arrows $\la \id_A,\cd,\id_A \ra$ in $\C\t$ are fixed by any permutation of the domain,
the formula (\ref{symexp}) gives
\eq
\label{symseqexp}
(\si\al)_A = \si\al_A
\eeq


\section{The monoidal closed structure of the category of symmetric multicategories}
\label{mon}

In the previous section we have studied the cartesian structure of $\Mlt$ and $\sMlt$.
We now consider a non-cartesian symmetric monoidal closed structure  
\[
(\sMlt,1^-,\otm_{BV})
\]
on symmetric multicategories; the symmetry hypothesis is necessary.
The \"Boardman-Vogt" tensor product $\otm_{BV}$ is
described in \cite{weiss}, \cite{moerdijk} and \cite{trova} and has its roots
in \cite{boardman} but also (for the one-sorted case) in the \"Kronecker" tensor product
of Lawvere theories of \cite{freyd}.
The set of objects of $\M\otm_{BV}\N$ (or simply $\M\otm\N$) is the product
of those of $\M$ and $\N$ while its arrows are generated by the arrows 
\[
f_X:\la A_1,X \ra,\cd,\la A_n,X \ra \to \la A,X \ra    \qv   h_A:\la A,X_1 \ra,\cd,\la A,X_m \ra \to \la A,X \ra
\]
(for any $f:A_1,\cd,A_n\to A$ in $\M$ and $X\in\N$, 
and any $h:X_1,\cd,X_m\to X$ in $\N$ and $A\in\M$ respectively)
with the obvious relations assuring that we have functors in $\sMlt$
\[
(-)_X:\M \to \M\otm\N  \qv (-)_A:\N \to \M\otm\N 
\]
and a commuting relation for any pair of arrows $f:A_1,\cd,A_n\to A$ in $\M$ 
and $h:X_1,\cd,X_m\to X$ in $\N$:
\[
f_X(h_{A_1},\cd,h_{A_n}) = \si h_A(f_{X_1},\cd,f_{X_m})
\]
where $\si$ is the obvious permutation.
For instance, if $f:A,B,C\to D$ and $h:X,Y\to Z$, 
one has the following equality in $\M\otm\N$:
\[
\xymatrix@R=2.5pc@C=0.1pc{
\la A,X\ra \ar@/_/@{-}[dr]    &\la A,Y\ra  \ar@{-}[d]  &\la B,X\ra \ar@/_/@{-}[dr]
&\la B,Y\ra \ar@{-}[d] & \la C,X\ra \ar@/_/@{-}[dr]   & \la C,Y\ra\ar@{-}[d] \\
                                                  & h_A \ar[d]                     &&       h_B \ar[d]                                & & h_C \ar[d]                     &     \\
                                                  & \la A,Z\ra  \ar@{-}@/_/[drr]   &&  \la B,Z\ra \ar@{-}[d]      && \la C,Z\ra  \ar@{-}@/^/[dll]  &       \\
                                                  &                                      &&  f_Z \ar[d]                 &                                     &        =            \\
                                                  &                                     & &  \la D,Z\ra                                     &                                     &
}
\xymatrix@R=2pc@C=0.1pc{
\la A,X\ra \ar@{-}[d]    &\la A,Y\ra  \ar@{-}[drr]  &\la B,X\ra \ar@{-}[dl]
&\la B,Y\ra \ar@{-}[dr] & \la C,X\ra \ar@{-}[dll]   & \la C,Y\ra\ar@{-}[d] \\
\la A,X\ra \ar@/_/@{-}[dr]    &\la B,X\ra  \ar@{-}[d]  &\la C,X\ra \ar@/^/@{-}[dl]
&\la A,Y\ra \ar@/_/@{-}[dr] & \la B,Y\ra \ar@{-}[d]   & \la C,Y\ra\ar@/^/@{-}[dl] \\
                                                  & f_X \ar[d]                     &&                                     & f_Y \ar[d]                     &     \\
                                                  & \la D,X\ra  \ar@{-}@/_/[drr]   &&        & \la D,Y\ra  \ar@{-}@/^/[dl]  &       \\
                                                  &                                      &&  h_D \ar[d]                 &                                     &                    \\
                                                  &                                      &&  \la D,Z\ra                                     &                                     &
}
\]

The associated internal hom $[\M,\N]$ has functors $\M\to\N$ as objects while an arrow
$\al:F_1,\cd,F_n\to F$ consists of arrows $\al_A:F_1A,\cd,F_nA\to FA$ ($A\in\M$) 
such that the following commuting condition holds for any arrow $f:A_1,\cd,A_m\to A$ in $\M$:
\eq
\label{mnat}
Ff(\al_{A_1},\cd,\al_{A_m}) = \si \al_A(F_1f,\cd,F_nf)
\eeq
where $\si$ is the obvious permutation.
The composition and the symmetric structure are defined pointwise.

For instance, an arrow $\al:F,G\to H$ is a family of arrows $\al_A:FA,GA\to HA$ ($A\in\M$)
such that for any arrow, say $f:A,B,C\to D$, the following equality holds:
\[
\xymatrix@R=2.5pc@C=0.3pc{
FA \ar@/_/@{-}[dr]    & GA  \ar@{-}[d]  & FB \ar@/_/@{-}[dr] && GB \ar@/^/@{-}[dl] & FC \ar@{-}[d]   & GC \ar@/^/@{-}[dl] \\
                                                  & \al_A \ar[d]                     &&       \al_B \ar[d]                                & & \al_C \ar[d]                     &     \\
                                                  & HA  \ar@{-}@/_/[drr]   &&  HB \ar@{-}[d]      && HC  \ar@{-}@/^/[dll]  &       \\
                                                  &                                      &&  Hf \ar[d]                 &                                     &                   \\
                                                  &                                     & &  HD                                     &                                     &
}
\xymatrix@R=2.2pc@C=1.5pc{ 
&\\            
              &  &         \\
              &  &           \\
              & = & 
}
\xymatrix@R=2pc@C=0.3pc{
FA \ar@{-}[d]    & GA  \ar@{-}[drrr]  & FB \ar@{-}[dl] && GB \ar@{-}[dr] & FC \ar@{-}[dlll]   & GC\ar@{-}[d] \\
FA \ar@/_/@{-}[dr]  & FB  \ar@{-}[d]  & FC \ar@/^/@{-}[dl] && GA \ar@/_/@{-}[dr] & GB \ar@{-}[d]   & GC \ar@/^/@{-}[dl] \\
                                                  & Ff \ar[d]                     &&                                     && Gf \ar[d]                          \\
                                                  & FD  \ar@{-}@/_/[drr]   &&        && GD  \ar@{-}@/^/[dll]  &       \\
                                                  &                                      &&  \al_D \ar[d]                 &                                     &                    \\
                                                  &                                      &&  HD                                     &                                     &
}
\]
\begin{remark}
\label{armon}
In particular, if $M,N,L:1\t \to \M$ are monoids on the objects $A$, $B$ and $C$ respectively,
an arrow $\al:M,N\to L$ in $[1\t,\M]$ is an arrow $\al_\star:A,B\to C$ in $\M$ such that
\[
\xymatrix@R=2pc@C=0.2pc{
\\
& A \ar@/_/@{-}[dr]  & &  B\ar@/^/@{-}[dl]    &  & A \ar@/_/@{-}[dr] & & B \ar@/^/@{-}[dl]      \\
&& \al_\star \ar[d] && && \al_\star \ar[d]  \\
&& C \ar@/_/@{-}[drr] &&              & & C \ar@/^/@{-}[dll] \\
&&                         &&  m^2_C \ar[d]                           &&                           \\
&&                        &&  C                                       &&
}
\xymatrix@R=2pc@C=1pc{
\\
\\
              &  &                       \\
              & = & 
}
\xymatrix@R=2pc@C=0.2pc{
& A \ar@{-}[d] & &   B \ar@{-}[drr] & & A \ar@{-}[dll] & & B \ar@{-}[d]    \\
& A \ar@/_/@{-}[dr]  & &  A\ar@/^/@{-}[dl]    &  & B \ar@/_/@{-}[dr] & & B \ar@/^/@{-}[dl]      \\
&&  m^2_A\ar[d] && && m^2_B\ar[d]   \\
&& A \ar@/_/@{-}[drr] &&              & & B \ar@/^/@{-}[dll] \\
&&                         &&  \al_\star \ar[d]                           &&                           \\
&&                        &&  C                                       &&
}
\]
(where for instance $m^2_A = Mm^2$ is the image of the unique 2-ary arrow of $1\t$)
and
\[
\xymatrix@R=2pc@C=0.2pc{
\\
&&                         &&  m^0_C \ar[d]                           &&                           \\
&&                        &&  C                                       &&
}
\xymatrix@R=2pc@C=1pc{
\\
              &  &                       \\
              & = & 
}
\xymatrix@R=2pc@C=0.2pc{
&&  m^0_A\ar[d] && && m^0_B\ar[d]   \\
&& A \ar@/_/@{-}[drr] &&              & & B \ar@/^/@{-}[dll] \\
&&                         &&  \al_\star \ar[d]                           &&                           \\
&&                        &&  C                                       &&
}
\]
\end{remark}

\begin{remark}
The multicategory $[\M,\N]$ is studied in \cite{tronin}, where its arrows are called 
natural multitransformations (of multifunctors).
Note that 
\[
[\M,\N]_- \iso \sMlt(\M,\N)
\]
where $\sMlt$ is given the usual 2-category structure.
Furthermore,
\[
[\C^-,\D^-] \iso (\D\,^\C)^-
\]
\end{remark}

\begin{prop}
For any $\C\in\Cat$ and $\N\in\sMlt$,
\[
[\C^-,\N] \iso \N^{\C\t}   \qv  \C^-\otm\N \iso \C\t\tm\N
\]
In particular, $\C^-\otm 1\t \iso \C\t$ and $\C^-\otm\D^- \iso (\C\tm\D)^-$.
\end{prop}
\pf
Since $\sMlt(\C^-,\N) \iso \Cat(\C,\N_-) \iso \Cat((\C\t)_-,\N_-)$,
the objects of $[\C^-,\N]$ and of $\N^{\C\t}$ are both the functors $\C\to\N_-$.

Since the condition (\ref{mnat}) relative to an arrow $\al$ in $[\M,\N]$ reduces,
for the unary $f$ in $\C^-$, to the condition (\ref{seqexp2}) relative to an arrow $\al$ in $\N^{\C\t}$,
also the arrows coincide.
Furthermore, by (\ref{symseqexp}), the symmetric structure is the pointwise one in both cases.

The second isomorphism then follows by adjunction.
\epf

\begin{corol}
The 2-category $\sMlt$ is tensored  and cotensored (as a $\Cat$-enriched category).
The tensor and the cotensor of $\M$ by $\C\in\Cat$ are given respectively by $\C\t\tm\M$
and $\M^{\C\t}$.
\end{corol}
\pf
\[
\sMlt(\C\t\tm\M,\N) \iso \sMlt(\C^-\otm\M,\N) \iso \sMlt(\C^-,[\M,\N]) \iso 
\]
\[
\iso \Cat(\C,[\M,\N]_-) \iso \Cat(\C,\sMlt(\M,\N))
\]

\[
\sMlt(\N,\M^{\C\t}) \iso \sMlt(\N,[\C^-,\M]) \iso \sMlt(\C^-,[\N,\M]) \iso 
\]
\[
\iso \Cat(\C,[\N,\M]_-) \iso \Cat(\C,\sMlt(\N,\M))
\]
(In fact, the proposition holds true in the non-symmetric context too.)
\epf



\subsection{Characterizations of sequential multicategories}

We are now in a position to characterize sequential multicategories in a more abstract way.
Following \cite{tronin}, we define the \"center" of $\M\in\sMlt$ as the full sub-multicategory
(in fact an operad) of $[\M,\M]$ generated by the object $\id_\M:\M\to\M$.
Thus, an $n$-ary arrow $\al$ in the center of $\M$ consists of $n$-ary arrows $\al_A:A,\cd,A\to A$ ($A\in\M$) 
such that for each $f:A_1,\cd,A_m\to A$
\[
f(\al_{A_1},\cd,\al_{A_m}) = \si \al_A(f,\cd,f)
\]
Accordingly, a \"central monoid" in $\M$ is a commutative monoid in the center of $\M$,
that is a functor $1\t\to[\M,\M]$  which restricts along $1^-\to 1\t$ to the name of the identity $1^-\to[\M,\M]$.
\begin{remark}
 By Hilton-Eckmann arguments, central monoids (when they exist) are unique
and unital magmas in the center of $\M$ are in fact central monoids.
These arguments are best developed in an operad (a one object symmetric multicategory), 
where two arrows $f$ and $g$ (of arity $n$ and $m$) are said to commute if
\[
f(g,\cd,g) = \si g(f,\cd,f)
\]
(where $\si$ is the usual canonical permutation on $nm$ symbols).
Clearly, two 0-ary arrows commute iff they are equal.
If two unital magmas $m$ and $n$ have commuting units ($m^0 = n^0$)
and multiplications ($m^2$ and $n^2$), then they coincide:
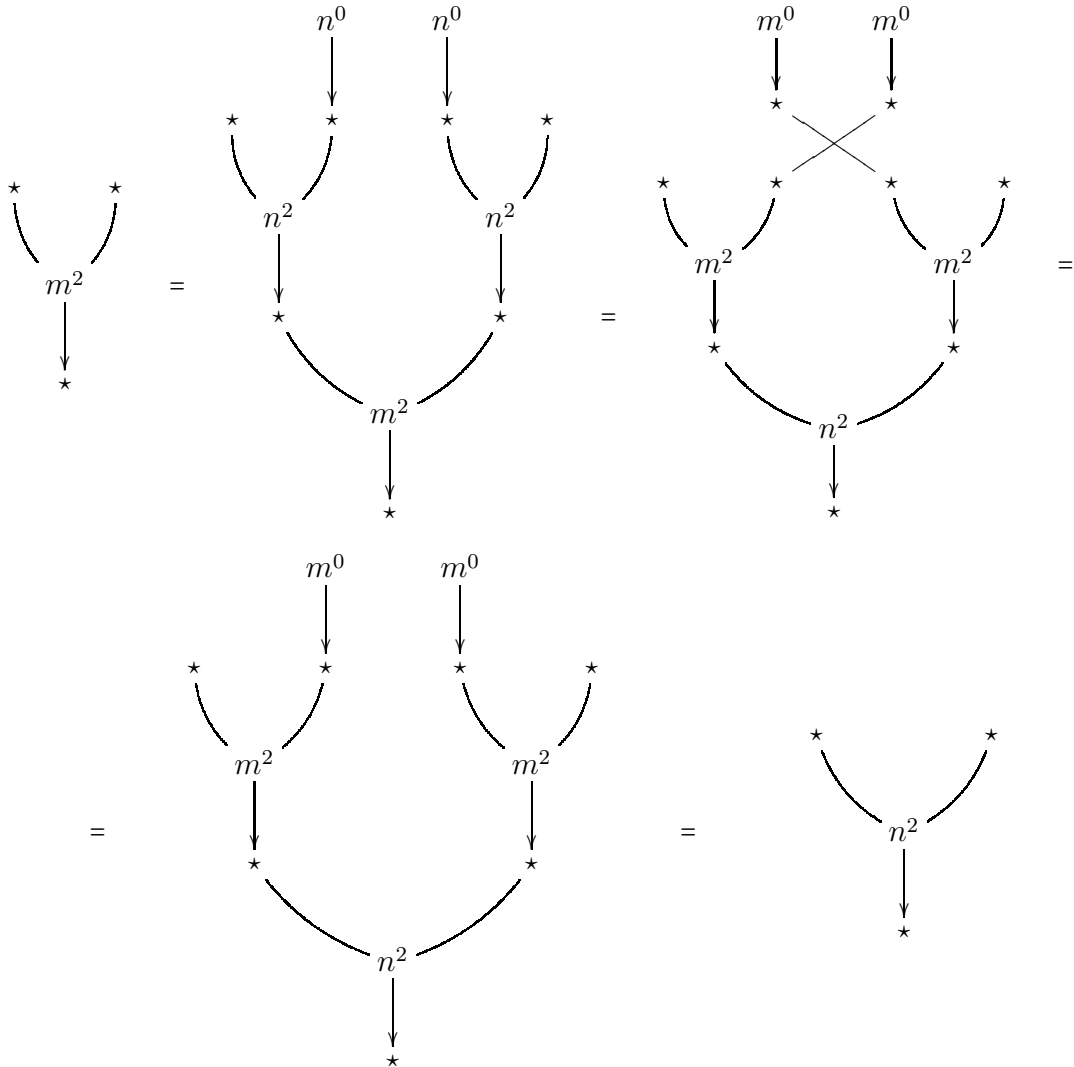
\[
\xymatrix@R=2pc@C=0.2pc{
\\
\\
                                                  & \star \ar@{-}@/_/[dr]   &                                         & \star \ar@{-}@/^/[dl]  &                         \\
                                                  &                                            &  m^2 \ar[d]                 &                                     & &    =               \\
                                                  &                                            &  \star                              &                                     &
}
\xymatrix@R=2pc@C=0.2pc{
&                                                   &                                     &   n^0 \ar[d]                 &&     n^0 \ar[d]                                                           \\
&\star  \ar@/_/@{-}[dr]                &                                     &   \star  \ar@/^/@{-}[dl]   && \star\ar@{-}@/_/[dr]              && \star\ar@/^/@{-}[dl]           \\
&                                                 & n^2 \ar[d]                     &&                                       & & n^2 \ar[d]                     &     \\
&                                                  & \star \ar@{-}@/_/[drr]   &&                                       &  & \star \ar@{-}@/^/[dll]  &&&                =         \\
&                                                  &                                      &&  m^2 \ar[d]                 &                                     &                    \\
&                                                  &                                     & &  \star                                    &                                     &
}
\xymatrix@R=1.5pc@C=0.2pc{
&                                                   &                                     &   m^0 \ar[d]                 &&     m^0 \ar[d]                                                           \\
&                                                  &                                      &  \star \ar@{-}[drr]                      && \star    \ar@{-}[dll]                                                          \\
& \star \ar@/_/@{-}[dr]                &                                     &   \star  \ar@/^/@{-}[dl]   && \star\ar@{-}@/_/[dr]      && \star\ar@/^/@{-}[dl]           \\
&                                                  & m^2 \ar[d]                     &&                                       & & m^2 \ar[d]                    && &    = \\
&                                                  & \star \ar@{-}@/_/[drr]   &&                                       &  & \star \ar@{-}@/^/[dll]  &                         \\
&                                                  &                                      &&  n^2 \ar[d]                 &                                     &                    \\
&                                                  &                                     & &  \star                                    &                                     &
}
\]
\[
\xymatrix@R=2.2pc@C=1.5pc{ 
&\\            
              &  &         \\
              &  &           \\
              & = & 
}
\xymatrix@R=2pc@C=0.5pc{
                                                   &                                     &   m^0 \ar[d]                 &&     m^0 \ar[d]                                                           \\
\star  \ar@/_/@{-}[dr]                &                                     &   \star  \ar@/^/@{-}[dl]   && \star\ar@{-}@/_/[dr]              && \star \ar@/^/@{-}[dl]           \\
                                                  & m^2 \ar[d]                     &&                                       & & m^2 \ar[d]                     &     \\
                                                  & \star \ar@{-}@/_/[drr]   &&                                       &  & \star \ar@{-}@/^/[dll]  &                         \\
                                                  &                                      &&  n^2 \ar[d]                 &                                     &                    \\
                                                  &                                     & &  \star                                    &                                     &
}
\xymatrix@R=2.2pc@C=1.5pc{   
     & \\          
              &  &         \\
              &  &           \\
              & = & 
}
\xymatrix@R=2pc@C=0.5pc{
&  \\
&\\
                                                  & \star \ar@{-}@/_/[drr]   &&                                       &  & \star \ar@{-}@/^/[dll]  &                         \\
                                                  &                                      &&  n^2 \ar[d]                 &                                     &                    \\
                                                  &                                     & &  \star                                    &                                     &
}
\]
Similarly one sees that, in the above hypothesis, $m=n$ is commutative:
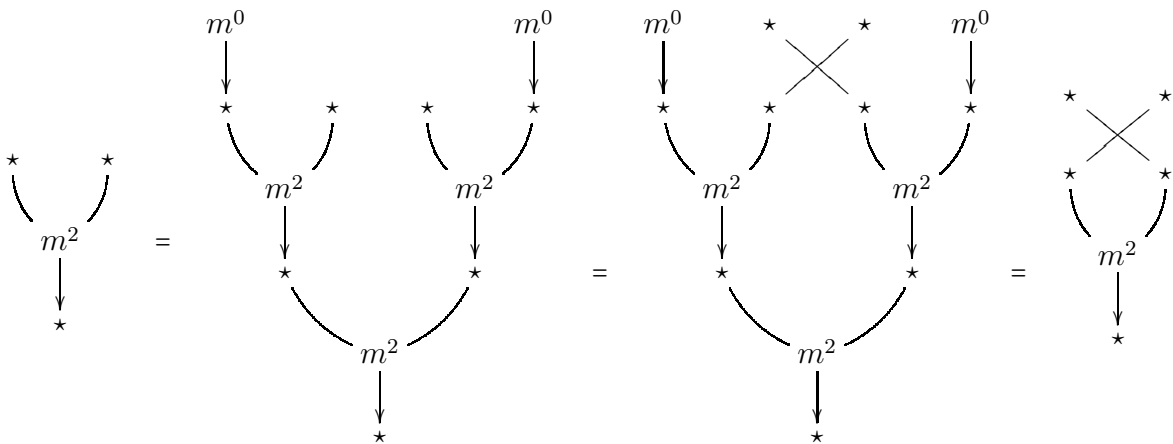
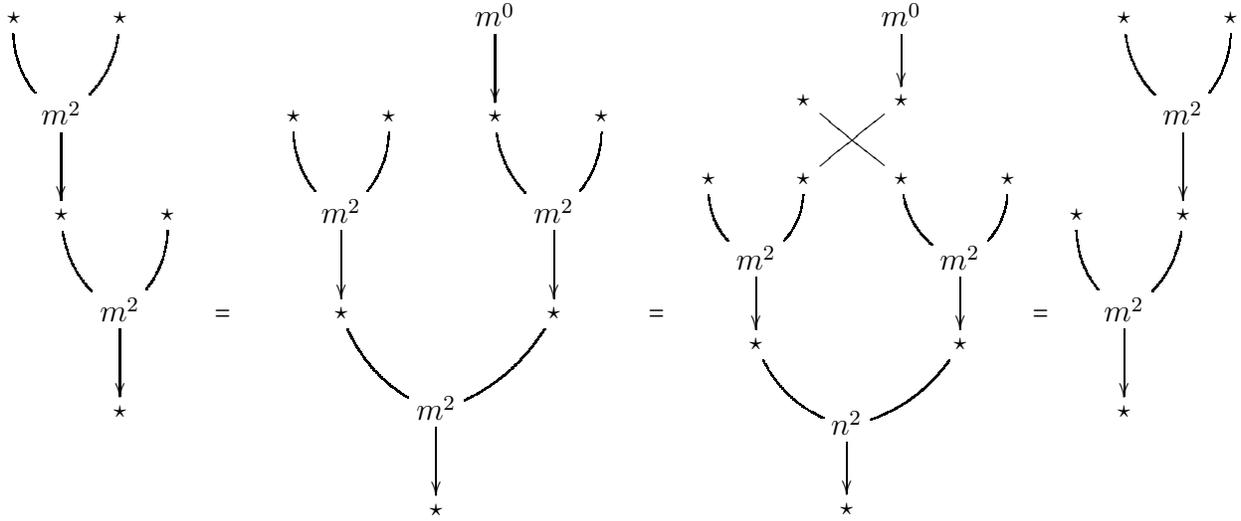
\[
\xymatrix@R=1.5pc@C=0.1pc{
\\
\\
                                                  & \star \ar@{-}@/_/[dr]   &                                         & \star \ar@{-}@/^/[dl]  &                         \\
                                                  &                                            &  m^2 \ar[d]                 &                                     & &    =               \\
                                                  &                                            &  \star                              &                                     &
}
\xymatrix@R=1.5pc@C=0.1pc{
&                                                          m^0 \ar[d]                            &&    &&  &&     m^0 \ar[d]                                                           \\
&\star  \ar@/_/@{-}[dr]          &                                     &   \star  \ar@/^/@{-}[dl]   && \star\ar@{-}@/_/[dr]              && \star\ar@/^/@{-}[dl]           \\
&                                                 & m^2 \ar[d]                     &&                                       & & m^2 \ar[d]                     &     \\
&                                                  & \star \ar@{-}@/_/[drr]   &&                                       &  & \star \ar@{-}@/^/[dll]  &&&                =         \\
&                                                  &                                      &&  m^2 \ar[d]                 &                                     &                    \\
&                                                  &                                     & &  \star                                    &                                     &
}
\xymatrix@R=1.5pc@C=0.1pc{
&      m^0 \ar[d]       &                                      &  \star \ar@{-}[drr]                      && \star    \ar@{-}[dll]        &&        m^0 \ar[d]               \\
& \star \ar@/_/@{-}[dr]                &                                     &   \star  \ar@/^/@{-}[dl]   && \star\ar@{-}@/_/[dr]      && \star\ar@/^/@{-}[dl]           \\
&                                                  & m^2 \ar[d]                     &&                                       & & m^2 \ar[d]                    &&  \\
&                                                  & \star \ar@{-}@/_/[drr]   &&                                       &  & \star \ar@{-}@/^/[dll]  &&  =                 \\
&                                                  &                                      &&  m^2 \ar[d]                 &                                     &                    \\
&                                                  &                                     & &  \star                                    &                                     &
}
\xymatrix@R=1.5pc@C=0.1pc{
\\
                                                  & \star \ar@{-}[drr]   &                                         & \star \ar@{-}[dll]  &                         \\
                                                  & \star \ar@{-}@/_/[dr]   &                                         & \star \ar@{-}@/^/[dl]  &                         \\
                                                  &                                            &  m^2 \ar[d]                 &                                     & &                   \\
                                                  &                                            &  \star                              &                                     &
}
\]
and associative:
\[
\xymatrix@R=2pc@C=0.1pc{
\star  \ar@/_/@{-}[dr]                &                                     &   \star  \ar@/^/@{-}[dl]   &&              &          \\
                                                 & m^2 \ar[d]                     &&                                       & &                      &     \\
                                                  & \star \ar@{-}@/_/[dr]   &&                               \star \ar@{-}@/^/[dl]  &         \\
                                                  &                                      &  m^2 \ar[d]                 &                                     &&      =              \\
                                                 &                                     &   \star                                    &                                     &
}
\xymatrix@R=2pc@C=0.1pc{
&                                                   &                                     &                   &&     m^0 \ar[d]                                                           \\
&\star  \ar@/_/@{-}[dr]                &                                     &   \star  \ar@/^/@{-}[dl]   && \star\ar@{-}@/_/[dr]              && \star\ar@/^/@{-}[dl]           \\
&                                                 & m^2 \ar[d]                     &&                                       & & m^2 \ar[d]                     &     \\
&                                                  & \star \ar@{-}@/_/[drr]   &&                                       &  & \star \ar@{-}@/^/[dll]  &&&                =         \\
&                                                  &                                      &&  m^2 \ar[d]                 &                                     &                    \\
&                                                  &                                     & &  \star                                    &                                     &
}
\xymatrix@R=1.5pc@C=0.1pc{
&                                                   &                                     &                  &&     m^0 \ar[d]                                                           \\
&                                                  &                                      &  \star \ar@{-}[drr]                      && \star    \ar@{-}[dll]                                                          \\
& \star \ar@/_/@{-}[dr]                &                                     &   \star  \ar@/^/@{-}[dl]   && \star\ar@{-}@/_/[dr]      && \star\ar@/^/@{-}[dl]           \\
&                                                  & m^2 \ar[d]                     &&                                       & & m^2 \ar[d]                     \\
&                                                  & \star \ar@{-}@/_/[drr]   &&                                       &  & \star \ar@{-}@/^/[dll]  &                         \\
&                                                  &                                      &&  n^2 \ar[d]                 &                                     &                    \\
&                                                  &                                     & &  \star                                    &                                     &
}
\xymatrix@R=2pc@C=0.1pc{
&&\star  \ar@/_/@{-}[dr]                &                                     &   \star  \ar@/^/@{-}[dl]             \\
&&                                                 & m^2 \ar[d]                     &&                                                     \\
                                                  & \star \ar@{-}@/_/[dr]   &&                               \star \ar@{-}@/^/[dl]  &         \\
 =                                               &                                      &  m^2 \ar[d]                 &                                     &              \\
                                                 &                                     &   \star                                    &                                     &
}
\]
\end{remark}

The functor $1^-\to 1\t$ induces canonical functors 
\[
[1\t,\M]\to\M      \qv    \M\to 1\t\otm\M
\]
The first one forgets the monoid structure while the second one inserts the arrow $f\in\M$ in $1\t\otm\M$
as $f_\star$ (where $\star$ is the unique object of $1\t$).


\begin{prop}
\label{char}
The following are equivalent for $\M\in\sMlt$:
\begin{enumerate}
\item
$\M$ is sequential, that is $\M\iso\C\t$, for a $\C\in\Cat$.
\item
$\M\iso(\M_-)\t$.
\item
$\M$ has a central monoid $1\t\to[\M,\M]$.
\item
$\M$ has a central unital magma.
\item
$[1\t,\M]\to\M$ is an iso
(that is $\M$ is orthogonal to $1^-\to 1\t$ in an enriched sense).
\item
$\M\to 1\t\otm \M$ is an iso.
\end{enumerate}
\end{prop}
\pf
Since $(\C\t)_- =\C$, (1) and (2) are equivalent.
As mentioned above, the equivalence of (3) and (4) follows from a Hilton-Eckmann argument.
If $[1\t,\M]\to\M$ has an inverse $\M\to [1\t,\M]$, one gets (by the closed structure of $\sMlt$)
a central monoid $1\t\to[\M,\M]$. Thus $(5) \imp (3)$, and similarly $(6) \imp (3)$.

Suppose now that $\M$ has a central monoid, given by $n$-ary arrows $m_A^n:A,\cd,A\to A$, 
$A\in\M$ and $n=0,1,2\cd$, such that for any $f:A_1,\cd,A_m\to A$
\[
f(m^n_{A_1},\cd,m^n_{A_m}) = \si m^n_A(f,\cd,f)
\]
(which implies in particular that any unary arrow is a monoid morphism).
We will define identity-on-objects functors 
\[ 
F:\M\to(\M_-)\t   \qv   G:(\M_-)\t \to \M
\]
which are each other's inverse, thus proving that $(3) \imp (2)$.

For an arrow $f:A_1,\cd,A_n \to A$ in $\M$, let $Ff = \la f_1,\cd,f_n \ra$ in $(\M_-)\t$ be given by 
\[
f_i  = f(m^0_{A_1},m^0_{A_2},\cd,m^0_{A_{i-1}}, A_i,m^0_{A_{i+1}},\cd,m^0_{A_{n-1}},m^0_{A_n})
\]
that is, $f_i$ is obtained by substituting the 0-ary arrow $m^0$ (the monoid unit) in any object of the 
domain but the $i$-th one.

In the other direction, if $f_i:A_i\to A$ ($i=1,\cd,n$), let $G\la f_1,\cd,f_n \ra = m_A^n(f_1,\cd,f_n)$.
Then one easily checks that $F$ and $G$ are indeed functor and that they are each other's inverse.
For instance, for 2-ary arrows we have
\[
\xymatrix@R=2pc@C=0.5pc{
                                                   &                                     &   m^0_B \ar[d]                 &&     m^0_A \ar[d]                                                           \\
A  \ar@/_/@{-}[dr]                &                                     &   B  \ar@/^/@{-}[dl]   && A\ar@{-}@/_/[dr]              && B\ar@/^/@{-}[dl]           \\
                                                  & f \ar[d]                     &&                                       & & f \ar[d]                     &     \\
                                                  & C \ar@{-}@/_/[drr]   &&                                       &  & C \ar@{-}@/^/[dll]  &                         \\
                                                  &                                      &&  m^2_C \ar[d]                 &                                     &                    \\
                                                  &                                     & &  C                                    &                                     &
}
\xymatrix@R=2.2pc@C=1.5pc{ 
&\\            
              &  &         \\
              &  &           \\
              & = & 
}
\xymatrix@R=2pc@C=0.5pc{
                                                   &                                     &   m^0_B \ar[d]                 &&     m^0_A \ar[d]                                                           \\
                                                   &                                      &  B \ar@{-}[drr]                      && A    \ar@{-}[dll]                                                          \\
A  \ar@/_/@{-}[dr]                &                                     &   A  \ar@/^/@{-}[dl]   && B\ar@{-}@/_/[dr]              && B\ar@/^/@{-}[dl]           \\
                                                  & m^2_A \ar[d]                     &&                                       & & m^2_B \ar[d]                     &     \\
                                                  & A \ar@{-}@/_/[drr]   &&                                       &  & B \ar@{-}@/^/[dll]  &                         \\
                                                  &                                      &&  f \ar[d]                 &                                     &                    \\
                                                  &                                     & &  C                                    &                                     &
}
\xymatrix@R=2.2pc@C=1.5pc{
&\\             
              &  &         \\
              &  &           \\
              & = & 
}
\]
\[
\xymatrix@R=2.2pc@C=1.5pc{ 
&\\            
              &  &         \\
              &  &           \\
              & = & 
}
\xymatrix@R=2pc@C=0.5pc{
                                                   &                                     &   m^0_A \ar[d]                 &&     m^0_B \ar[d]                                                           \\
A  \ar@/_/@{-}[dr]                &                                     &   A  \ar@/^/@{-}[dl]   && B\ar@{-}@/_/[dr]              && B\ar@/^/@{-}[dl]           \\
                                                  & m^2_A \ar[d]                     &&                                       & & m^2_B \ar[d]                     &     \\
                                                  & A \ar@{-}@/_/[drr]   &&                                       &  & B \ar@{-}@/^/[dll]  &                         \\
                                                  &                                      &&  f \ar[d]                 &                                     &                    \\
                                                  &                                     & &  C                                    &                                     &
}
\xymatrix@R=2.2pc@C=1.5pc{   
     & \\          
              &  &         \\
              &  &           \\
              & = & 
}
\xymatrix@R=2pc@C=0.5pc{
&  \\
&\\
                                                  & A \ar@{-}@/_/[drr]   &&                                       &  & B \ar@{-}@/^/[dll]  &                         \\
                                                  &                                      &&  f \ar[d]                 &                                     &                    \\
                                                  &                                     & &  C                                    &                                     &
}
\]
and in the other direction:
\[
\xymatrix@R=2pc@C=0.5pc{
                                                   &                                         &&                                   &&     m^0_B \ar[d]                                                           \\
                                                  &  A  \ar@{-}[d]               &&                                     && B\ar@{-}[d]              &          \\
                                                  & f_1 \ar[d]                     &&                                       & & f_2 \ar[d]                     &     \\
                                                  & C \ar@{-}@/_/[drr]   &&                                       &  & C \ar@{-}@/^/[dll]  &                         \\
                                                  &                                      &&  m^2_C \ar[d]                 &                                     &                    \\
                                                  &                                     & &  C                                    &                                     &
}
\xymatrix@R=2.2pc@C=1.5pc{ 
&\\            
              &  &         \\
              &  &           \\
              & = & 
}
\xymatrix@R=2pc@C=0.5pc{
                                                   &                                         &&                                   &&                                                              \\
                                                  &  A  \ar@{-}[d]               &&                                     &&               &          \\
                                                  & f_1 \ar[d]                     &&                                       & & m^0_C \ar[d]                     &     \\
                                                  & C \ar@{-}@/_/[drr]   &&                                       &  & C \ar@{-}@/^/[dll]  &                         \\
                                                  &                                      &&  m^2_C \ar[d]                 &                                     &                    \\
                                                  &                                     & &  C                                    &                                     &
}
\xymatrix@R=2.2pc@C=1.5pc{
&\\             
              &  &         \\
              &  &           \\
              & = & 
}
\xymatrix@R=2pc@C=0.5pc{
& \\
                                                   &                                                \\
                                                  &  A  \ar@{-}[d]                        \\
                                                  & f_1 \ar[d]                                  \\
                                                  & C                       
}
\]

To prove that $(1) \imp (5)$, note that $m_A = \la \id_A,\cd,\id_A \ra$ is the unique 
monoid structure on $A$ in $\C\t$. Indeed, if for instance $m^2_A = \la f, g \ra$, then
$f = m^2_A(\la\,\ra,\id_A) = \id_A$.
Furthermore, $m_A$ clearly commutes with any arrow in $\C\t$.

Finally, we prove that $(1) \imp (6)$. 
Note first that for any $\M\in\sMlt$, by the definition of the Boardman-Vogt tensor product,
$1\t\otm \M$ has the same objects as $\M$, while its arrows are generated by the $f_\star$ 
(with $f$ in $\M$) and by the $m^n_A:A,\cd,A\to A$ ($A\in\M$),
where $m^n$ is the unique $n$-ary arrow in $1\t$.
These $m^n_A$ form a commutative monoid on $A$ and commute with all the $f_\star$.
Thus $m_A$ is a central monoid in $1\t\otm\M$ and (since we already know that $(3) \imp (1)$)
$1\t\otm\M$ is sequential.
Now suppose that $\M$ is itself sequential. To show that $1\t\otm\M \iso \M$ it is enough to
compare them with sequential categories $\L$:
\[    
\begin{array}{c}
1\t\otm\M \to \L \\ \hline
\M \to [1\t,\L] \\ \hline
\M \to \L
\end{array}
\]
(since we already know that $(1) \imp (5)$).
\epf

\begin{corol}
If $\L$ is sequential and $\M\in\sMlt$ then $\M\otm\L$, $[\M,\L]$ and $[\L,\M]$ are also sequential.
\end{corol}
\pf

\[
(\M\otm \L)\otm 1\t \iso \M\otm(\L\otm 1\t) \iso \M\otm\L
\]

\[
[ 1\t , [\M , \L]] \iso [\M , [1\t , \L]] \iso [\M , \L] 
\]

\[
[ 1\t , [\L , \M]] \iso [1\t\otm\L , \M] \iso [\L , \M] 
\]

\epf

\begin{remark}
Note that the central monoid in $[\M,\L]$ is given by $(m_F)_A = m_{FA}$, 
while that in $[\L,\M]$ is given by $(m_F)_A = Fm_A$ (where $m_A$ is the central
monoid in $\L$).
\end{remark}

\begin{corol}
\label{mon2}
The functors $1\t\otm - $ and $[1\t,-]$ give respectively a reflection and a coreflection
of $\sMlt$ in $\Seq$.
\end{corol}
\pf
If $\L$ is sequential, then for any $\M\in\sMlt$ we have natural isomorphisms
\[    
\begin{array}{c}
1\t\otm\M \to \L \\ \hline
\M \to [1\t,\L] \\ \hline
\M \to \L
\end{array}
\qv
\begin{array}{c}
\L \to [1\t,\M] \\ \hline
1\t\otm\L \to \M \\ \hline
\L \to \M
\end{array}
\]
\epf

\begin{corol}
\label{lradj}
$(-)\t:\Cat\to\sMlt$ has both a left and a right adjoint.
\end{corol}
\pf
Just compose the above (co)reflection with an equivalence $\Cat\simeq\Seq$;
for instance:
\[
(1\t\otm - )_-  \adj (-)\t \adj [1\t,-]_- \iso \cMon(-)
\]
\epf

\begin{remark}
The above results can be placed in a more general perspective.
Consider a monoidal category $(\C,\I,\otm)$ with a terminal object $1\in\C$.
Then we have a (unique) monoid structure on $1$ (with respect to $\otm$) 
and thus also a monad $1\otm(-)$ on $\C$.
For instance, for $\C= (\Set,0,+)$ one gets the monad for pointed sets.

If $\C$ is closed, we have $1\otm(-) \adj [1,-]$, so that a $[1,-]$ is a comonad 
on $\C$ with the same algebras.
Suppose furthermore that $1\otm 1 \to 1$ is an iso. Then the (co)monad is idempotent
and the isomorphic categories of algebras isolate a full subcategory $\D\subseteq\C$ which
is both reflective and coreflective.

When $\C = \sMlt$ we get $\D\iso\Seq\simeq\Cat$.
For another instance, consider the poset $\C$ of the subsets of a monoid $M$
with the usual tensor (truth values convolution on the discrete $M$).
Then $\D$ is given by the monoid ideals of $M$.
(Note that if $\C$ is cartesian or $\I\iso 1$ one gets $\D\iso\C$, while for pointed sets or abelian groups
$\D\subseteq\C$ is the zero object inclusion.)
\end{remark}

\begin{corol}
\[
(\C\tm\D)\t \iso \C\t\otm \D\t  \qv (\D^{\,\C})\t \iso [\C\t,\D\t]
\]
\end{corol}
(Recall however that $1\t \not\iso \I = 1^-$.)

\pf
\[
\C\t\otm \D\t \iso (\C^-\otm 1\t)\otm\D\t  \iso \C^-\otm (1\t \otm\D\t) \iso \C^-\otm\D\t \iso \C\t\tm\D\t 
\]

\[ 
[\C\t , \D\t ] \iso [\C^-\otm 1\t ,\D\t ] \iso [\C^- , [1\t ,\D\t]] \iso [\C^- , \D\t] \iso \D\t^{\,\C\t} 
\]
\epf

\begin{remark}
In the same way one sees that, for any $\M\in\sMlt$,
\[
\C\t\otm\M \iso \C\t\tm(\M\otm1\t)  \qv [\C\t,\M] \iso [1\t,\M]^{\,\C\t}  \qv [\M,\C\t] \iso \C\t^{\,\M\otm 1\t}
\]
\end{remark}

\begin{corol}
By restricting $\otm$ and $[-,-]$ to $\Seq\subset\sMlt$ one gets a monoidal closed
structure $(\Seq,1\t,\otm)$ which is cartesian and coincides essentially with the
cartesian closed structure on $\Cat$.
\epf
\end{corol}


\subsection{Sequentiality and representability}

In the rest of this section we investigate some consequence of the previous results,
under the hypothesis that some of the multicategories involved are representable.
We begin by recalling what are the sequential representable categories.

\begin{prop}
\label{repseq}
The full sub-multicategory $\sRep\cap\Seq\subset\sMlt$ is equivalent to the category 
of cocartesian monoidal categories (those in which $\I$ is initial and $\otm$ gives sums) 
and to the category of categories with finite sums (and all functors).
\end{prop}
\pf
Indeed, a preuniversal arrow $A_1,\cd,A_n\to A$ in $\C\t$ amounts
to a representation for the functor $\C(A_1,-)\tm\cd\tm\C(A_n,-)$,
that is to a universal cone $A_i\to A_1+\cd + A_n$.
Furthermore, universal cones are closed with respect to composition
and all functors between finite sum categories are lax monoidal
functor between the corresponding cocartesian monoidal categories.
\epf

Next, there is the fact that representability is inherited pointwise
by $[\M,\N]$ from $\N$.

\begin{prop}
\label{rep}
If $\N$ is representable so it is also $[\M,\N]$, for any $\M\in\sMlt$.
\end{prop}
\pf
For simplicity, we concentrate on 2-ary arrows. 
Given $F,G:\M\to\N$ we want to find a universal arrow $u:F,G\to H$ in $[\M,\N]$.
Let $u_A:FA,GA\to HA$ ($A\in\M$) be universal in $\N$.
Since we want to extend the assignment $A\mapsto HA$ to a functor $H:\M\to\N$
in such a way that the $u_A$ become an arrow in $[\M,\N]$, the effect of $H$ on arrows
is forced by the naturality conditions.
Namely, for $f:A_1,\cd,A_n\to A$ in $\M$, the condition (\ref{mnat}) 
\[
Hf(u_{A_1},\cd,u_{A_n}) = \si u_A(Ff,Gf)
\]
defines (by the universality of the $u_A$) a unique arrow $Hf$ in $\N$.
For instance, given the arrow $f:A,B\to C$ in $\M$, $Hf$ is defined by
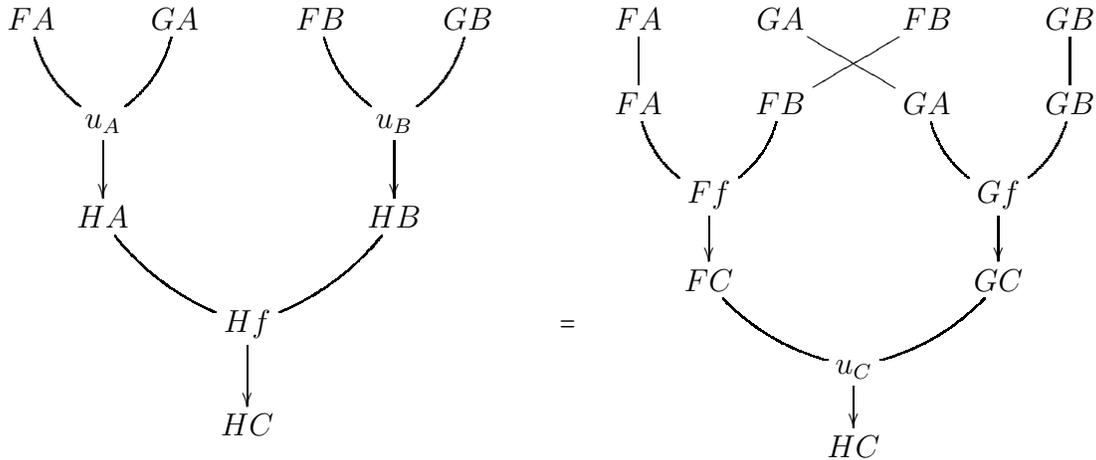
\[
\xymatrix@R=2pc@C=0.2pc{
& FA \ar@/_/@{-}[dr]  & &  GA\ar@/^/@{-}[dl]    &  & FB \ar@/_/@{-}[dr] & & GB \ar@/^/@{-}[dl]      \\
&&  u_A\ar[d] && && u_B\ar[d]   \\
&& HA \ar@/_/@{-}[drr] &&              & & HB \ar@/^/@{-}[dll] \\
&&                         &&  Hf \ar[d]                           &&&&&&          =                 \\
&&                        &&  HC                                       &&
}
\xymatrix@R=1.5pc@C=0.2pc{
& FA \ar@{-}[d]  & &  GA\ar@{-}[drr]    &  & FB \ar@{-}[dll] & & GB \ar@{-}[d]      \\
& FA \ar@/_/@{-}[dr]  & &  FB\ar@/^/@{-}[dl]    &  & GA \ar@/_/@{-}[dr] & & GB \ar@/^/@{-}[dl]      \\
&&  Ff\ar[d] && && Gf\ar[d]   \\
&& FC \ar@/_/@{-}[drr] &&              & & GC \ar@/^/@{-}[dll] \\
&&                         &&  u_C\ar[d]                           &&                           \\
&&                        &&  HC                                       &&
}
\]
Now, drawing the apposite diagrams, it is routine to check that $H$ is indeed a functor.
For instance, let us show that $H$ preserves composition of 2-ary arrows:
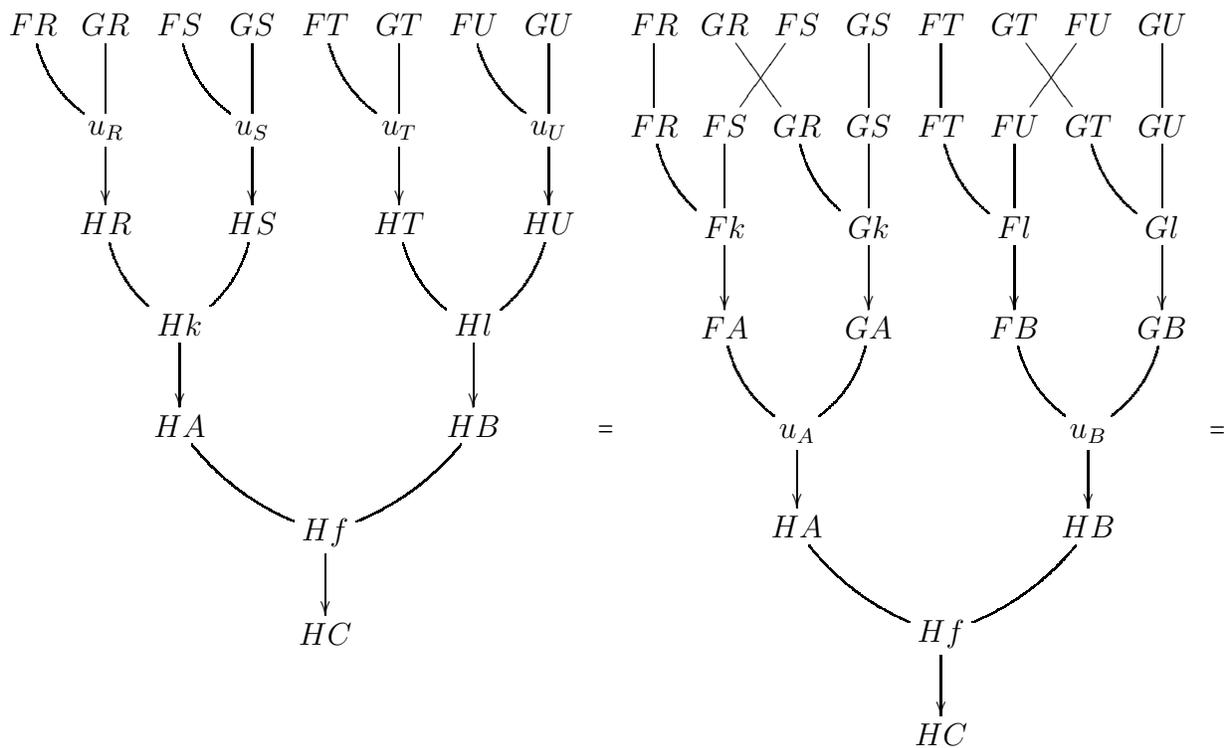
\[
\xymatrix@R=2pc@C=0.2pc{
FR \ar@/_/@{-}[dr] & GR \ar@{-}[d]    &   FS \ar@/_/@{-}[dr] & GS \ar@{-}[d]    & FT \ar@/_/@{-}[dr] & GT \ar@{-}[d]    
& FU \ar@/_/@{-}[dr] & GU \ar@{-}[d]                     \\
& u_R \ar[d]  & &  u_S \ar[d]    &  & u_T \ar[d] & & u_U \ar[d]      \\
& HR \ar@/_/@{-}[dr]  & &  HS\ar@/^/@{-}[dl]    &  & HT \ar@/_/@{-}[dr] & & HU \ar@/^/@{-}[dl]      \\
&&  Hk\ar[d] && && Hl\ar[d]  \\
&& HA \ar@/_/@{-}[drr] &&              & & HB \ar@/^/@{-}[dll] && = \\
&&                         &&  Hf \ar[d]                           &&                           \\
&&                        &&  HC                                       &&
}
\xymatrix@R=2pc@C=0.2pc{
FR \ar@{-}[d] & GR \ar@{-}[dr]    &   FS \ar@{-}[dl] & GS \ar@{-}[d]    & FT \ar@{-}[d] & GT \ar@{-}[dr]    
& FU \ar@{-}[dl] & GU \ar@{-}[d]                     \\
FR \ar@/_/@{-}[dr] & FS \ar@{-}[d]    &   GR \ar@/_/@{-}[dr] & GS \ar@{-}[d]    & FT \ar@/_/@{-}[dr] & FU \ar@{-}[d]    
& GT \ar@/_/@{-}[dr] & GU \ar@{-}[d]                     \\
& Fk \ar[d]  & &  Gk \ar[d]    &  & Fl \ar[d] & & Gl \ar[d]      \\
& FA \ar@/_/@{-}[dr]  & &  GA\ar@/^/@{-}[dl]    &  & FB \ar@/_/@{-}[dr] & & GB \ar@/^/@{-}[dl]      \\
&&  u_A\ar[d] && && u_B\ar[d]  && = \\
&& HA \ar@/_/@{-}[drr] &&              & & HB \ar@/^/@{-}[dll] \\
&&                         &&  Hf \ar[d]                           &&                           \\
&&                        &&  HC                                       &&
}
\]
\[
\xymatrix@R=2pc@C=0.2pc{
FR \ar@{-}[d] & GR \ar@{-}[dr]    &   FS \ar@{-}[dl] & GS \ar@{-}[d]    & FT \ar@{-}[d] & GT \ar@{-}[dr]    
& FU \ar@{-}[dl] & GU \ar@{-}[d]                     \\
FR \ar@/_/@{-}[dr] & FS \ar@{-}[d]    &   GR \ar@/_/@{-}[dr] & GS \ar@{-}[d]    & FT \ar@/_/@{-}[dr] & FU \ar@{-}[d]    
& GT \ar@/_/@{-}[dr] & GU \ar@{-}[d]                     \\
& Fk \ar[d]  & &  Gk \ar[d]    &  & Fl \ar[d] & & Gl \ar[d]      \\
& FA \ar@{-}[d]  & &  GA\ar@{-}[drr]    &  & FB \ar@{-}[dll] & & GB \ar@{-}[d]      \\
& FA \ar@/_/@{-}[dr]  & &  FB\ar@/^/@{-}[dl]    &  & GA \ar@/_/@{-}[dr] & & GB \ar@/^/@{-}[dl]      \\
= &&  Ff\ar[d] && && Gf\ar[d]  && = \\
&& FC \ar@/_/@{-}[drr] &&              & & GC \ar@/^/@{-}[dll] \\
&&                         &&  u_C\ar[d]                           &&                           \\
&&                        &&  HC                                       &&
}
\xymatrix@R=2pc@C=0.2pc{
FR \ar@{-}[d] & GR \ar@{-}[dr]    &   FS \ar@{-}[dl] & GS \ar@{-}[d]    & FT \ar@{-}[d] & GT \ar@{-}[dr]    
& FU \ar@{-}[dl] & GU \ar@{-}[d]                     \\
FR \ar@{-}[d] & FS \ar@{-}[d]    &   GR \ar@{-}[drr] & GS \ar@{-}[drr]    & FT \ar@{-}[dll] & FU \ar@{-}[dll]    
& GT \ar@{-}[d] & GU \ar@{-}[d]                     \\
FR \ar@/_/@{-}[dr] & FS \ar@{-}[d]    &   FT \ar@/_/@{-}[dr] & FU \ar@{-}[d]    & GR \ar@/_/@{-}[dr] & GS \ar@{-}[d]    
& GT \ar@/_/@{-}[dr] & GU \ar@{-}[d]                     \\
& Fk \ar[d]  & &  Fl \ar[d]    &  & Gk \ar[d] & & Gl \ar[d]      \\
& FA \ar@/_/@{-}[dr]  & &  FB\ar@/^/@{-}[dl]    &  & GA \ar@/_/@{-}[dr] & & GB \ar@/^/@{-}[dl]      \\
&&  Ff\ar[d] && && Gf\ar[d]  && = \\
&& FC \ar@/_/@{-}[drr] &&              & & GC \ar@/^/@{-}[dll] \\
&&                         &&  u_C\ar[d]                           &&                           \\
&&                        &&  HC                                       &&
}
\]
\[
\xymatrix@R=2pc@C=0.2pc{
FR \ar@{-}[d] & GR \ar@{-}[drrr]    &   FS \ar@{-}[dl] & GS \ar@{-}[drr]    & FT \ar@{-}[dll] & GT \ar@{-}[dr]    
& FU \ar@{-}[dlll] & GU \ar@{-}[d]                     \\
FR \ar@/_/@{-}[dr] & FS \ar@{-}[d]    &   FT \ar@/_/@{-}[dr] & FU \ar@{-}[d]    & GR \ar@/_/@{-}[dr] & GS \ar@{-}[d]    
& GT \ar@/_/@{-}[dr] & GU \ar@{-}[d]                     \\
& Fk \ar[d]  & &  Fl \ar[d]    &  & Gk \ar[d] & & Gl \ar[d]      \\
& FA \ar@/_/@{-}[dr]  & &  FB\ar@/^/@{-}[dl]    &  & GA \ar@/_/@{-}[dr] & & GB \ar@/^/@{-}[dl]      \\
= &&  Ff\ar[d] && && Gf\ar[d]  && = \\
&& FC \ar@/_/@{-}[drr] &&              & & GC \ar@/^/@{-}[dll] \\
&&                         &&  u_C\ar[d]                           &&                           \\
&&                        &&  HC                                       &&
}
\xymatrix@R=2pc@C=0.2pc{
\\
FR \ar@/_/@{-}[dr] & GR \ar@{-}[d]    &   FS \ar@/_/@{-}[dr] & GS \ar@{-}[d]    & FT \ar@/_/@{-}[dr] & GT \ar@{-}[d]    
& FU \ar@/_/@{-}[dr] & GU \ar@{-}[d]                     \\
& u_R \ar[d]  & &  u_S \ar[d]    &  & u_T \ar[d] & & u_U \ar[d]      \\
& HR \ar@/_/@{-}[drrr]  & &  HS\ar@/_/@{-}[dr]    &  & HT \ar@/^/@{-}[dl] & & HU \ar@/^/@{-}[dlll]      \\
&&                         &&  Hf(k,l) \ar[d]                           &&                           \\
&&                        &&  HC                                       &&
}
\]
Now, it is clear that if $t:F,G\to L$ factors in $[\M,\N]$ as $t = ul$ then $l$ is uniquely defined by $t_A = u_Al_A$.
Then, to show that $u: F,G\to H$ is universal we need to check that $l_A:HA\to LA$ so defined
is indeed a map in $[\M,\N]$, that is (for 2-ary arrows) that $Lf(l_A,l_B)=l_CHf$, for any $f:A,B\to C$:
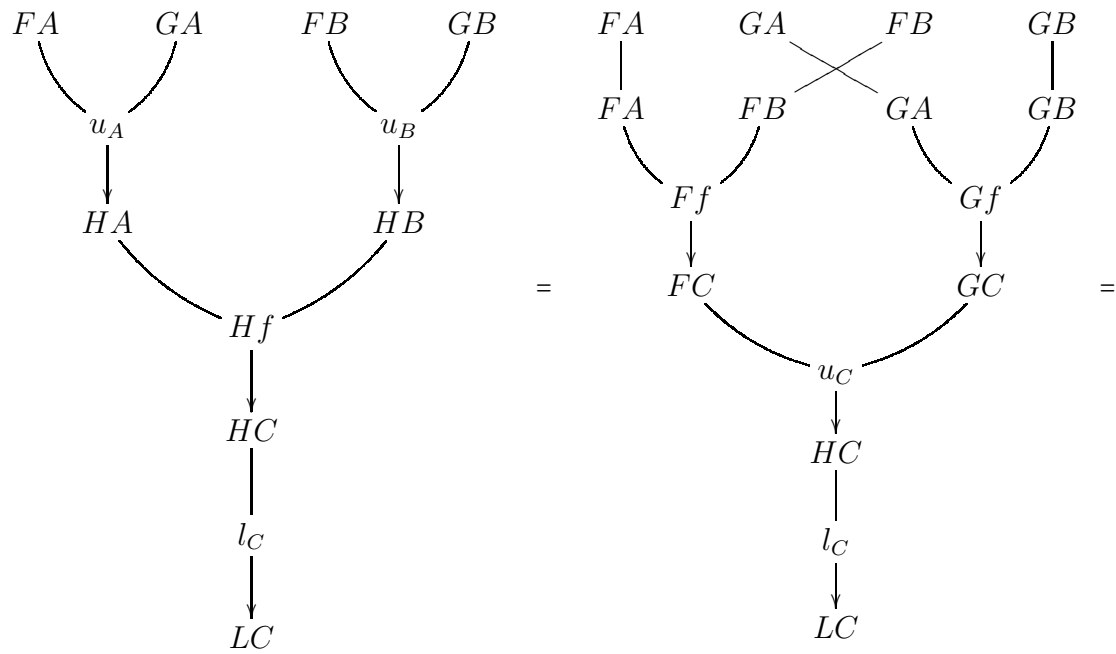
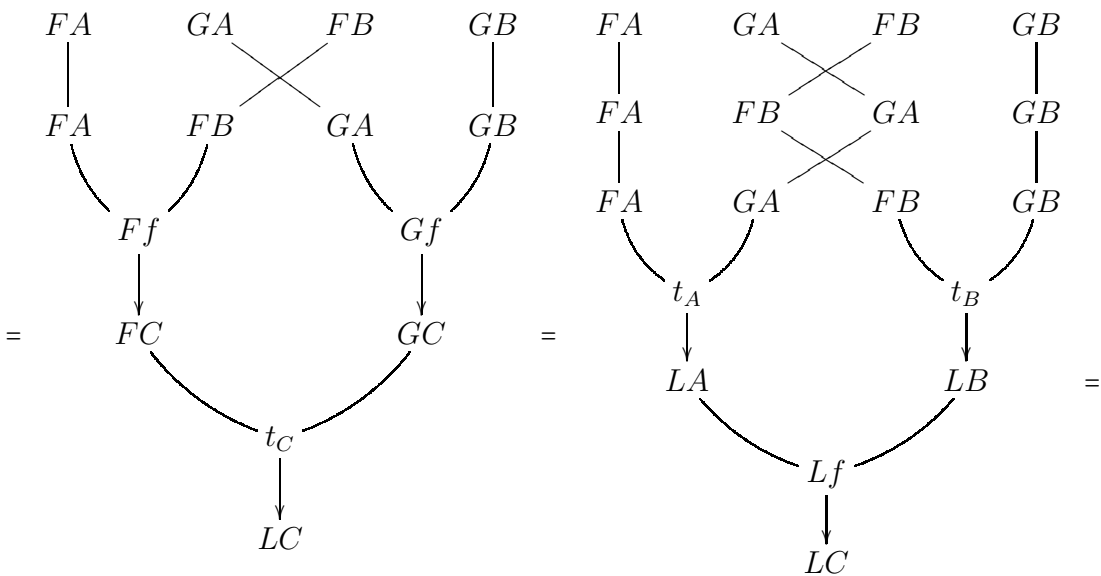
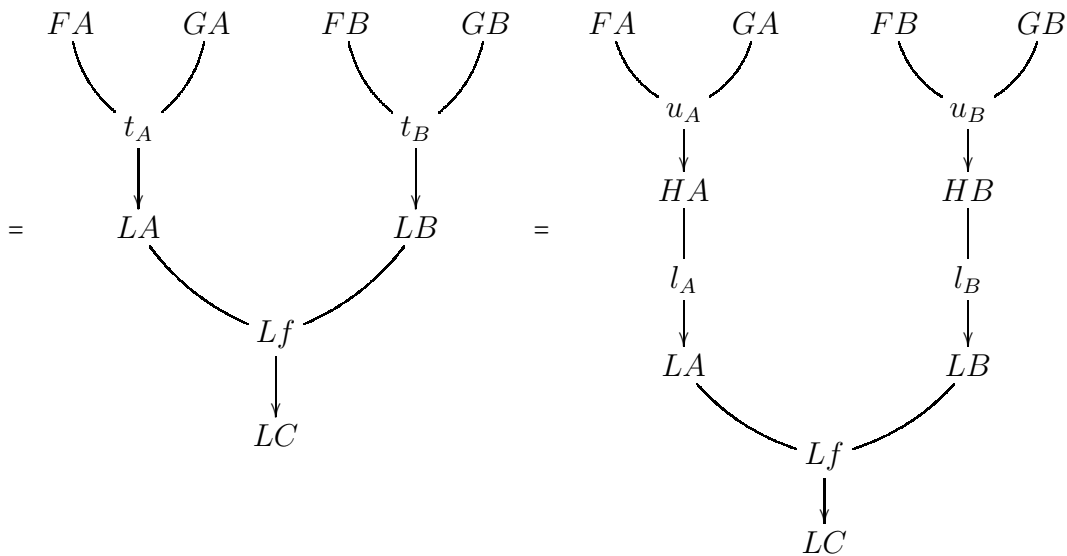
\[
\xymatrix@R=2pc@C=0.2pc{
& FA \ar@/_/@{-}[dr]  & &  GA \ar@/^/@{-}[dl]    &  & FB \ar@/_/@{-}[dr] & & GB \ar@/^/@{-}[dl]      \\
&&  u_A \ar[d] && && u_B \ar[d]  \\
&& HA \ar@/_/@{-}[drr] &&              & & HB \ar@/^/@{-}[dll] &&  \\
&&                         &&  Hf \ar[d]                           &&                           \\
&&                        &&  HC \ar@{-}[d]                                    &&            \\
&&                        &&  l_C \ar[d]                                     &&            \\
&&                        &&  LC                                       &&            \\
}
\xymatrix@R=1.5pc@C=0.2pc{
&& FA \ar@{-}[d]  & &  GA \ar@{-}[drr]    &  & FB \ar@{-}[dll] & & GB \ar@{-}[d]      \\
&& FA \ar@/_/@{-}[dr]  & &  FB \ar@/^/@{-}[dl]    &  & GA \ar@/_/@{-}[dr] & & GB \ar@/^/@{-}[dl]      \\
&&&  Ff \ar[d] && && Gf \ar[d]  \\
= &&& FC \ar@/_/@{-}[drr] &&              & & GC \ar@/^/@{-}[dll] && = \\
&&&                         &&  u_C \ar[d]                           &&                           \\
&&&                        &&  HC \ar@{-}[d]                                    &&            \\
&&&                        &&  l_C \ar[d]                                     &&            \\
&&&                        &&  LC                                       &&            \\
}
\]
\[
\xymatrix@R=2pc@C=0.2pc{
& FA \ar@{-}[d]  & &  GA \ar@{-}[drr]    &  & FB \ar@{-}[dll] & & GB \ar@{-}[d]      \\
& FA \ar@/_/@{-}[dr]  & &  FB \ar@/^/@{-}[dl]    &  & GA \ar@/_/@{-}[dr] & & GB \ar@/^/@{-}[dl]      \\
&&  Ff \ar[d] && && Gf \ar[d]  \\
= && FC \ar@/_/@{-}[drr] &&              & & GC \ar@/^/@{-}[dll] && = \\
&&                         &&  t_C \ar[d]                           &&                           \\
&&                        &&  LC                                     &&         \\   
}
\xymatrix@R=1.6pc@C=0.2pc{
& FA \ar@{-}[d]  & &  GA \ar@{-}[drr]    &  & FB \ar@{-}[dll] & & GB \ar@{-}[d]      \\
& FA \ar@{-}[d]  & &  FB \ar@{-}[drr]    &  & GA \ar@{-}[dll] & & GB \ar@{-}[d]      \\
& FA \ar@/_/@{-}[dr]  & &  GA \ar@/^/@{-}[dl]    &  & FB \ar@/_/@{-}[dr] & & GB \ar@/^/@{-}[dl]      \\
&&  t_A \ar[d] && && t_B \ar[d]  \\
&& LA \ar@/_/@{-}[drr] &&              & & LB \ar@/^/@{-}[dll] && = \\
&&                         &&  Lf \ar[d]                           &&                           \\
&&                        &&  LC                                     &&         \\   
}
\]
\[
\xymatrix@R=2pc@C=0.2pc{
& FA \ar@/_/@{-}[dr]  & &  GA \ar@/^/@{-}[dl]    &  & FB \ar@/_/@{-}[dr] & & GB \ar@/^/@{-}[dl]      \\
&&  t_A \ar[d] && && t_B \ar[d]  \\
= && LA \ar@/_/@{-}[drr] &&              & & LB \ar@/^/@{-}[dll] && = \\
&&                         &&  Lf \ar[d]                           &&                           \\
&&                        &&  LC                                     &&         \\   
}
\xymatrix@R=1.5pc@C=0.2pc{
& FA \ar@/_/@{-}[dr]  & &  GA \ar@/^/@{-}[dl]    &  & FB \ar@/_/@{-}[dr] & & GB \ar@/^/@{-}[dl]      \\
&&  u_A \ar[d] && && u_B \ar[d]  \\
&& HA \ar@{-}[d] &&              & & HB \ar@{-}[d]  \\
&& l_A \ar[d] &&              & & l_B \ar[d]  \\
&& LA \ar@/_/@{-}[drr] &&              & & LB \ar@/^/@{-}[dll]  \\
&&                         &&  Lf \ar[d]                           &&                           \\
&&                        &&  LC                                     &&         \\   
}
\]
\epf

\begin{remark}
\label{monten}
In particular, for $\M = 1\t$ one gets the usual monoidal structure on $\cMon(\N) = [1\t,\N]_-$.
More specifically, let $M,N:1\t \to \N$ be monoids on the objects $A$ and $B$ respectively,
and let $u:A,B\to A\otm B$ be universal in $\N$. 
Then $u$ is universal in $[1\t,\N]$ as well, when $A\otm B$ is given the monoid structure such that
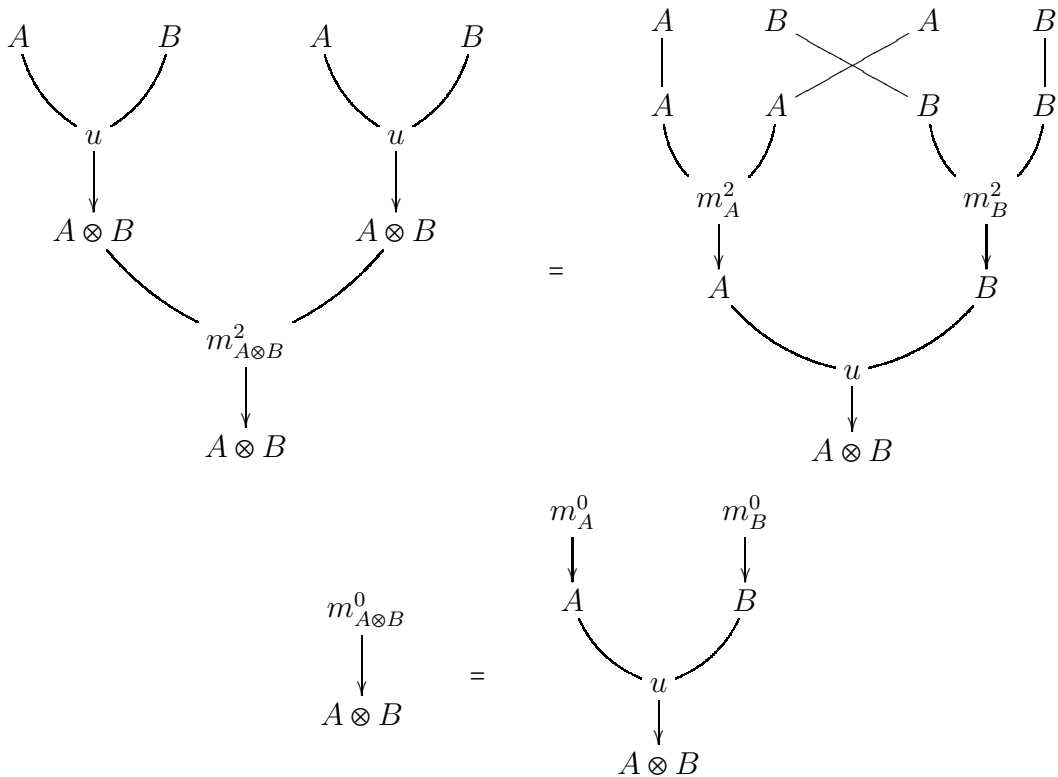
\[
\xymatrix@R=2pc@C=0.2pc{
\\
& A \ar@/_/@{-}[dr]  & &  B\ar@/^/@{-}[dl]    &  & A \ar@/_/@{-}[dr] & & B \ar@/^/@{-}[dl]      \\
&&  u\ar[d] && && u\ar[d]  \\
&& A\otm B \ar@/_/@{-}[drr] &&              & & A\otm B \ar@/^/@{-}[dll] \\
&&                         &&  m^2_{A\otm B} \ar[d]                           &&                           \\
&&                        &&  A\otm B                                       &&
}
\xymatrix@R=2pc@C=1pc{
\\
\\
\\
              &  &                       \\
              & = & 
}
\xymatrix@R=1.5pc@C=0.2pc{
\\
& A \ar@{-}[d] & &   B \ar@{-}[drr] & & A \ar@{-}[dll] & & B \ar@{-}[d]    \\
& A \ar@/_/@{-}[dr]  & &  A\ar@/^/@{-}[dl]    &  & B \ar@/_/@{-}[dr] & & B \ar@/^/@{-}[dl]      \\
&&  m^2_A\ar[d] && && m^2_B\ar[d]   \\
&& A \ar@/_/@{-}[drr] &&              & & B \ar@/^/@{-}[dll] \\
&&                         &&  u \ar[d]                           &&                           \\
&&                        &&  A\otm B                                       &&
}
\]
\[
\xymatrix@R=2pc@C=0.2pc{
\\
&& \,\, m^0_{A\otm B} \ar[d]       \\
&&  A\otm B                                     
}
\xymatrix@R=2pc@C=1pc{
\\
              &  &                       \\
              & = & 
}
\xymatrix@R=1.5pc@C=0.2pc{
m^0_A\ar[d] && m^0_B\ar[d]   \\
A \ar@/_/@{-}[dr] && B \ar@/^/@{-}[dl] \\
&  u \ar[d]                                 \\
&  A\otm B                                 
}
\]
Similarly, a 0-ary $u:\,\,\to \I$ universal in $\N$ is universal also in $[1\t,\N]$, 
when $\I$ is given the monoid structure such that $m^2_\I(u,u) = u = m^0_\I$.
\end{remark}

Now, from Corollary \ref{mon2} and propositions \ref{repseq} and \ref{rep} we get

\begin{corol}
\label{corefsum}
The coreflection of $\sMlt$ in $\Seq$ restricts to $\sRep \subset \sMlt$.
Thus the commutative monoids construction $\cMon(-) = [1\t, -]_-$
gives the coreflection of symmetric monoidal categories in the cocartesian ones.
\epf
\end{corol}
(In \cite{fox} it is stated a similir result, but with strict in place of lax monoidal functors.)
In particular,

\begin{corol}
If $\M$ is representable, then the category $\cMon(\M)$ of commutative monoids
in $\M$, with the tensor product of Remark \ref{monten}, is cocartesian.
\epf
\end{corol}


So, cocartesian symmetric monoidal categories $(\C,\I,\otm)$ can be characterized 
as those of the form $\C\iso\cMon(\D)$ for a symmetric monoidal $(\D,\I,\otm)$,
or also as those for which $\cMon(\C)\to\C$ is an iso.
In fact, propositions \ref{char} and \ref{repseq} yield a more effective characterization, 
namely as those monoidal categories with a central monoid (or unital magma).
The latter amounts to a monoid (or unital magma) $m_A$ on each object $A\in\C$,
which commutes with both unary and universal arrows
(since these generate the corresponding multicategory $\C_\otm$).
The first condition says that each arrow is a monoid (or unital magma) morphism, 
that is that the $m^0_A$ and the $m^2_A$ are natural. 
The second one says that the monoid (unital magma) structure $m_{A\otm B}$ on $A\otm B$ 
is the expected one, related as in Remark \ref{monten} to $m_A$ and $m_B$ 
(and similarly for $m_\I$).
Thus, we find again the following piece of folklore:

\begin{prop}
A symmetric monoidal category $(\C,\I,\otm)$ is cocartesian 
iff it has a natural monoid (or unital magma) $m_A$ ($A\in\C$)
such that $m_{A\otm B}$ is related in the usual way with $m_A$ and $m_B$.
\epf
\end{prop}
(In fact, the weaker conditions $m_{A\otm B}(u_{A,B}(A,m^0_B),u_{A,B}(m^0_A,B)) = u_{A,B}$ 
are sufficient to ensure that $\C$ is cocartesian, 
as pointed out by Jeff Egger on the categories mailing list, 
and as it follows from inspection of the proof of Proposition \ref{char}.)


\section{Cartesian multicategories and preadditive categories}
\label{fp}

The algebra of abstract operations encoded in a (small) multicategory $\M$ can be seen as
an algebraic theory, whose category of models is $\Mlt(\M,\Set_\tm)$.
The algebraic theories which can be so represented are those which can be
given by equations between terms with the same variables in the same order
(named \"strongly regular" in \cite{leinster}).
Notably, the terminal multicategory is the theory for monoids.
In order to allow permutation (exchange) of variables (as for commutative monoids)
one considers symmetric multicategories instead.

More general algebraic theories, such as the theory of groups,
require laws which can be expressed in multicategories of the form $\C_\tm$
(for a finite product category $\C$) where one can furthermore duplicate or delete 
the variables.
For instance, given an operation $f:A\tm B \tm A \to D$ in $\Set$, one gets another
operation $\si f: B\tm C \tm A \to D$ by
\[ \si f(b,c,a) = f(a,b,a) \]
It is convenient to see the domain of $f$ as the family of sets ${\bf 3} = \{1,2,3\}\to\obj\Set$
given by $1\mapsto A,2\mapsto B,3\mapsto A$, that of $\si f$ as the family of sets ${\bf 3} \to\obj\Set$ 
given by $1\mapsto B,2\mapsto C,3\mapsto A$ and 
$\si:1\mapsto 3,2\mapsto 1,3\mapsto 3$ as a family morphism 
(that is a mapping ${\bf 3} \to \bf 3$ over $\obj\Set$).

More generally, if $A_1,\cd , A_n$ and $B_1,\cd , B_m$ are families of objects in $\C$
and $\si:{\bf n} \to \bf m$ is a map of families in the above sense, one gets (for each $C\in\C$) a mapping
\[ \si(-) : \C(A_1\tm\cd\tm A_n,C) \to \C(B_1\tm\cd\tm B_m,C) \]
by precomposing with the map $\ov\si : B_1\tm\cd\tm B_m \to A_1\tm\cd\tm A_n$
such that $p_i\ov\si = q_{\si i}$.
As particular instances of $\ov\si$ we find projections and diagonals.

Abstracting from this, one gets the notion of cartesian multicategory defined below
and the 2-category $\fpMlt$ of cartesian multicategories.
 
Thus the categories $\Mlt$, $\sMlt$ and $\fpMlt$ can be seen as the
\"doctrines" of strongly regular, linear and algebraic theories
respectively, and the latter has the same expressive power of (many-sorted) Lawvere theories.
Since the obvious forgetful functors have left adjoints (see \cite{gould})
each doctrine \"contains" the weaker ones.
For instance, one has an algebraic theory for commutative monoids by taking the
free algebraic theory on $1\t\in\sMlt$, which is $\bf N\t\in\fpMlt$ 
(the sequential multicategory on the monoid of natural numbers
with the cartesian structure given by zero and addition).


\subsection{Cartesian multicategories}

Cartesian multicategories have appeared in various guises in the literature;
when considered in relation with deduction theory  they are usually called \"Gentzen multicategories"
(see for instance \cite{lambek}).
The definition that we give here follows essentially \cite{boardman} and
\cite{fiore2} (see also \cite{gould}, especially for the case of operads).

Let $\bf N$ be the full subcategory of $\Set$ which has objects ${\bf 0, 1, 2}, \cdots $, with
${\bf n} = \{ 1, 2, \cdots , n \}$ (so that ${\bf 0} = \emptyset$) and consider, 
for a multicategory $\M$, the obvious comma category ${\bf N}/\obj\M$. 
Now, the domain of an $n$-ary arrow in $\M$ is in fact a mapping $\alpha:{\bf n} \to \obj\M$,
that is an object of ${\bf N}/\obj\M$; thus, for any fixed codomain object $A\in\M$
we have a mapping
\[
\obj ({\bf N}/\obj\M) \to \obj\Set  \qv  \alpha \mapsto \M(\alpha;A) = \M(A_1,\cd,A_n; A)
\]
To give a cartesian structure on $\M$ means to extend these mappings to functors
\[
{\bf N}/\obj\M \to \Set
\] 
in a way that is compatible with composition.

To illustrate such functoriality and compatibility conditions 
and to do calculations in cartesian multicategories, it is convenient to make use of 
a graphical calculus (which in fact we have already used in the particular case of symmetric multicategories).
For instance, if $f:X,Y,X\to U$ is an arrow in $\M$ and $\si$ is the following mapping in ${\bf N}/\obj\M$
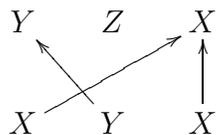
\[
\xymatrix@R=2pc@C=1.5pc{
Y            & Z  &   X                     \\
X \ar[urr] &        Y   \ar[ul]          & X \ar[u] 
}
\]
we have the arrow $\si f:Y,Z,X\to U$ which we denote by
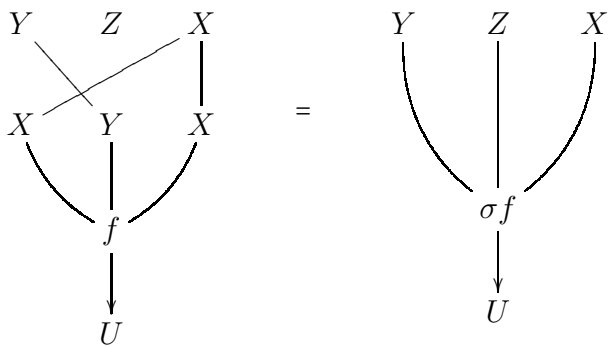
\[
\xymatrix@R=2pc@C=1.5pc{
Y  \ar@{-}[dr]          & Z  &   X  \ar@{-}[d]  \ar@{-}[dll]                   \\
X \ar@/_/@{-}[dr] &        Y   \ar@{-}[d]          & X \ar@/^/@{-}[dl] \\
                         &  f \ar[d]                           &                           \\
                         &  U                                       &
}
\xymatrix@R=2pc@C=1.5pc{
              &  &                       \\
              & = & 
}
\xymatrix@R=2pc@C=1.5pc{
Y  \ar@/_1pc/@{-}[ddr]     &   Z \ar@{-}[dd]   &   X \ar@/^1pc/@{-}[ddl]                     \\
                                                 &                                &                           \\
                                                 &  \si f\ar[d]       &                            \\
                                                 &  U                            &
}
\]

Functoriality (that is the fact that mappings over $\obj\M$ act on hom-sets) 
\[ 
\rho(\si f) = (\rho\si) f
\]
is illustrated by the fact that we can compose mappings over arrows. For instance,
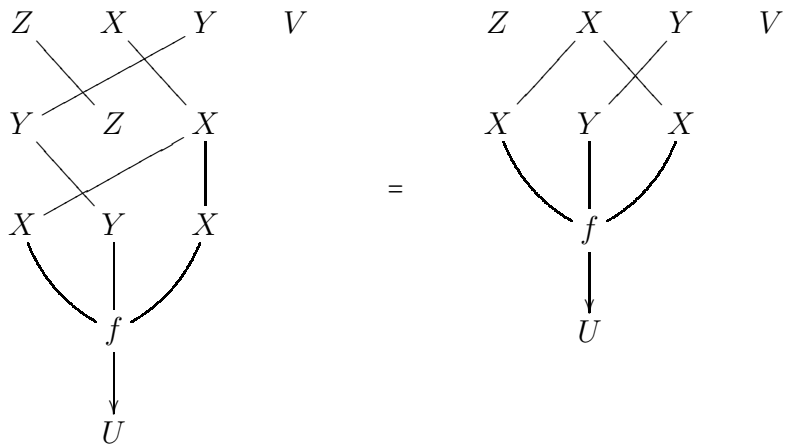
\[
\xymatrix@R=2pc@C=1.5pc{
Z  \ar@{-}[dr]          & X \ar@{-}[dr]  &   Y \ar@{-}[dll]   & V               \\
Y  \ar@{-}[dr]          & Z  &   X  \ar@{-}[d]  \ar@{-}[dll]                   \\
X \ar@/_/@{-}[dr] &        Y   \ar@{-}[d]          & X \ar@/^/@{-}[dl] \\
                         &  f \ar[d]                           &                           \\
                         &  U                                       &
}
\xymatrix@R=2pc@C=1.5pc{
              &  &                       \\
              &&  \\
              & = & 
}
\xymatrix@R=2pc@C=1.5pc{
Z           & X\ar@{-}[dl] \ar@{-}[dr]    &   Y  \ar@{-}[dl]    & V            \\
X \ar@/_/@{-}[dr] &        Y   \ar@{-}[d]          & X \ar@/^/@{-}[dl] \\
                         &  f \ar[d]                           &                           \\
                         &  U                                       &
}
\]

As for the compatibility conditions with respect to composition, there are two of them.
The first one is pretty obvious: 
\[ 
f(\si_1 f_1,\cdots, \si_n f_n) = (\si_1+\cdots+ \si_n)f(f_1,\cdots, f_n)
\]
that is, composing $f$ with arrows $f_i$ acted upon by $\si_i$ is the same as 
composing $f$ with the $f_i$ and then acting on it with the obvious \"sum" 
of maps in ${\bf N}/\obj\M$.
Thus, diagrams such as this one
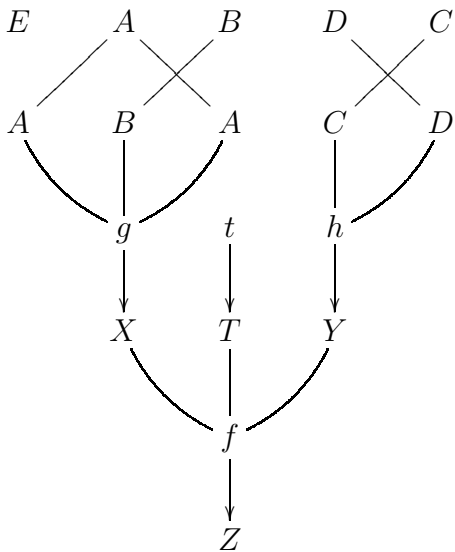
\[
\xymatrix@R=2pc@C=2pc{
E                & A \ar@{-}[dl]\ar@{-}[dr] &   B  \ar@{-}[dl]   & D\ar@{-}[dr]  & C\ar@{-}[dl]           \\
A  \ar@/_/@{-}[dr]                & B \ar@{-}[d] &   A  \ar@/^/@{-}[dl]   & C\ar@{-}[d]  & D\ar@/^/@{-}[dl]           \\
                 & g \ar[d]              &       t \ar[d]                   & h \ar[d]                            \\
                 & X \ar@{-}@/_/[dr]    &          T   \ar@{-}[d]           & Y \ar@{-}@/^/[dl]                            \\
                 &                               &  f \ar[d]                              &      &                    \\
                 &        &  Z                                          &      &
}
\]
can be interpreted in an unambiguous way in a cartesian multicategory.

The second compatibility condition concerns composition in the case when it is $f$ 
that is acted upon by a mapping $\si$:
\[ 
(\si f)(f_1,\cdots, f_n) = \si'(f(f_{\si 1},\cdots, f_{\si n}))
\]
where $\si'$ is a suitably defined map in ${\bf N}/\obj\M$, which is graphically obvious.
For instance,
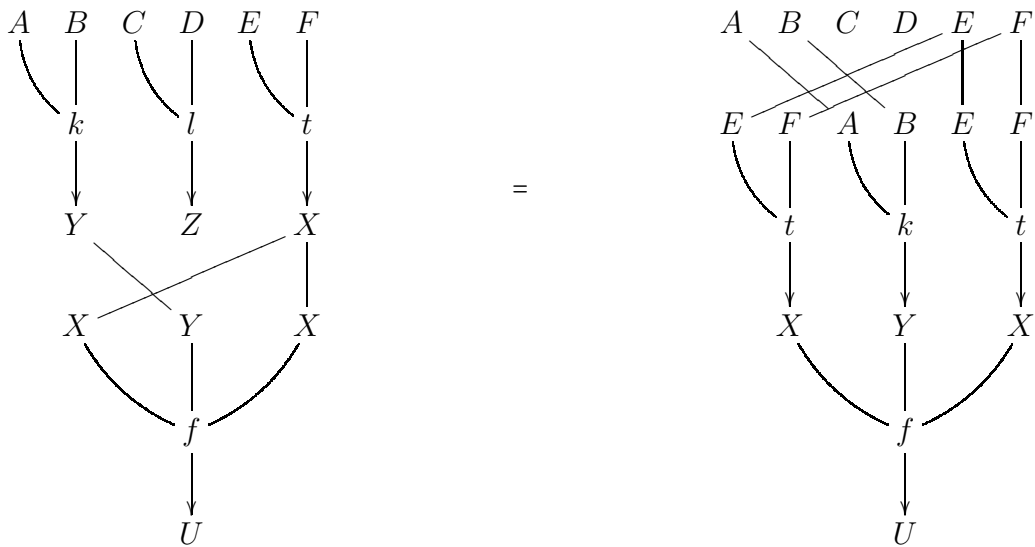
\[
\xymatrix@R=2pc@C=0.5pc{
A \ar@/_/@{-}[dr]  &   B\ar@{-}[d]   & C\ar@/_/@{-}[dr] & D\ar@{-}[d]  & E \ar@/_/@{-}[dr] &  F\ar@{-}[d]      \\
&  k\ar[d] && l \ar[d]&& t\ar[d]  \\
& Y   \ar@{-}[drr]         &&  Z  &&   X  \ar@{-}[d]  \ar@{-}[dllll]                   \\
& X \ar@/_/@{-}[drr] &&        Y   \ar@{-}[d]         & & X \ar@/^/@{-}[dll] \\
&                         &&  f \ar[d]                           &&                           \\
&                        &&  U                                       &&
}
\xymatrix@R=2pc@C=5pc{
              &  &                       \\
              &  &                       \\
              & = & 
}
\xymatrix@R=2pc@C=0.5pc{
A \ar@{-}[drr]  &   B\ar@{-}[drr]   &  C & D &  E\ar@{-}[d] \ar@{-}[dllll] &  F\ar@{-}[d] \ar@{-}[dllll]    \\
E \ar@/_/@{-}[dr]  &   F\ar@{-}[d]   & A\ar@/_/@{-}[dr]  &  B\ar@{-}[d]  & E\ar@/_/@{-}[dr] &  F\ar@{-}[d]  &      \\
&  t\ar[d] && k \ar[d]&& t\ar[d]  \\
& X \ar@/_/@{-}[drr] &&        Y   \ar@{-}[d]         & & X \ar@/^/@{-}[dll] \\
&                         &&  f \ar[d]                           &&                           \\
&                        &&  U                                       &&
}
\]


A functor $F:\M\to\M'$ induces an obvious functor  ${\bf N}/\obj\M \to {\bf N}/\obj\M'$,
and maps of cartesian multicategories are of course those functors which 
commute with the actions of ${\bf N}/\obj\M$ and ${\bf N}/\obj\M'$.
For instance, for $\si$ as above,

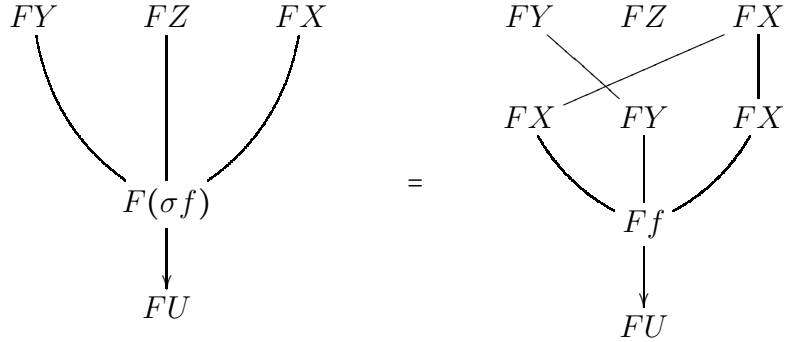
\[
\xymatrix@R=2pc@C=1.5pc{
FY  \ar@/_1pc/@{-}[ddr]     &   FZ \ar@{-}[dd]   &   FX \ar@/^1pc/@{-}[ddl]                     \\
                                                 &                                &                           \\
                                                 &  F(\si f)\ar[d]       &                            \\
                                                 &  FU                            &
}
\xymatrix@R=2pc@C=1.5pc{
\\
              &  &                       \\
              & = & 
}
\xymatrix@R=2pc@C=1.5pc{
FY  \ar@{-}[dr]          & FZ  &   FX  \ar@{-}[d]  \ar@{-}[dll]                   \\
FX \ar@/_/@{-}[dr] &        FY   \ar@{-}[d]          & FX \ar@/^/@{-}[dl] \\
                         &  Ff \ar[d]                           &                           \\
                         &  FU                                       &
}
\]

By defining 2-cells as in $\sMlt$, we so obtain the 2-category $\fpMlt$ of cartesian multicategories, 
with the obvious forgetful functor $U:\fpMlt \to \sMlt$.


\begin{remark}
\label{cart0}

An important special case are the actions of constant shape $\bf n \to \bf 1$, 
giving \"contraction" mappings on hom-sets:
\[
\gamma_n:\M(A,\cd,A;B)\to \M(A;B)
\]
If $\M$ is promonoidal, its cartesian structure is determined by its symmetric structure
and by contractions.
Indeed, any $\si: \bf m \to \bf n$ can be written as a bijection followed by a monotone $\si': \bf m \to \bf n$,
and the latter is a sum of constant mappings $\si'_1,\cd,\si'_n$.
Now, to know $\si' f$ it is enough to write $f$ as the composite $f'(f_1,\cd,f_n)$ following
the pattern of the $\si'_n$ (using promonoidality) and then to contract the $f_i$ 
(using the first compatibility condition).
For instance,
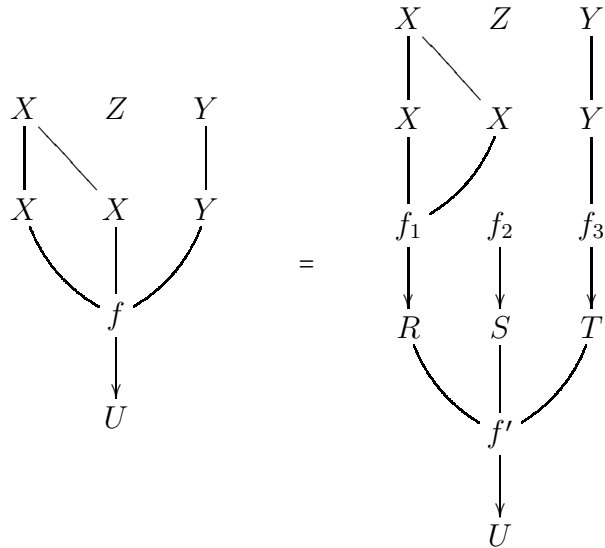
\[
\xymatrix@R=2pc@C=1.5pc{
\\
X  \ar@{-}[d]\ar@{-}[dr]          & Z  &   Y  \ar@{-}[d]                     \\
X \ar@/_/@{-}[dr] &        X   \ar@{-}[d]          & Y \ar@/^/@{-}[dl] \\
                         &  f \ar[d]                           &                           \\
                         &  U                                       &
}
\xymatrix@R=2pc@C=1.5pc{
\\
\\
              &  &                       \\
              & = & 
}
\xymatrix@R=2pc@C=1.5pc{
X  \ar@{-}[d]\ar@{-}[dr]          & Z  &   Y  \ar@{-}[d]                     \\
X  \ar@{-}[d]          & X\ar@/^/@{-}[dl]  &   Y  \ar@{-}[d]                     \\
f_1 \ar[d] &     f_2 \ar[d]             & f_3 \ar[d]                \\
R \ar@/_/@{-}[dr] &     S   \ar@{-}[d]      & T \ar@/^/@{-}[dl]                      \\
                         &  f' \ar[d]                           &                           \\
                         &  U                                       &
}
\]
\end{remark}


\subsection{Algebraic products}

Let $\M$ be a cartesian multicategory.
An algebraic product of $A,B\in\M$ consists of an object $C$ along with maps
$p:C\to A$, $q:C\to B$ and $u:A,B\to C$ such that
\eq
\label{algp1}
\xymatrix@R=2pc@C=0.5pc{
                                                   &                                         &&             C \ar@{-}[dll] \ar@{-}[drr]        &                                \\
                                                  &  C  \ar@{-}[d]               &&                                     && C\ar@{-}[d]              &          \\
                                                  & p \ar[d]                     &&                                       & &q \ar[d]                     &     \\
                                                  & A \ar@{-}@/_/[drr]   &&                                       &  & B \ar@{-}@/^/[dll]  &                         \\
                                                  &                                      &&  u \ar[d]                 &                                     &                    \\
                                                  &                                     & &  C                                    &                                     &
}
\xymatrix@R=2.2pc@C=1.5pc{ 
&\\            
              &  &         \\
              &  &           \\
              & = & 
}
\xymatrix@R=2pc@C=0.5pc{
& \\
                                                   &                                                \\
                                                  &  C  \ar@{-}[d]                        \\
                                                  & \id \ar[d]                                  \\
                                                  & C                       
}
\eeq

\eq
\label{algp2}
\xymatrix@R=2pc@C=0.2pc{
                                                  & A \ar@{-}@/_/[drr]   &&                                       &  & B \ar@{-}@/^/[dll]                           \\
                                                  &                                      &&  u \ar[d]                 &                                                         \\
                                                  &                                     & &  C    \ar@{-}[d]                                                               \\
                                                  &                                     &&    p \ar[d]                                            \\
                                                  &                                    &&    A                                     
}
\xymatrix@R=2.2pc@C=1pc{             
              &  &         \\
              &  &           \\
              & = & 
}
\xymatrix@R=2pc@C=0.5pc{
                                                  &  A  \ar@{-}[d]  &&  B                                             \\
                                                  &  A  \ar@{-}[d]                        \\
                                                  & \id \ar[d]                                  \\
                                                  & A                       
}
\qq \qq
\xymatrix@R=2pc@C=0.2pc{
                                                  & A \ar@{-}@/_/[drr]   &&                                       &  & B \ar@{-}@/^/[dll]                          \\
                                                  &                                      &&  u \ar[d]                 &                                                       \\
                                                  &                                     & &  C    \ar@{-}[d]                                                               \\
                                                  &                                     &&    q \ar[d]                                            \\
                                                  &                                    &&    B                                     
}
\xymatrix@R=2.2pc@C=1pc{ 
              &  &         \\
              &  &           \\
              & = & 
}
\xymatrix@R=2pc@C=0.5pc{
                                                  A   &&  B \ar@{-}[d]                    \\
                                                  &&  B  \ar@{-}[d]                        \\
                                                  && \id \ar[d]                                  \\
                                                  && B                       
}
\eeq

The definition easily extends to any finite family $A_1,\cd,A_n$ of objects.
In particular, for $n=0$, an algebraic product of the empty family 
is an object $C$ with a 0-ary arrow $u\in(\,\,;C)$ such that
\[
\xymatrix@R=2pc@C=0.5pc{
                                                  &                                       &&     C                                                         \\
                                                  &                                      &&  u \ar[d]                 &                                                         \\
                                                  &                                     & &  C                                                                                                                    
}
\xymatrix@R=2.2pc@C=1pc{  
               &  &           \\
              & = & 
}
\xymatrix@R=2pc@C=0.5pc{                                                 
                                                  &  C  \ar@{-}[d]                        \\
                                                  & \id \ar[d]                                  \\
                                                  & C                       
}
\]
(On the left, the 0-ary $u$ is acted upon by a map of shape $\bf 0 \to 1$.)

We say that the cartesian multicategory $\M$ has algebraic products 
if any finite family of objects in $\M$ has an algebraic product.


\subsection{Universal products}

Let $\M$ be a cartesian multicategory.
A universal product of $A,B\in\M$ consists of an object $C$ along with maps
$p:C\to A$ and $q:C\to B$ such that any pair of arrows $f:X_1,\cd,X_n\to A$
and $g:X_1,\cd,X_n\to B$ with the same domain factors uniquely as $f=pt$ and $g=qt$.

The definition of course extends to any family of objects.
In particular, a universal product of the empty family is an object $C$ such that
each hom-set $\M(X_1,\cd,X_n, C)$ has a unique element.

\begin{prop}
\label{cart1}
For any finite family $A_1,\cd,A_n$ in a cartesian multicategory $\M$, the following are equivalent:
\begin{enumerate}
\item
$A_1,\cd,A_n$ has an algebraic product.
\item
$A_1,\cd,A_n$ has a universal product.
\item
$A_1,\cd,A_n$ is the domain of a preuniversal arrow.
\end{enumerate}
\end{prop}
\pf
We prove the case $n=2$, the other ones being similar.
To prove that (1) implies (3),
we show that in an algebraic product $p:C\to A$, $q:C\to B$ and $u:A,B\to C$, 
$u$ is preuniversal;
that is, that any arrow $f:A,B \to D$ factors uniquely as $f=tu$.
Indeed, if such a $t:C\to D$ does exist, we have 
\[
\xymatrix@R=1.5pc@C=0.5pc{
                                                  \\          \\
                                                  &                                     & &  C  \ar@{-}[d]                                   &                                     &           \\
                                                  &                                     & &  t    \ar[d]                                &                                     &     \\
                                                  &                                     & &  D                                   &                                     &
}
\xymatrix@R=2.2pc@C=1.5pc{ 
&\\            
              &  &         \\
              &  &           \\
              & = & 
}
\xymatrix@R=1.5pc@C=0.5pc{
                                                   &                                         &&             C \ar@{-}[dll] \ar@{-}[drr]        &                                \\
                                                  &  C  \ar@{-}[d]               &&                                     && C\ar@{-}[d]              &          \\
                                                  & p \ar[d]                     &&                                       & &q \ar[d]                     &     \\
                                                  & A \ar@{-}@/_/[drr]   &&                                       &  & B \ar@{-}@/^/[dll]  &                         \\
                                                  &                                      &&  u \ar[d]                 &                                     &                    \\
                                                  &                                     & &  C  \ar@{-}[d]                                   &                                     &           \\
                                                  &                                     & &  t    \ar[d]                                &                                     &     \\
                                                  &                                     & &  D                                   &                                     &
}
\xymatrix@R=2.2pc@C=1.5pc{ 
&\\            
              &  &         \\
              &  &           \\
              & = & 
}
\xymatrix@R=1.5pc@C=0.5pc{
                                                   &                                         &&             C \ar@{-}[dll] \ar@{-}[drr]        &                                \\
                                                  &  C  \ar@{-}[d]               &&                                     && C\ar@{-}[d]              &          \\
                                                  & p \ar[d]                     &&                                       & &q \ar[d]                     &     \\
                                                  & A \ar@{-}@/_/[drr]   &&                                       &  & B \ar@{-}@/^/[dll]  &                         \\
                                                  &                                      &&  f \ar[d]                 &                                     &                    \\
                                                  &                                     & &  D                                   &                                     &
} 
\]
thus showing unicity.
Furthermore, this $t$ gives indeed the desired factorization:
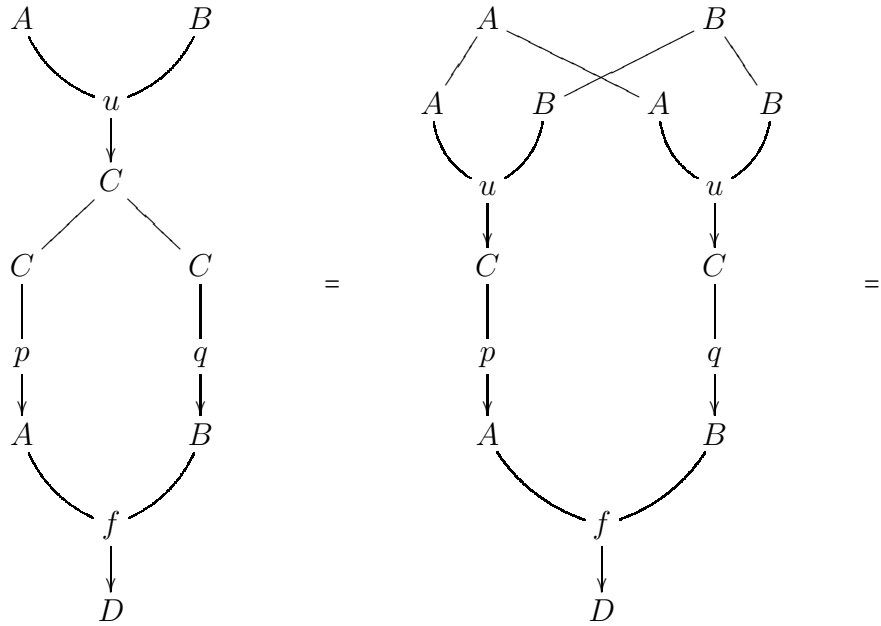
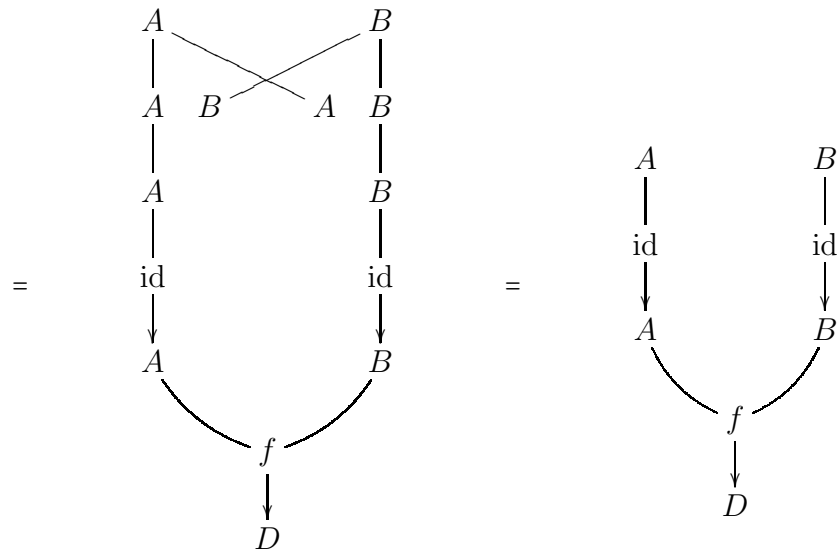
\[
\xymatrix@R=1.5pc@C=0.5pc{
                                                  & A \ar@{-}@/_/[drr]   &&                                       &  & B \ar@{-}@/^/[dll]  &                         \\
                                                  &                                      &&  u \ar[d]                 &                                     &                    \\
                                                  &                                     & &  C  \ar@{-}[dll] \ar@{-}[drr]                          &                                     &           \\
                                                 &  C  \ar@{-}[d]               &&                                     && C\ar@{-}[d]              &          \\
                                                  & p \ar[d]                     &&                                       & &q \ar[d]                     &     \\
                                                  & A \ar@{-}@/_/[drr]   &&                                       &  & B \ar@{-}@/^/[dll]  &                         \\
                                                  &                                      &&  f \ar[d]                 &                                     &                    \\
                                                  &                                     & &  D                                   &                                     &
}
\xymatrix@R=2.2pc@C=1.5pc{ 
&\\            
              &  &         \\
              &  &           \\
              & = & 
}
\xymatrix@R=1.5pc@C=0.5pc{
                                                                           &A\ar@{-}[dl]\ar@{-}[drrr]         &&&&       B   \ar@{-}[dr]\ar@{-}[dlll]             &    \\
                                                   A \ar@{-}@/_/[dr]   && B \ar@{-}@/^/[dl]        && A \ar@{-}@/_/[dr] && B \ar@{-}@/^/[dl]     \\
                                                  &  u\ar[d]                        &&                                     && u\ar[d]                                      &           \\
                                                 &  C  \ar@{-}[d]               &&                                     && C\ar@{-}[d]              &          \\
                                                  & p \ar[d]                     &&                                        &&q \ar[d]                     &     \\
                                                  & A \ar@{-}@/_/[drr]   &&                                       && B \ar@{-}@/^/[dll]  &                         \\
                                                  &                                      &&  f \ar[d]                 &                                     &                    \\
                                                  &                                     & &  D                                   &                                     &
}
\xymatrix@R=2.2pc@C=1.5pc{ 
&\\            
              &  &         \\
              &  &           \\
              & = & 
}
\]

\[
\xymatrix@R=2.2pc@C=1.5pc{ 
&\\            
              &  &         \\
              &  &           \\
              & = & 
}
\xymatrix@R=1.5pc@C=0.5pc{
                                                                           &A\ar@{-}[d]\ar@{-}[drrr]         &&&&       B   \ar@{-}[d]\ar@{-}[dlll]             &    \\
                                                  & A \ar@{-}[d]   & B                  && A                      & B \ar@{-}[d]     \\
                                                 &  A  \ar@{-}[d]               &&                                     && B\ar@{-}[d]              &          \\
                                                  & \id \ar[d]                     &&                                     &&\id \ar[d]                     &     \\
                                                  & A \ar@{-}@/_/[drr]   &&                                       && B \ar@{-}@/^/[dll]  &                         \\
                                                  &                                      &&  f \ar[d]                 &                                     &                    \\
                                                  &                                     & &  D                                   &                                     &
}
\xymatrix@R=2.2pc@C=1.5pc{ 
&\\            
              &  &         \\
              &  &           \\
              & = & 
}
\xymatrix@R=1.5pc@C=0.5pc{
&\\
&\\
                                                  &  A  \ar@{-}[d]               &&                                     && B\ar@{-}[d]              &          \\
                                                  & \id \ar[d]                     &&                                     &&\id \ar[d]                     &     \\
                                                  & A \ar@{-}@/_/[drr]   &&                                       && B \ar@{-}@/^/[dll]  &                         \\
                                                  &                                      &&  f \ar[d]                 &                                     &                    \\
                                                  &                                     & &  D                                   &                                     &
}
\]

Similarly, to prove that (1) implies (2)
we show that in an algebraic product $p:C\to A$, $q:C\to B$ and $u:A,B\to C$, 
$p$ and $q$ are the projection for a universal product.
Suppose that we have arrows (unary, for graphical simplicity) $f:D\to A$ and $f:D\to B$.
If they factor as $f=pt$ and $g=qt$, then $t$ is equal to
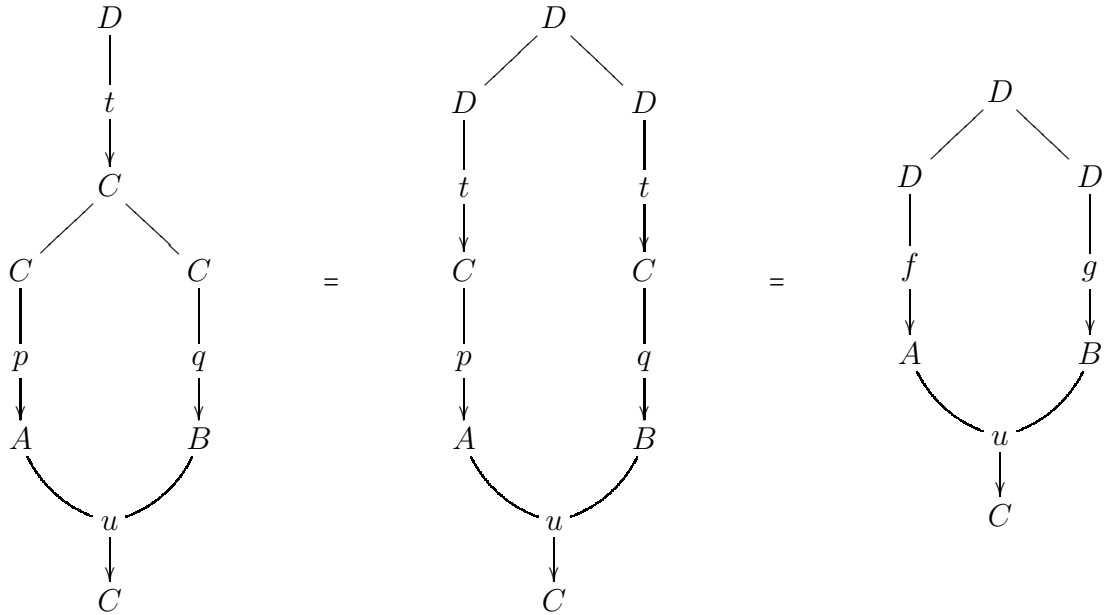
\[
\xymatrix@R=1.5pc@C=0.5pc{
                                                  &                                     & &  D  \ar@{-}[d]                                  &                                     &               \\
                                                  &                                     & &  t    \ar[d]                                &                                     &     \\
                                                   &                                         &&             C \ar@{-}[dll] \ar@{-}[drr]        &                                \\
                                                  &  C  \ar@{-}[d]               &&                                     && C\ar@{-}[d]              &          \\
                                                  & p \ar[d]                     &&                                       & &q \ar[d]                     &     \\
                                                  & A \ar@{-}@/_/[drr]   &&                                       &  & B \ar@{-}@/^/[dll]  &                         \\
                                                  &                                      &&  u \ar[d]                 &                                     &                    \\
                                                  &                                     & &  C                                    &                                     &           
}
\xymatrix@R=2.2pc@C=1.5pc{ 
&\\            
              &  &         \\
              &  &           \\
              & = & 
}
\xymatrix@R=1.5pc@C=0.5pc{
                                                   &                                         &&             D \ar@{-}[dll] \ar@{-}[drr]        &                                \\
                                                  &   D  \ar@{-}[d]          &&                       &&                          D  \ar@{-}[d]           &               \\
                                                  &   t    \ar[d]                      &&                  &&                       t    \ar[d]                       \\
                                                  &  C  \ar@{-}[d]               &&                                     && C\ar@{-}[d]              &          \\
                                                  & p \ar[d]                     &&                                       & &q \ar[d]                     &     \\
                                                  & A \ar@{-}@/_/[drr]   &&                                       &  & B \ar@{-}@/^/[dll]  &                         \\
                                                  &                                      &&  u \ar[d]                 &                                     &                    \\
                                                  &                                     & &  C                                    &                                     &           
}
\xymatrix@R=2.2pc@C=1.5pc{ 
&\\            
              &  &         \\
              &  &           \\
              & = & 
}
\xymatrix@R=1.5pc@C=0.5pc{
&\\
                                                   &                                         &&             D \ar@{-}[dll] \ar@{-}[drr]        &                                \\
                                                  &   D  \ar@{-}[d]          &&                       &&                          D  \ar@{-}[d]           &               \\
                                                  &   f    \ar[d]                      &&                  &&                       g    \ar[d]                       \\
                                                  & A \ar@{-}@/_/[drr]   &&                                       &  & B \ar@{-}@/^/[dll]  &                         \\
                                                  &                                      &&  u \ar[d]                 &                                     &                    \\
                                                  &                                     & &  C                                    &                                     &           
}
\]
and this $t$ gives indeed the desired factorization:
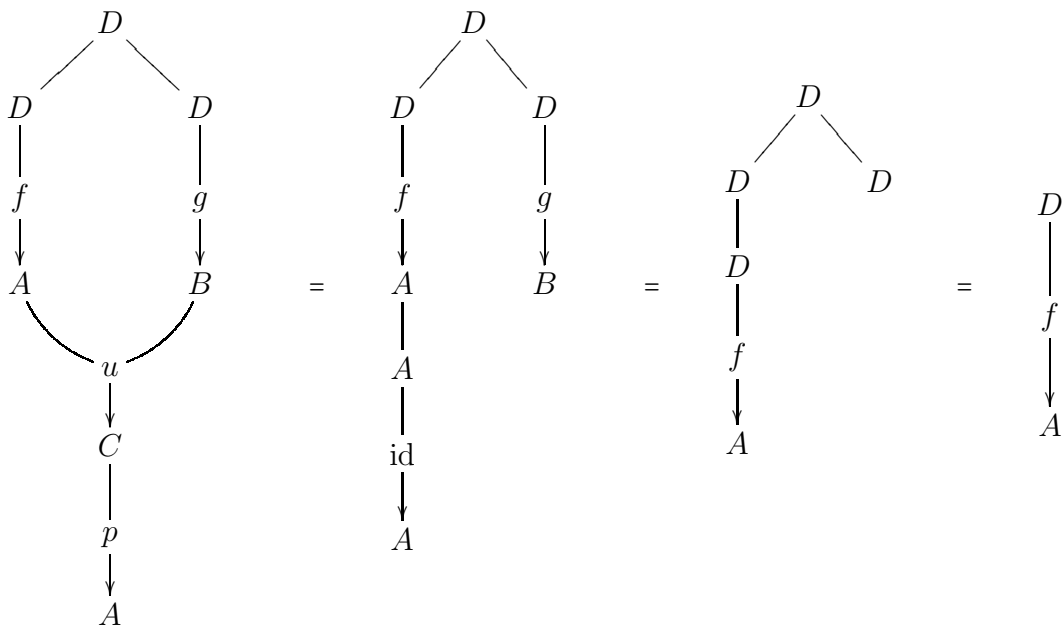
\[
\xymatrix@R=1.5pc@C=0.5pc{
                                                                                            &&             D \ar@{-}[dll] \ar@{-}[drr]        &                                \\
                                                     D  \ar@{-}[d]          &&                       &&                          D  \ar@{-}[d]           &               \\
                                                     f    \ar[d]                      &&                  &&                       g    \ar[d]                       \\
                                                   A \ar@{-}@/_/[drr]   &&                                       &  & B \ar@{-}@/^/[dll]  &                         \\
                                                                                        &&  u \ar[d]                 &                                     &                    \\
                                                                                       & &  C  \ar@{-}[d]                                &                                     &     \\      
                                                                                       & &  p  \ar[d]                                &                                     &          \\ 
                                                                                       & &  A                                &                                     &           
}
\xymatrix@R=2.2pc@C=1pc{ 
\\            
                &         \\
                &           \\
              & = & 
}
\xymatrix@R=1.5pc@C=0.2pc{
                                                                                            &&             D \ar@{-}[dll] \ar@{-}[drr]        &                                \\
                                                     D  \ar@{-}[d]          &&                       &&                          D  \ar@{-}[d]           &               \\
                                                     f    \ar[d]                      &&                  &&                       g    \ar[d]                       \\
                                                   A \ar@{-}[d]   &&                                       &  & B   &                         \\
                                                                                         A  \ar@{-}[d]                                &                                     &     \\      
                                                                                         \id  \ar[d]                                &                                     &          \\ 
                                                                                         A                                &                                     &           
}
\xymatrix@R=2.2pc@C=1pc{ 
\\            
                &         \\
                &           \\
             &  = & 
}
\xymatrix@R=1.5pc@C=0.2pc{
&\\
                                                                                            &&             D \ar@{-}[dll] \ar@{-}[drr]                                       \\
                                                     D  \ar@{-}[d]          &&                       &&                          D                         \\
                                                     D  \ar@{-}[d]          &&                       &&                                                   \\
                                                     f    \ar[d]                      &&                  &&                                             \\
                                                   A    &&                                       &                          
}
\xymatrix@R=2.2pc@C=1pc{ 
\\            
                &         \\
                &           \\
             &  = & 
}
\xymatrix@R=2.2pc@C=0.5pc{ 
\\ 
\\           
                D  \ar@{-}[d]         \\
                 f    \ar[d]             \\
                A 
}
\]
In the other direction, to show that (3) implies (1), suppose that $u:A,B\to C$ is preuniversal
and define $p:C\to A$ and $q:C\to B$ as the unique arrows fulfilling 
equations (\ref{algp2}) in the definition of algebraic product.
We so get indeed an algebraic product for $A$ and $B$, 
since condition (\ref{algp1}) follows again from the preuniversality of $u$:
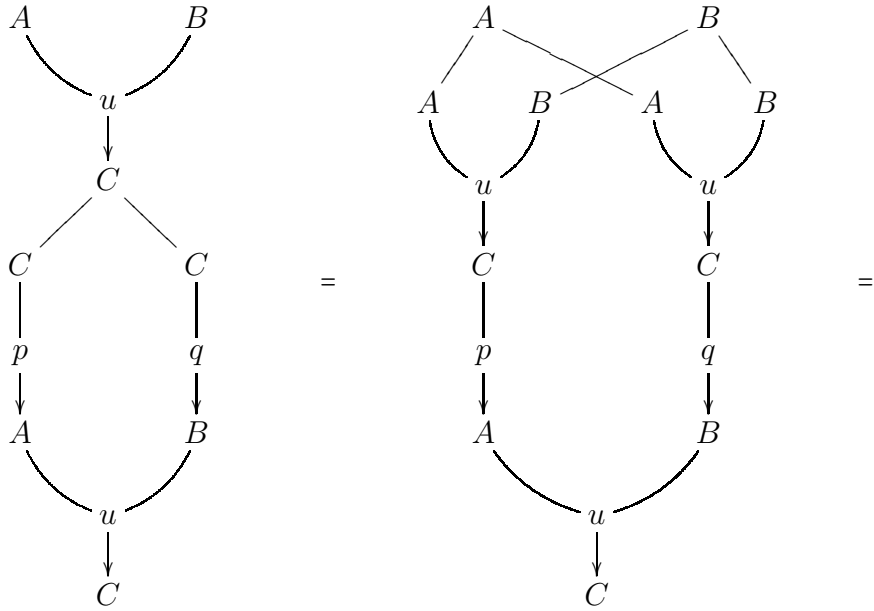
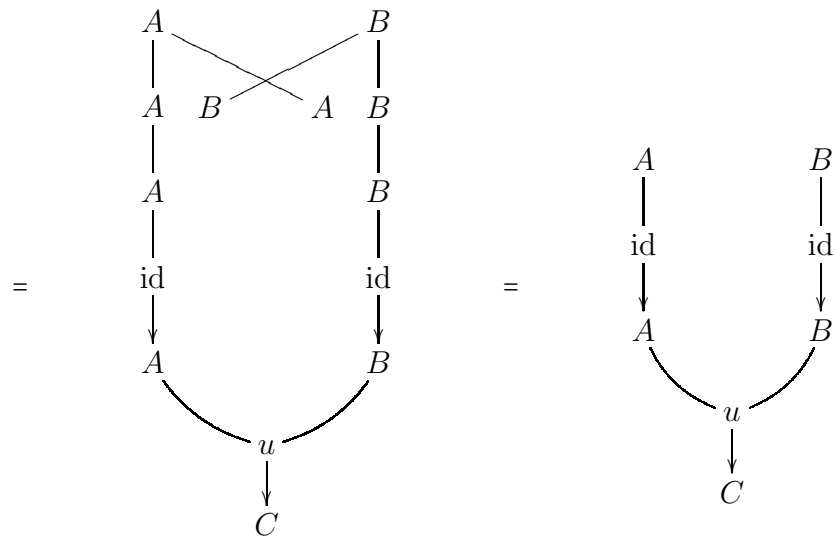
\[
\xymatrix@R=1.5pc@C=0.5pc{
                                                  & A \ar@{-}@/_/[drr]   &&                                       &  & B \ar@{-}@/^/[dll]  &                         \\
                                                  &                                      &&  u \ar[d]                 &                                     &                    \\
                                                  &                                     & &  C  \ar@{-}[dll] \ar@{-}[drr]                          &                                     &           \\
                                                 &  C  \ar@{-}[d]               &&                                     && C\ar@{-}[d]              &          \\
                                                  & p \ar[d]                     &&                                       & &q \ar[d]                     &     \\
                                                  & A \ar@{-}@/_/[drr]   &&                                       &  & B \ar@{-}@/^/[dll]  &                         \\
                                                  &                                      &&  u \ar[d]                 &                                     &                    \\
                                                  &                                     & &  C                                   &                                     &
}
\xymatrix@R=2.2pc@C=1.5pc{ 
&\\            
              &  &         \\
              &  &           \\
              & = & 
}
\xymatrix@R=1.5pc@C=0.5pc{
                                                                           &A\ar@{-}[dl]\ar@{-}[drrr]         &&&&       B   \ar@{-}[dr]\ar@{-}[dlll]             &    \\
                                                   A \ar@{-}@/_/[dr]   && B \ar@{-}@/^/[dl]        && A \ar@{-}@/_/[dr] && B \ar@{-}@/^/[dl]     \\
                                                  &  u\ar[d]                        &&                                     && u\ar[d]                                      &           \\
                                                 &  C  \ar@{-}[d]               &&                                     && C\ar@{-}[d]              &          \\
                                                  & p \ar[d]                     &&                                        &&q \ar[d]                     &     \\
                                                  & A \ar@{-}@/_/[drr]   &&                                       && B \ar@{-}@/^/[dll]  &                         \\
                                                  &                                      &&  u \ar[d]                 &                                     &                    \\
                                                  &                                     & &  C                                   &                                     &
}
\xymatrix@R=2.2pc@C=1.5pc{ 
&\\            
              &  &         \\
              &  &           \\
              & = & 
}
\]

\[
\xymatrix@R=2.2pc@C=1.5pc{ 
&\\            
              &  &         \\
              &  &           \\
              & = & 
}
\xymatrix@R=1.5pc@C=0.5pc{
                                                                           &A\ar@{-}[d]\ar@{-}[drrr]         &&&&       B   \ar@{-}[d]\ar@{-}[dlll]             &    \\
                                                  & A \ar@{-}[d]   & B                  && A                      & B \ar@{-}[d]     \\
                                                 &  A  \ar@{-}[d]               &&                                     && B\ar@{-}[d]              &          \\
                                                  & \id \ar[d]                     &&                                     &&\id \ar[d]                     &     \\
                                                  & A \ar@{-}@/_/[drr]   &&                                       && B \ar@{-}@/^/[dll]  &                         \\
                                                  &                                      &&  u \ar[d]                 &                                     &                    \\
                                                  &                                     & &  C                                   &                                     &
}
\xymatrix@R=2.2pc@C=1.5pc{ 
&\\            
              &  &         \\
              &  &           \\
              & = & 
}
\xymatrix@R=1.5pc@C=0.5pc{
&\\
&\\
                                                  &  A  \ar@{-}[d]               &&                                     && B\ar@{-}[d]              &          \\
                                                  & \id \ar[d]                     &&                                     &&\id \ar[d]                     &     \\
                                                  & A \ar@{-}@/_/[drr]   &&                                       && B \ar@{-}@/^/[dll]  &                         \\
                                                  &                                      &&  u \ar[d]                 &                                     &                    \\
                                                  &                                     & &  C                                   &                                     &
}
\]
Finally, to show that (2) implies (1), assume that $p:C\to A$ and $q:C\to B$ are a universal product for $A$ and $B$
and let $u:A,B\to C$ be the unique arrow fulfilling equations (\ref{algp2}) in the definition of algebraic product.
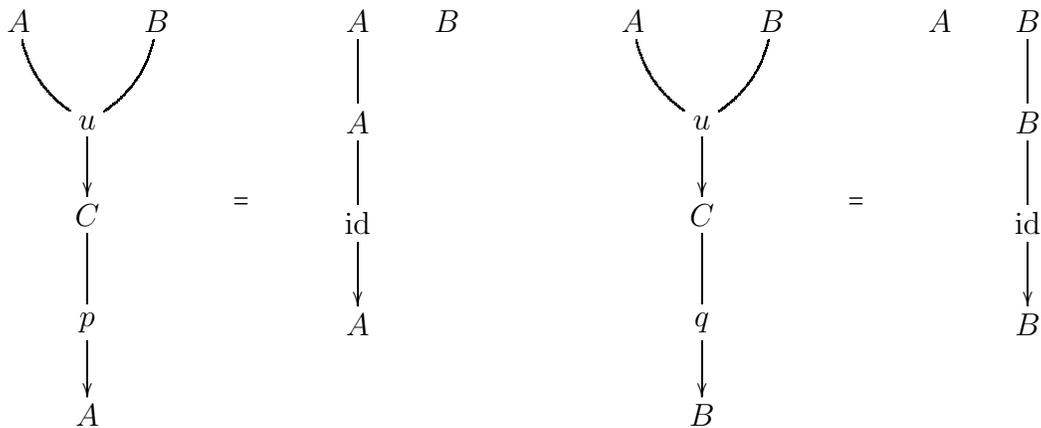
\[
\xymatrix@R=2pc@C=0.2pc{
                                                  & A \ar@{-}@/_/[drr]   &&                                       &  & B \ar@{-}@/^/[dll]                           \\
                                                  &                                      &&  u \ar[d]                 &                                                         \\
                                                  &                                     & &  C    \ar@{-}[d]                                                               \\
                                                  &                                     &&    p \ar[d]                                            \\
                                                  &                                    &&    A                                     
}
\xymatrix@R=2.2pc@C=1pc{             
              &  &         \\
              &  &           \\
              & = & 
}
\xymatrix@R=2pc@C=0.5pc{
                                                  &  A  \ar@{-}[d]  &&  B                                             \\
                                                  &  A  \ar@{-}[d]                        \\
                                                  & \id \ar[d]                                  \\
                                                  & A                       
}
\qq \qq
\xymatrix@R=2pc@C=0.2pc{
                                                  & A \ar@{-}@/_/[drr]   &&                                       &  & B \ar@{-}@/^/[dll]                          \\
                                                  &                                      &&  u \ar[d]                 &                                                       \\
                                                  &                                     & &  C    \ar@{-}[d]                                                               \\
                                                  &                                     &&    q \ar[d]                                            \\
                                                  &                                    &&    B                                     
}
\xymatrix@R=2.2pc@C=1pc{ 
              &  &         \\
              &  &           \\
              & = & 
}
\xymatrix@R=2pc@C=0.5pc{
                                                  A   &&  B \ar@{-}[d]                    \\
                                                  &&  B  \ar@{-}[d]                        \\
                                                  && \id \ar[d]                                  \\
                                                  && B                       
}
\]
That we so get indeed an algebraic product follows now from the universality of projections, since:
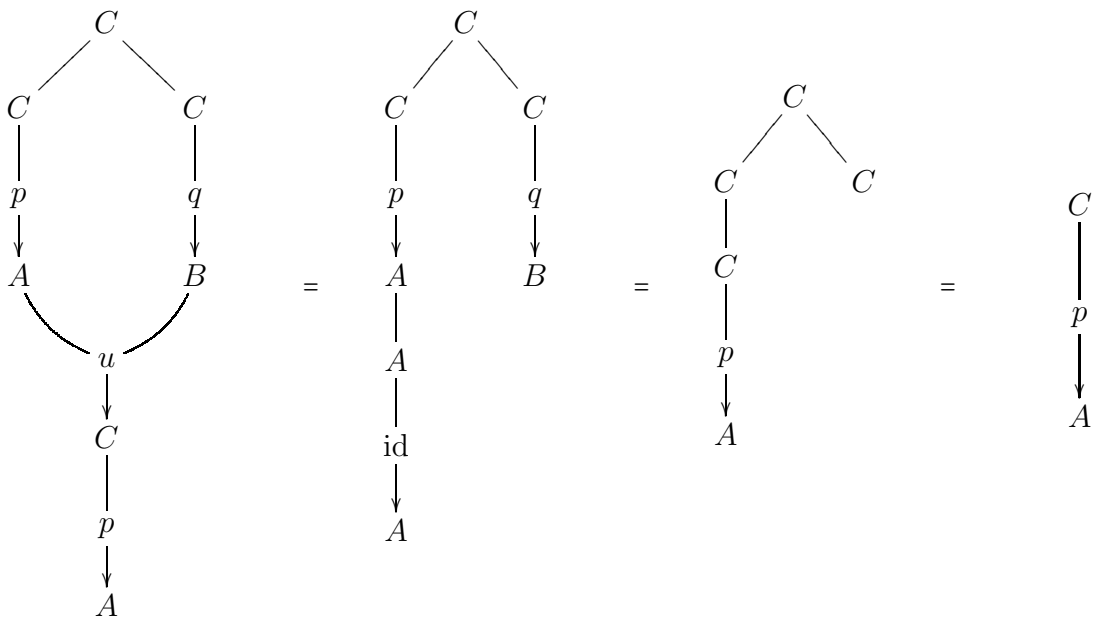
\[
\xymatrix@R=1.5pc@C=0.5pc{
                                                                                            &&             C \ar@{-}[dll] \ar@{-}[drr]        &                                \\
                                                     C  \ar@{-}[d]          &&                       &&                          C  \ar@{-}[d]           &               \\
                                                     p    \ar[d]                      &&                  &&                       q    \ar[d]                       \\
                                                   A \ar@{-}@/_/[drr]   &&                                       &  & B \ar@{-}@/^/[dll]  &                         \\
                                                                                        &&  u \ar[d]                 &                                     &                    \\
                                                                                       & &  C  \ar@{-}[d]                                &                                     &     \\      
                                                                                       & &  p  \ar[d]                                &                                     &          \\ 
                                                                                       & &  A                                &                                     &           
}
\xymatrix@R=2.2pc@C=1pc{ 
&\\            
              &  &         \\
              &  &           \\
              & = & 
}
\xymatrix@R=1.5pc@C=0.2pc{
                                                                                            &&             C \ar@{-}[dll] \ar@{-}[drr]        &                                \\
                                                     C  \ar@{-}[d]          &&                       &&                          C  \ar@{-}[d]           &               \\
                                                     p    \ar[d]                      &&                  &&                       q    \ar[d]                       \\
                                                   A \ar@{-}[d]   &&                                       &  & B   &                         \\
                                                                                         A  \ar@{-}[d]                                &                                     &     \\      
                                                                                         \id  \ar[d]                                &                                     &          \\ 
                                                                                         A                                &                                     &           
}
\xymatrix@R=2.2pc@C=1pc{ 
&\\            
              &  &         \\
              &  &           \\
              & = & 
}
\xymatrix@R=1.5pc@C=0.2pc{
&\\
                                                                                            &&             C \ar@{-}[dll] \ar@{-}[drr]                                       \\
                                                     C  \ar@{-}[d]          &&                       &&                          C                         \\
                                                     C  \ar@{-}[d]          &&                       &&                                                   \\
                                                     p    \ar[d]                      &&                  &&                                             \\
                                                   A    &&                                       &                          
}
\xymatrix@R=2.2pc@C=1pc{ 
&\\            
              &  &         \\
              &  &           \\
              & = & 
}
\xymatrix@R=2.2pc@C=1pc{ 
\\
&\\            
              & C \ar@{-}[d] &         \\
              & p \ar[d] &           \\
              & A & 
}
\]
and similarly with $q$ in place of $p$.
\epf


Suppose now that $p:C\to A$, $q:C\to B$, $u:A,B\to C$ and $l:E\to C$, $r:E\to D$, $v:C,D\to E$
are algebraic products in $\M$.
One may wonder if their \"composite" 
\[
pl:E\to A \qv ql:E\to A \qv r:E\to D \qv v(u,D):A,B,D \to E
\]
is again an algebraic product for $A$, $B$ and $D$.
The answer is affirmative.

\begin{prop}
\label{cart1b}
Algebraic products in a cartesian multicategory are closed with respect to composition.
\end{prop}
\pf
We give the proof in the above case of composition of binary products.
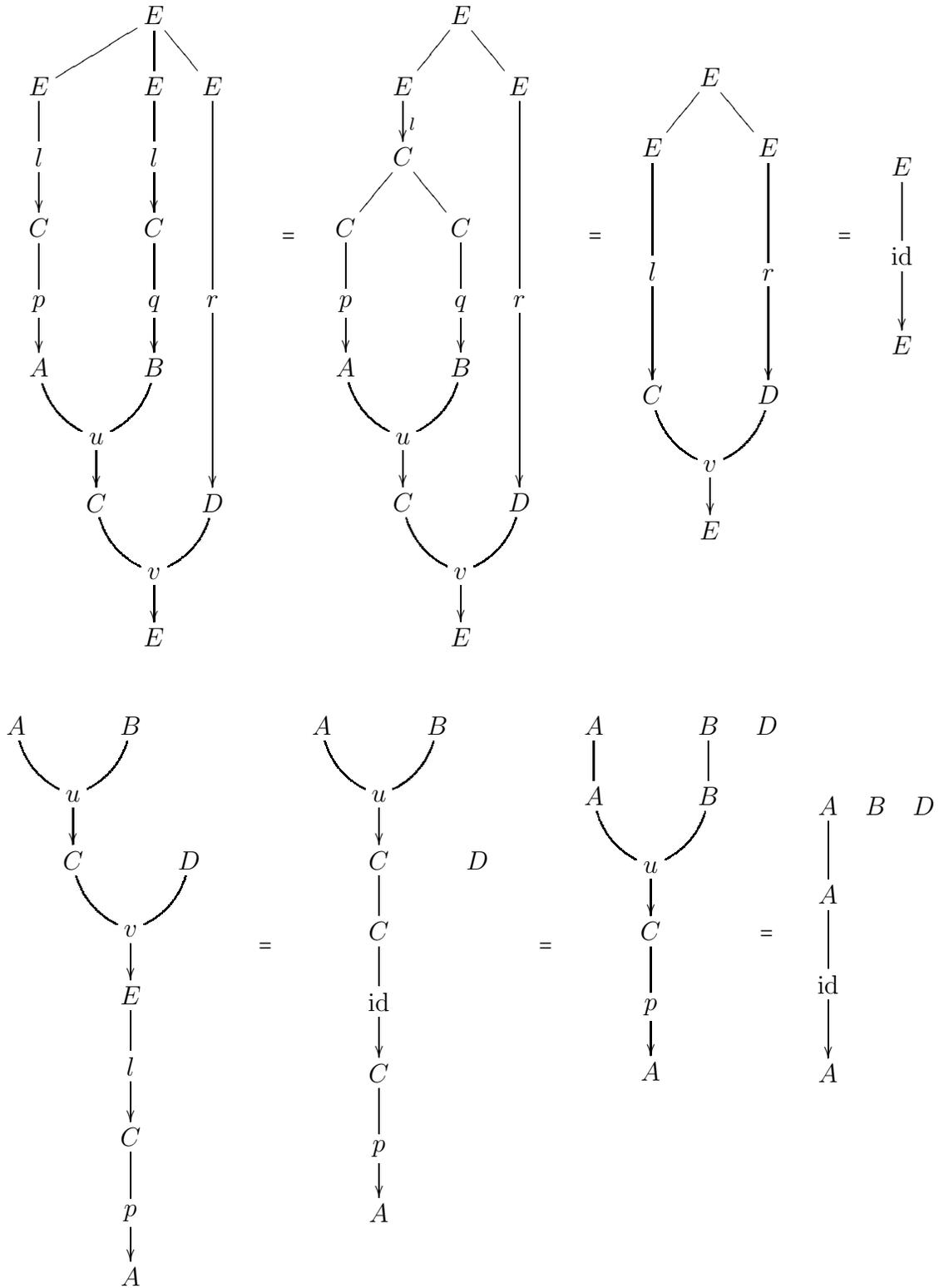
\[
\xymatrix@R=1.5pc@C=0.2pc{
                                                                                            & & & &          E \ar@{-}[dllll] \ar@{-}[d]  \ar@{-}[drr]              \\
                                                     E  \ar@{-}[d]          &&                       &&                          E  \ar@{-}[d]      &&       E   \ar@{-}[ddd]    \\
                                                     l    \ar[d]                      &&                  &&                       l    \ar[d]                       \\
                                                     C  \ar@{-}[d]          &&                       &&                          C  \ar@{-}[d]           &               \\
                                                     p    \ar[d]                      &&                  &&                       q    \ar[d]       &&       r    \ar[ddd]            \\
                                                   A \ar@{-}@/_/[drr]   &&                                       &  & B \ar@{-}@/^/[dll]  &                         \\
                                                                                        &&  u \ar[d]                 &                                     &                    \\
                                                             & &  C \ar@{-}@/_/[drr]   &&                                       &  & D \ar@{-}@/^/[dll]  &                         \\
                                                                                       && &&  v \ar[d]                 &                                     &                    \\
                                                                                       && &&  E                 &                                     &                    
                                                                                       }
\xymatrix@R=2.2pc@C=0.5pc{ 
&\\            
              &  &         \\
              &  &           \\
              & = & 
}
\xymatrix@R=1.5pc@C=0.2pc{
                                                                                            & & & &          E \ar@{-}[dll]   \ar@{-}[drr]              \\
                                                                              &&               E  \ar[d]^l      &&                 &&       E   \ar@{-}[ddd]    \\
                                                                          &&          C \ar@{-}[dll] \ar@{-}[drr]      &&                                            \\
                                                     C  \ar@{-}[d]          &&                       &&                          C  \ar@{-}[d]           &               \\
                                                     p    \ar[d]                      &&                  &&                       q    \ar[d]       &&       r    \ar[ddd]            \\
                                                   A \ar@{-}@/_/[drr]   &&                                       &  & B \ar@{-}@/^/[dll]  &                         \\
                                                                                        &&  u \ar[d]                 &                                     &                    \\
                                                             & &  C \ar@{-}@/_/[drr]   &&                                       &  & D \ar@{-}@/^/[dll]  &                         \\
                                                                                       && &&  v \ar[d]                 &                                     &                    \\
                                                                                       && &&  E                 &                                     &                    
                                                                                       }
\xymatrix@R=2.2pc@C=0.5pc{ 
&\\            
              &  &         \\
              &  &           \\
              & = & 
}
\xymatrix@R=1.5pc@C=0.2pc{
\\
                                                                                             & &          E \ar@{-}[dll]   \ar@{-}[drr]              \\
                                                                                             E  \ar@{-}[dd]      &&                 &&       E   \ar@{-}[dd]    \\
                                                                                        &&                                            \\
                                                                                  l  \ar[dd]         &&                              &&       r    \ar[dd]            \\
                                                                                        &&                   &                                     &                    \\
                                                               C \ar@{-}@/_/[drr]   &&                                       &  & D \ar@{-}@/^/[dll]  &                         \\
                                                                                        &&  v \ar[d]                 &                                     &                    \\
                                                                                        &&  E                 &                                     &                    
                                                                                       }
\xymatrix@R=2.2pc@C=0.5pc{ 
&\\            
              &  &         \\
              &  &           \\
              & = & 
}
\xymatrix@R=2.2pc@C=0.5pc{ 
\\
\\
E \ar@{-}[d]          \\            
      \rm{id}\ar[d]   \\
              E           
}
\]
\[
\xymatrix@R=1.5pc@C=0.2pc{
                                                                                                       A \ar@{-}@/_/[drr]   &&                                       &  & B \ar@{-}@/^/[dll]  &          \\
                                                                                        &&  u \ar[d]                 &                                     &                    \\
                                                             & &  C \ar@{-}@/_/[drr]   &&                                       &  & D \ar@{-}@/^/[dll]  &                         \\
                                                                                       && &&  v \ar[d]                 &                                     &                    \\
                                                                                       && &&  E   \ar@{-}[d]              &                                     &              \\      
                                                                                       && &&                       l    \ar[d]       &            \\
                                                                                       && &&                       C  \ar@{-}[d]           &     \\         
                                                                                       && &&                       p    \ar[d]       &            \\
                                                                                       && &&                       A             &             
}
\xymatrix@R=2.2pc@C=0.5pc{ 
&\\            
              &  &         \\
              &  &           \\
              & = & 
}
\xymatrix@R=1.5pc@C=0.2pc{
                                                                                                       A \ar@{-}@/_/[drr]   &&                                       &  & B \ar@{-}@/^/[dll]  &          \\
                                                                                       &&  u \ar[d]                 &                                     &                    \\
                                                                                       &&  C \ar@{-}[d]   &&                                       & D   &                         \\
                                                                                       &&  C   \ar@{-}[d]              &                                     &             \\      
                                                                                       &&                       \rm id    \ar[d]       &            \\
                                                                                       &&                       C  \ar@{-}[d]           &     \\         
                                                                                       &&                       p    \ar[d]       &            \\
                                                                                       &&                       A             &             
}
\xymatrix@R=2.2pc@C=0.2pc{ 
&\\            
              &  &         \\
              &  &           \\
              & = & 
}
\xymatrix@R=1.5pc@C=0.2pc{
                                                                                                       A \ar@{-}[d]   &&                                         &  & B \ar@{-}[d]  &&    D      \\
                                                                                                       A \ar@{-}@/_/[drr]   &&                                       &  & B \ar@{-}@/^/[dll]  &          \\
                                                                                       &&  u \ar[d]                 &                                     &                    \\
                                                                                       &&  C \ar@{-}[d]   &&                                       &  &            =                \\
                                                                                       &&                       p    \ar[d]       &            \\
                                                                                       &&                       A             &             
}
\xymatrix@R=2.2pc@C=0.5pc{ 
\\
&A \ar@{-}[d]  & B & D       \\ 
&A \ar@{-}[d]          \\            
&      \rm{id}\ar[d]   \\
&              A           
}
\]
(and similarly for $q$ in place of $p$).
\[
\xymatrix@R=1.5pc@C=0.2pc{
                                                                                                       A \ar@{-}@/_/[drr]   &&                                       &  & B \ar@{-}@/^/[dll]  &          \\
                                                                                        &&  u \ar[d]                 &                                     &                    \\
                                                             & &  C \ar@{-}@/_/[drr]   &&                                       &  & D \ar@{-}@/^/[dll]  &                         \\
                                                                                       && &&  v \ar[d]                 &                                     &                    \\
                                                                                       && &&  E   \ar@{-}[d]              &                                     &              \\      
                                                                                       && &&                       r    \ar[d]       &            \\
                                                                                       && &&                       D             &              
}
\xymatrix@R=2.2pc@C=0.5pc{ 
&\\            
              &  &         \\
              &  &           \\
              & = & 
}
\xymatrix@R=1.5pc@C=0.2pc{
                                                                                                       A \ar@{-}@/_/[drr]   &&                                       &  & B \ar@{-}@/^/[dll]  &          \\
                                                                                       &&  u \ar[d]                 &                                     &                    \\
                                                                                       &&  C    &&                                       &   D \ar@{-}[d]  &                         \\
                                                                                       && &&                   &                               D  \ar@{-}[d]               \\
                                                                                       && &&                     &  \rm id    \ar[d]       &            \\
                                                                                       && &&                     &    D             &              
}
\xymatrix@R=2.2pc@C=0.5pc{ 
&\\            
              &  &         \\
              &  &           \\
              & = & 
}
\xymatrix@R=2.2pc@C=0.5pc{ 
\\
&A & B & D  \ar@{-}[d]      \\ 
&&& D \ar@{-}[d]          \\            
&&&      \rm{id}\ar[d]   \\
&&&              D           
}
\]

\epf

\begin{corol}
\label{cart2}
If a cartesian multicategory $\M$ is representable then $\M_-$ has finite products.
If $\M$ has algebraic products, then it is representable.
Any functor $F:\M\to\N$ in $\fpMlt$ preserves algebraic products, preuniversal arrows 
and universal products that exist in $\M$.
\end{corol}
\pf
The first part follows directly from Proposition \ref{cart1}, since universal products in $\M$
are in particular products in $\M_-$.
For the second part, by Proposition \ref{cart1} every finite family is the domain of
a preuniversal arrow and by Proposition \ref{cart1b} these arrows are closed
with respect to composition. 
Since algebraic products are clearly preserved by functors in $\fpMlt$,
the third part follows again from Proposition \ref{cart1}. 
\epf


\subsection{Finite product categories as cartesian multicategories}

As sketched at the beginning of this section, finite product categories give rise to  
cartesian multicategories.
In fact, if we denote by $\fpCat$ the category of finite product categories and 
finite product preserving functors, we have:
\begin{prop}  \label{cart3}
The category $\fpRep$ of representable cartesian multicategories is equivalent to $\fpCat$.
\end{prop}
\pf
Given a finite product category $\C$, let us fix a universal cone $p_i:A_1\tm \cd \tm A_n \to A_i$
for any finite family $A_1,\cd,A_n$.
Then we define $\C_\tm \in\fpRep$ by
\[ 
\C_\tm(A_1,\cd,A_n;A) \iso \C(A_1\tm \cd \tm A_n;A)
\]
with the obvious composition and the following cartesian structure.
If $f:A_1,\cd, A_n \to A$, then we define $\si f =f\ov\si:B_1,\cd,B_m\to A $,
where $\ov\si: B_1\tm\cd\tm B_m \to A_1\tm\cd\tm A_n$
is the map such that $p_i\ov\si = q_{\si i}$.
Conversely, given a representation on $\M\in\fpRep$
we get, as in Proposition \ref{cart1}, a universal cone in $\M_-$ for any finite family $A_1,\cd,A_n$.
Now, it is straightforward to check that $(\C_\tm)_- \iso \C$ and that,
in the other direction, we have an isomorphism of multicategories $(\M_-)_\tm \iso \M$ 
which takes an arrow $f:A_1\tm\cd\tm A_n \to A$ to $fu:A_1,\cd,A_n \to A$ 
(where $u$ is the corresponding universal arrow in the representation).
To show that this isomorphism respects the cartesian structure, by Remark \ref{cart0}
we only need to check that it respects permutations and contractions,
which is also straightforward.

Both the constructions $\C_\tm$ and $\M_-$ extend to functors (recall Corollary \ref{cart2})
giving the desired equivalence $\fpCat \simeq \fpRep$.
\epf

\begin{corol}
\label{cart4}
For a symmetric multicategory $\M$ the following are equivalent:
\begin{enumerate}
\item
$\M$ is cartesian and has algebraic products.
\item
$\M$ is cartesian and has finite universal products.
\item
$\M$ is cartesian and representable.
\item
$\M$ is cartesian monoidal, that is $\M\iso\C_\tm$ for a $\C\in\fpCat$.
\end{enumerate}
\end{corol}
\pf
The equivalence of the first three conditions follows from Proposition \ref{cart1} and Corollary \ref{cart2},
while the equivalence of the last two follows from Proposition \ref{cart3}.
\epf


\subsection{Preadditive categories as cartesian multicategories}

We have just seen that the representable cartesian multicategories
arise from finite product categories.
Presently, we show that the sequential cartesian multicategories 
arise from preadditive categories.
Note that, following other authors, we omit the prefix \"semi" which is sometimes 
used to refer to commutative monoids (rather than abelian groups) enrichments.
Thus, the category $\preAdd$  of preadditive categories and additive functors
is the category $\cMon$-$\Cat$ of categories enriched in commutative monoids.
Similarly, $\Add$ is the full subcategory of those preadditive categories 
with (bi)products.

\begin{prop}  
\label{cart5}
The category $\fpSeq$ of sequential cartesian multicategories is equivalent to $\preAdd$.
\end{prop}
\pf
Giving a preadditive structure on a category $\C$
amounts to giving a cartesian structure on $\C\t$.
Indeed, suppose that $\C\in\preAdd$ and let $f=\la f_1,\cd,f_n\ra:A_1,\cd,A_n\to A$
be an arrow in $\C\t$. We define $\si f=\la g_1,\cd,g_m\ra:B_1,\cd,B_m\to A$ by
\[
g_i = \Sigma_{\si k=i} f_k
\]
and it is straightforward to check that (because of distributivity of composition)
this is indeed a cartesian structure on $\C\t$.

Conversely, if $\C\t$ has a cartesian structure the contraction mappings
give a preadditive structure on $\C$
and the processes are easily seen to be each other's inverse. 
(Again, by Remark \ref{cart0}
it is enough to check that the cartesian structures on $\C\t$ coincide on contractions.)
Thus, the equivalence $(-)_- \adj (-)\t : \Cat \to \Seq$ lifts to the desired equivalence
$\preAdd\simeq\fpSeq$.
\epf

\begin{remark}
\label{cart9}
By propositions \ref{repseq} and \ref{cart3},
we know that a sequential (respectively, cartesian) multicategory $\M$ is representable
iff $\M_-$ has finite sums (respectively, products) iff $\M$ is represented by sums
(respectively, products).
Thus, Proposition \ref{cart5} makes it precise the idea that preadditive categories 
are those categories in which finite sums and finite products coincide, whenever they do exist
(see also Corollary \ref{cart8}).
\end{remark}

\begin{remark}
\label{free}
For a category $\C\in\Cat$, one has
\[
(F\C)\t \iso F(\C\t)
\]
where the $F$ on the the left is the usual free preadditive category functor, 
while that on the right is the free cartesian multicategory functor.
That is, $F\adj U:\fpMlt\to\sMlt$ restricts, on sequential multicategories,
to the usual adjunction $F\adj U:\preAdd\to\Cat$.
\end{remark}


\subsection{Comparison with classical algebraic products}

It is well-known (see for instance \cite{maclane}) that in a preadditive category $\C$ the existence of
finite products and of finite sums are equivalent (to the additivity of $\C$), since both of them
amount to the \"algebraic" (rather than universal) notion of biproduct.
We are now in a position to show that this is just a particular case 
of Corollary \ref{cart4}.

\begin{corol}
\label{cart8}
For a category $\C$ the following are equivalent:
\begin{enumerate}
\item
$\C$ is preadditive and has algebraic biproducts.
\item
$\C$ is preadditive and has finite products.
\item
$\C$ is preadditive and has finite sums.
\item
$\C$ has both finite products and finite sums and they coincide.
\end{enumerate}
\end{corol}
\pf
First observe that the definition of an algebraic product of $A,B\in\M$ 
becomes, for $\M \iso\C\t$ ($\C$ preadditive), an object $C$ along with maps
$p:C\to A$, $q:C\to B$, $i:A\to C$ and $j:B\to C$ such that
\[
\xymatrix@R=1.5pc@C=0.5pc{
                                                   &                                         &&             C \ar@{-}[dll] \ar@{-}[drr]        &                                \\
                                                  &  C  \ar@{-}[d]               &&                                     && C\ar@{-}[d]              &          \\
                                                  & p \ar[d]                     &&                                       & &q \ar[d]                     &     \\
                                                  & A \ar_i[ddrr]   &&                                       &  & B \ar^j[ddll]  &                         \\
                                                  &                                      &&                   &                                     &                    \\
                                                  &                                     & &  C                                    &                                     &
}
\xymatrix@R=2.2pc@C=1.5pc{             
              &  &         \\
              &  &           \\
              & = & 
}
\xymatrix@R=2pc@C=0.5pc{
                                                   &                                                \\
                                                  &  C  \ar@{-}[d]                        \\
                                                  & \id \ar[d]                                  \\
                                                  & C                       
}
\]
which, recalling that contractions in $\C\t$ are given by sums of arrows in $\C$, translates as 
\[
ip+jq=\id_C
\]
and such that
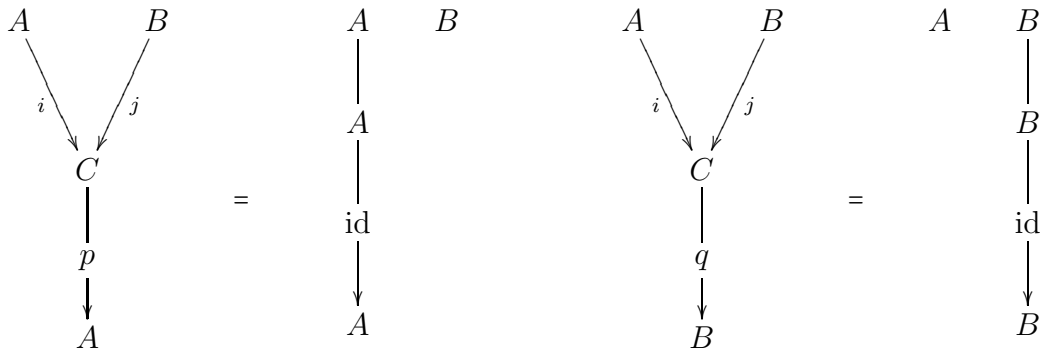
\[
\xymatrix@R=1.5pc@C=0.2pc{
                                                  & A \ar_i[ddrr]   &&                                       &  & B \ar^j[ddll]                           \\
                                                  &                                                        &                                                         \\
                                                  &                                     & &  C    \ar@{-}[d]                                                               \\
                                                  &                                     &&    p \ar[d]                                            \\
                                                  &                                    &&    A                                     
}
\xymatrix@R=2.2pc@C=1pc{             
              &  &         \\
              &  &           \\
              & = & 
}
\xymatrix@R=2pc@C=0.5pc{
                                                  &  A  \ar@{-}[d]  &&  B                                             \\
                                                  &  A  \ar@{-}[d]                        \\
                                                  & \id \ar[d]                                  \\
                                                  & A                       
}
\qq \qq
\xymatrix@R=1.5pc@C=0.2pc{
                                                  & A \ar_i[ddrr]   &&                                       &  & B \ar^j[ddll]                          \\
                                                  &                                      &&                 &                                                       \\
                                                  &                                     & &  C    \ar@{-}[d]                                                               \\
                                                  &                                     &&    q \ar[d]                                            \\
                                                  &                                    &&    B                                     
}
\xymatrix@R=2.2pc@C=1pc{ 
              &  &         \\
              &  &           \\
              & = & 
}
\xymatrix@R=2pc@C=0.5pc{
                                                  A   &&  B \ar@{-}[d]                    \\
                                                  &&  B  \ar@{-}[d]                        \\
                                                  && \id \ar[d]                                  \\
                                                  && B                       
}
\]
which translate as  
\[
pi = \id_A \qv pj = 0_{B,A}  \qv  qj = \id_B  \qv  qi = 0_{A,B}
\]
(The same argument can be of course repeated for any finite family of objects.)
Thus equations (\ref{algp1}) and (\ref{algp2}) become, for $\M\iso\C\t$, the classical
equations giving algebraic biproducts. 
Therefore, when $\M$ is sequential the first condition of Corollary \ref{cart4} 
becomes the first condition listed above.

Now, recalling Proposition \ref{cart5} and Remark \ref{cart9} it is easy
to see that also the other conditions of Corollary \ref{cart4} become, when $\M$ is sequential,
the conditions listed above which are therefore equivalent.
\epf


\subsection{The coreflection in preadditive categories}

\begin{prop}
\label{cart6}
If $\N$ is a cartesian multicategory then, for any $\M\in\sMlt$, 
$[\M,\N]$ is also cartesian in a natural way.
\end{prop}
\pf
The cartesian structure on $[\M,\N]$ is inherited pointwise from that of $\N$.
For instance, the contraction $\al' : F \to G$ of $\al : F, F \to G$ is given by
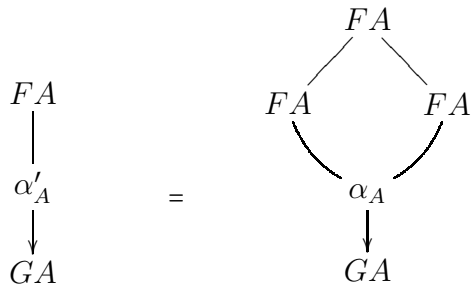
\[
\xymatrix@R=1.5pc@C=0.5pc{
\\
                                                                                        &  FA  \ar@{-}[d]                                           &           \\
                                                                                        &  \al'_A \ar[d]                                              &                    \\
                                                                                       &  GA                                                                 &
}
\xymatrix@R=2.2pc@C=1.5pc{ 
              &  &         \\
              &  &           \\
              & = & 
}
\xymatrix@R=1.5pc@C=0.5pc{
                                                                                        &  FA  \ar@{-}[dl] \ar@{-}[dr]                          &                                     &           \\
                                                   FA \ar@{-}@/_/[dr]   &                                         & FA \ar@{-}@/^/[dl]  &                         \\
                                                                                        &  \al_A \ar[d]                 &                                     &                    \\
                                                                                       &  GA                                   &                                     &
}
\]
Since the conditions for a cartesian multicategory are cleary inherited pointwise from $\N$,
we only have to check that $\al' : F \to G$ is indeed an arrow in $[\M,\N]$:
\[
\xymatrix@R=1.5pc@C=0.3pc{
                                                  & FA \ar@{-}@/_/[drr]   &&                                       &  & FB \ar@{-}@/^/[dll]  &                         \\
                                                  &                                      &&  Ff \ar[d]                 &                                     &                    \\
                                                  &                                     & &  FC  \ar@{-}[dll] \ar@{-}[drr]                          &                                     &           \\
                                                  & FC \ar@{-}@/_/[drr]   &&                                       &  & FC \ar@{-}@/^/[dll]  &                         \\
                                                  &                                      &&  \al_C \ar[d]                 &                                     &                    \\
                                                  &                                     & &  GC                                   &                                     &
}
\xymatrix@R=2.2pc@C=1pc{ 
&\\            
              &  &         \\
              &  &           \\
              & = & 
}
\xymatrix@R=1.5pc@C=0.3pc{
                                                                           &FA\ar@{-}[dl]\ar@{-}[drrr]         &&&&       FB   \ar@{-}[dr]\ar@{-}[dlll]             &    \\
                                                   FA \ar@{-}@/_/[dr]   && FB \ar@{-}@/^/[dl]        && FA \ar@{-}@/_/[dr] && FB \ar@{-}@/^/[dl]     \\
                                                  &  Ff\ar[d]                        &&                                     && Ff\ar[d]                                      &           \\
                                                  & FC \ar@{-}@/_/[drr]   &&                                       && FC \ar@{-}@/^/[dll]  &                         \\
                                                  &                                      &&  \al_C \ar[d]                 &                                     &                    \\
                                                  &                                     & &  GC                                   &                                     &
}
\xymatrix@R=2.2pc@C=1pc{ 
&\\            
              &  &         \\
              &  &           \\
              & = & 
}
\]

\[
\xymatrix@R=2.2pc@C=1pc{ 
&\\            
              &  &         \\
              &  &           \\
              & = & 
}
\xymatrix@R=1.5pc@C=0.3pc{
                                                                           &FA\ar@{-}[dl]\ar@{-}[drrr]         &&&&       FB   \ar@{-}[dr]\ar@{-}[dlll]             &    \\
                                                   FA \ar@{-}[d]            && FB \ar@{-}[drr]        && FA \ar@{-}[dll] && FB \ar@{-}[d]     \\
                                                   FA \ar@{-}@/_/[dr]   && FA \ar@{-}@/^/[dl]        && FB \ar@{-}@/_/[dr] && FB \ar@{-}@/^/[dl]     \\
                                                  &  \al_A\ar[d]                        &&                                     && \al_B\ar[d]                                      &           \\
                                                  & GA \ar@{-}@/_/[drr]   &&                                       && GB \ar@{-}@/^/[dll]  &                         \\
                                                  &                                      &&  Gf \ar[d]                 &                                     &                    \\
                                                  &                                     & &  GC                                   &                                     &
}
\xymatrix@R=2.2pc@C=0.3pc{ 
&\\            
              &  &         \\
              &  &           \\
              & = & 
}
\xymatrix@R=1.5pc@C=0.3pc{
                                                                           &FA\ar@{-}[dl]\ar@{-}[dr]         &&&&       FB   \ar@{-}[dr]\ar@{-}[dl]             &    \\
                                                   FA \ar@{-}@/_/[dr]   && FA \ar@{-}@/^/[dl]        && FB \ar@{-}@/_/[dr] && FB \ar@{-}@/^/[dl]     \\
                                                  &  \al_A\ar[d]                        &&                                     && \al_B\ar[d]                                      &           \\
                                                  & GA \ar@{-}@/_/[drr]   &&                                       && GB \ar@{-}@/^/[dll]  &                         \\
                                                  &                                      &&  Gf \ar[d]                 &                                     &                    \\
                                                  &                                     & &  GC                                   &                                     &
}
\]
Similarly, for weakenings we have:
\[
\xymatrix@R=1.5pc@C=0.3pc{
                                                   FA \ar@{-}@/_/[dr]   && FB \ar@{-}@/^/[dl]        && HA \ar@{-}@/_/[dr] && HB \ar@{-}@/^/[dl]     \\
                                                  &  Ff\ar[d]                        &&                                     && Hf\ar[d]                                      &           \\
                                                  & FC \ar@{-}[d]   &&                                       && HC                    &                         \\
                                                  & FC \ar@{-}[d]   &&                                       &&                    &                         \\
                                                  &         \al_C \ar[d]                 &                                     &                    \\
                                                  &            GC                                   &                                     &
}
\xymatrix@R=2.2pc@C=1pc{ 
&\\            
              &  &         \\
              &  &           \\
              & = & 
}
\xymatrix@R=1.5pc@C=0.3pc{
                                                         FA \ar@{-}[d]   && FB \ar@{-}[d]        &&  HA &&  HB    \\
                                                    FA \ar@{-}@/_/[dr]   && FB \ar@{-}@/^/[dl]        &&   &&      \\
                                                  &  Ff\ar[d]                        &&                                     &&                                       &           \\
                                                  & FC \ar@{-}[d]   &&                                       &&                    &                         \\
                                                  &         \al_C \ar[d]                 &                                     &                    \\
                                                  &            GC                                   &                                     &
}
\xymatrix@R=2.2pc@C=1pc{ 
&\\            
              &  &         \\
              &  &           \\
              & = & 
}
\]
\[
\xymatrix@R=2.2pc@C=1pc{ 
&\\            
              &  &         \\
              &  &           \\
              & = & 
}
\xymatrix@R=1.5pc@C=0.3pc{
                                                         FA \ar@{-}[d]   && FB \ar@{-}[d]        &&  HA &&  HB    \\
                                                     FA \ar@{-}[d]   && FB \ar@{-}[d]        &&   &&      \\
                                                   \al_A \ar[d]   && \al_B \ar[d]        &&   &&      \\
                                                    GA \ar@{-}@/_/[dr]   && GB \ar@{-}@/^/[dl]        &&   &&      \\
                                                  &  Gf\ar[d]                        &&                                     &&                                       &           \\
                                                  & GC    &&                                       &&                    &                         
}
\xymatrix@R=2.2pc@C=1pc{ 
&\\            
              &  &         \\
              &  &           \\
              & = & 
}
\xymatrix@R=1.5pc@C=0.5pc{
                                                         FA \ar@{-}[d]   &FB\ar@{-}[dr] & HA \ar@{-}[dl]        &   HB \ar@{-}[d]   \\
                                                         FA \ar@{-}[d]   &HA & FB \ar@{-}[d]        &   HB    \\
                                                     FA \ar@{-}[d]   && FB \ar@{-}[d]        &&   &&      \\
                                                   \al_A \ar[d]   && \al_B \ar[d]        &&   &&      \\
                                                    GA \ar@{-}@/_/[dr]   && GB \ar@{-}@/^/[dl]        &&   &&      \\
                                                  &  Gf\ar[d]                        &&                                     &&                                       &           \\
                                                  & GC    &&                                       &&                    &                         
}
\]

\epf

\begin{corollary}
\label{cart7}
The functor $[1\t,-]$ gives the coreflection of $\fpMlt$ in $\fpSeq$.
Thus, the preadditive categories are exactly those of the form $\cMon(\M)$, 
for a cartesian multicategory $\M$.
\end{corollary}
\pf
By Corollary \ref{mon2}, $[1\t,-]$ gives the coreflection of $\sMlt$ in $\Seq$
and, by Proposition \ref{cart6}, if $\M$ is cartesian so it is also $[1\t,\M]$.
Furthermore, arrows $\al:M_1,\cd,M_n\to M$ in $[1\t,\M]$ are arrows 
$\al_\star:A_1,\cd,A_n\to A$ in $\M$ between the underlying objects 
(satisfying the commutativity conditions of Remark \ref{armon}) 
and $\si \al$ is computed as $\si \al_\star$.

Thus, in the correspondence
\[
\frac{\C\t \to \M}{\C\t \to [1\t,\M]}
\]
the upper functor preserves the cartesian structure iff the lower one does so.
\epf

Since $[1\t,-]$ preserves representability (by Proposition \ref{rep}),
it also gives the coreflection of $\fpRep \simeq \fpCat$ in $\fpRep\cap\fpSeq \simeq \Add$.
\begin{corollary}
The commutative monoids construction gives the coreflection of $\fpCat$ in $\Add$.
Thus, the additive categories are exactly those of the form $\cMon(\C)$, 
for a finite product category $\C$.
\epf
\end{corollary}


\subsection{The multicategory of cartesian functors}

If $\M$ and $\N$ are cartesian multicategories, we denote by $[\M,\N]\fp$
the full sub-multicategory of $[\M,\N]$ whose objects are functors in $\fpMlt$.
Since the cartesian structure and sequentiality are inherited by $[\M,\N]$ from $\N$, as well as 
by full sub-multicategories, it follows that $[\M,\N]\fp$ has a natural cartesian structure
and that it is sequential whenever $\N$ is so.
More interestingly, also representability is inherited by $[\M,\N]\fp$ from $\N$.

\begin{prop}
If $\M$ and $\N$ are cartesian multicategories 
and $\N$ is representable, then $[\M,\N]\fp$ is also representable.
\end{prop}
\pf
Recall from Proposition \ref{rep} that given functors $F,G:\M \to \N$ and a family of universal arrows
$u_A:FA,GA\to HA$, one gets a functor $H:\M\to\N$ and an arrow $u:F,G\to H$
universal in $[\M,\N]$. 
Thus we only have to check that  if $F$ and $G$ are functors in $\fpMlt$, then so it is also $H$.
Let us show for instance that $H$ preserves the contraction $f':A\to B$ of $f:A,A\to B$.
\[
\xymatrix@R=1.5pc@C=0.3pc{
                                                  & FA \ar@{-}@/_/[drr]   &&                                       &  & GA \ar@{-}@/^/[dll]  &                         \\
                                                  &                                      &&  u_A \ar[d]                 &                                     &                    \\
                                                  &                                     & &  HA  \ar@{-}[dll] \ar@{-}[drr]                          &                                     &           \\
                                                  & HA \ar@{-}@/_/[drr]   &&                                       &  & HA \ar@{-}@/^/[dll]  &                         \\
                                                  &                                      &&  Hf \ar[d]                 &                                     &                    \\
                                                  &                                     & &  HB                                   &                                     &
}
\xymatrix@R=2.2pc@C=1pc{ 
&\\            
              &  &         \\
              &  &           \\
              & = & 
}
\xymatrix@R=1.5pc@C=0.3pc{
                                                                           &FA\ar@{-}[dl]\ar@{-}[drrr]         &&&&       GA   \ar@{-}[dr]\ar@{-}[dlll]             &    \\
                                                   FA \ar@{-}@/_/[dr]   && GA \ar@{-}@/^/[dl]        && FA \ar@{-}@/_/[dr] && GA \ar@{-}@/^/[dl]     \\
                                                  &  u_A\ar[d]                        &&                                     && u_A\ar[d]                                      &           \\
                                                  & HA \ar@{-}@/_/[drr]   &&                                       && HA \ar@{-}@/^/[dll]  &                         \\
                                                  &                                      &&  Hf \ar[d]                 &                                     &                    \\
                                                  &                                     & &  HB                                   &                                     &
}
\xymatrix@R=2.2pc@C=1pc{ 
&\\            
              &  &         \\
              &  &           \\
              & = & 
}
\]

\[
\xymatrix@R=2.2pc@C=1pc{ 
&\\            
              &  &         \\
              &  &           \\
              & = & 
}
\xymatrix@R=1.5pc@C=0.3pc{
                                                                           &FA\ar@{-}[dl]\ar@{-}[drrr]         &&&&       GA   \ar@{-}[dr]\ar@{-}[dlll]             &    \\
                                                   FA \ar@{-}[d]            && GA \ar@{-}[drr]        && FA \ar@{-}[dll] && GA \ar@{-}[d]     \\
                                                   FA \ar@{-}@/_/[dr]   && FA \ar@{-}@/^/[dl]        && GA \ar@{-}@/_/[dr] && GA \ar@{-}@/^/[dl]     \\
                                                  &  Ff\ar[d]                        &&                                     && Gf\ar[d]                                      &           \\
                                                  & FB \ar@{-}@/_/[drr]   &&                                       && GB \ar@{-}@/^/[dll]  &                         \\
                                                  &                                      &&  u_B \ar[d]                 &                                     &                    \\
                                                  &                                     & &  HB                                   &                                     &
}
\xymatrix@R=2.2pc@C=0.3pc{ 
&\\            
              &  &         \\
              &  &           \\
              & = & 
}
\xymatrix@R=1.5pc@C=0.3pc{
                                                                           &FA\ar@{-}[dl]\ar@{-}[dr]         &&&&       GA   \ar@{-}[dr]\ar@{-}[dl]             &    \\
                                                   FA \ar@{-}@/_/[dr]   && FA \ar@{-}@/^/[dl]        && GA \ar@{-}@/_/[dr] && GA \ar@{-}@/^/[dl]     \\
                                                  &  Ff \ar[d]                        &&                                     && Gf \ar[d]                                      &           \\
                                                  & FB \ar@{-}@/_/[drr]   &&                                       && GB \ar@{-}@/^/[dll]  &           =              \\
                                                  &                                      &&  u_B \ar[d]                 &                                     &                    \\
                                                  &                                     & &  HB                                   &                                     &
}
\]
\[
\xymatrix@R=2.2pc@C=0.3pc{             
              &  &         \\
              &  &           \\
              & = & 
}
\xymatrix@R=1.5pc@C=0.3pc{
                                                   &FA\ar@{-}[d]        &&&&       GA   \ar@{-}[d]             &    \\
                                                  &  Ff' \ar[d]                        &&                                     && Gf' \ar[d]                                      &           \\
                                                  & FB \ar@{-}@/_/[drr]   &&                                       && GB \ar@{-}@/^/[dll]  &                         \\
                                                  &                                      &&  u_B \ar[d]                 &                                     &                    \\
                                                  &                                     & &  HB                                   &                                     &
}
\xymatrix@R=2.2pc@C=0.3pc{ 
              &  &         \\
              &  &           \\
              & = & 
}
\xymatrix@R=1.5pc@C=0.3pc{
                                                  & FA \ar@{-}@/_/[drr]   &&                                       && GA \ar@{-}@/^/[dll]  &                         \\
                                                  &                                      &&  u_A \ar[d]                 &                                     &                    \\
                                                  &                                     & &  HA   \ar@{-}[d]                                 &                                     \\
                                                  &&&                                        Hf'   \ar[d]             &    \\
                                                  &                                     & &  HB                                  &                                     
}
\]

\epf

\begin{remark}
The following well-known facts can be seen as particular cases of the above proposition
when $\M$ is also representable or when $\M$ and $\N$ are sequential, respectively:
\begin{itemize}
\item
if $\C$ and $\D$ are finite product categories,
then $[\C,\D]\fp$ has (pointwise) finite products;
\item
if $\C$ is preadditive and $\D$ is additive,
then $[\C,\D]_{add}$ has a (pointwise) additive structure.
\end{itemize}
\end{remark}


\subsection{The monoidal closed structure of $\fpMlt$}

The cartesian multicategory $[\M,\N]\fp$ is in fact the internal hom for a symmetric
monoidal closed structure on $\fpMlt$.
The tensor product is the free cartesian multicategory generated by the same arrows and relations 
defining the monoidal structure on $\sMlt$. 
The unit is $F 1^-$, the free cartesian multicategory generated by $1^-$.

This monoidal structure restricts, on $\fpSeq\simeq\preAdd$, 
to the usual monoidal structure on $\cMon$-enriched categories.
Furthermore, the role played by $1\t$ in $\sMlt$ is played here by
$F 1\t \iso {\bf N}\t$, where $\bf N$ is the rig of natural numbers
(see Remark \ref{free}).
Indeed, we have the following analogous of Corollary \ref{mon2}.

\begin{prop}
The functors ${\bf N}\t\otm\fp - $ and $[{\bf N}\t,-]\fp$ give respectively a reflection and a coreflection
of $\fpMlt$ in $\fpSeq$.
\end{prop}
\pf
The coreflection formula follows from Corollary \ref{cart7}, 
since 
\[
[{\bf N}\t,\M]\fp \iso [F 1\t,\M]\fp \iso [1\t,\M]
\]
As for the reflection, ${\bf N}\t\otm\fp \M$ is sequential (since it has a central monoid) 
and, for any sequential $\L$, we have natural isomorphisms
\[    
\begin{array}{c}
{\bf N}\t\otm\fp\M \to \L \\ \hline
\M \to [{\bf N}\t,\L]\fp \\ \hline
\M \to \L
\end{array}
\]

\epf

Similar results concerning the reflection of Lawvere theories in annular theories 
are stated in \cite{freyd}, \cite{lawvere} and \cite{wraith}.
Note that, on the other hand, the coreflection $[{\bf N}\t,-]\fp$ can not work
in the single sorted context, since an operad may of course support several commutative monoids.

\begin{refs}

\bibitem[Adamek et al., 2011]{adamek} J. Adamek, J. Rosicky, E.M. Vitale (2011), {\em Algebraic theories: 
a categorical introduction to general algebra}, Cambridge tracts in mathematics, {\bf 184}.

\bibitem[Benabou, 2000]{benabou} J. Benabou (2000), {\em Distributors at work}, Lecture notes by T. Streicher of a course 
given at TU Darmstadt, 2000.

\bibitem[Boardman \& Vogt, 1973]{boardman} J. M. Boardman, R. M. Vogt (1973), {\em Homotopy invariant algebraic structures 
on topological spaces}, Lecture Notes in Mathematics, {\bf 347}, Springer-Verlag.

\bibitem[Day et al., 2005]{day} B. Day, E. Panchadcharam, R. Street (2005), {\em On centres and lax centres for promonoidal categories},
Colloque International  Charles Ehresmann : 100 ans Universite de Picardie Jules Verne, Amiens.


\bibitem[T. Fiore, 2005]{fiore2} T. Fiore (2006), {\em Pseudo limits, biadjoints, and pseudo algebras : 
Categorical foundations of conformal field theory}, Memoirs of the American Mathematical Society, 182(860), math.CT/04028298.

\bibitem[Fox, 1976]{fox} T. Fox (1976) Coalgebras and cartesian categories, {\em Communications in Algebra}, {\bf 4(7)}, 665-667.

\bibitem[Freyd, 1966]{freyd} P. Freyd (1966), Algebra valued functors in general and tensor products in particular, 
{\em Colloquium mathematicum} {\bf 14}, 89-106.

\bibitem[Gould, 2008]{gould} M. Gould (2008), {\em Coherence for categorified operadic theories}, PhD Thesis, math.CT/1002.0879.

\bibitem[Hermida, 2000]{hermida} C. Hermida (2000), Representable multicategories, {\em Advances in Math.}, {\bf 151}, 164-225.

\bibitem[Hermida, 2004]{hermida2} C. Hermida (2004), {\em Fibrations for abstract multicategories}, in: 
Galois Theory, Hopf Algebras and Semiabelian Categories, Fields inst. comm. AMS, 281-293.




\bibitem[Lambek, 1989]{lambek} J. Lambek (1989), {\em Multicategories revisited}, 
in: Contemporary Mathematics 92 (Amer. Math. Soc., Providence), 217-239.

\bibitem[Lawvere, 2004]{lawvere} F.W. Lawvere (2004), Functorial semantics of algebraic theories and some algebraic problems 
in the context of functorial semantics of algebraic theories, {\em Reprints in Theory and Appl. of Cat.} {\bf 5}, 1-121.

\bibitem[Leinster, 2003]{leinster} T. Leinster (2003), {\em Higher operads, higher categories}, Cambridge University Press, math.CT/0305049.

\bibitem[MacLane, 1971]{maclane} S. MacLane (1971), {\em Categories
for the working mathematician}, Springer-Verlag.

\bibitem[Moerdijk \& Weiss, 2007]{moerdijk} I. Moerdijk, I. Weiss (2007), Dendroidal sets, {\em Algebr. Geom. Topol.} {\bf 7}, 1441-1470.


\bibitem[Pisani, 2013]{pisani} C. Pisani (2013), {\em Some remarks on multicategories and additive categories}, math.CT/1304.3033.


\bibitem[Tronin, 2011]{tronin} S.N. Tronin (2011), Natural multitransformations of multifunctors, {\em Russian Mathematics},
{\bf 55(11)}, 49-60.

\bibitem[Trova, 2010]{trova} F. Trova (2010), {\em On the Geometric Realization of Dendroidal Sets},
Master Thesis, Leiden University, University of Padova.

\bibitem[Weiss, 2011]{weiss} I. Weiss (2011), {\em From operads to dendroidal sets}, in: Mathematical Foundations of Quantum Field 
and Perturbative String Theory, Proc. of Symp. in Pure Math. AMS, {\bf 83}, 21-70. 

\bibitem[Wraith, 1970]{wraith} G. Wraith (1970), Algebraic theories, Lecture Notes Series {\bf 22},
Matematisk Institut, Aarhus Universitet.
\end{refs}

\end{document}